\documentclass[twocolumn]{autart}    

\usepackage{graphicx}          
\usepackage{amsfonts,amssymb}
\usepackage{amsmath}
\usepackage{bm}
\usepackage{array}
\usepackage{multirow}
\usepackage{makecell}    

\usepackage[linesnumbered,ruled,vlined]{algorithm2e}

\usepackage[caption=false]{subfig}

\usepackage{stfloats}   

\allowdisplaybreaks[4]

\usepackage{booktabs}   
\usepackage{color}
\usepackage{diagbox}

\usepackage{natbib} 

\definecolor{rev}{RGB}{191,0,64} 

\begin{document}

\begin{frontmatter}

\title{Collaborative-Online-Learning-Enabled Distributionally Robust Motion Control for Multi-Robot Systems\thanksref{footnoteinfo}}               

\thanks[footnoteinfo]{This paper was not presented at any IFAC meeting. (Corresponding authors: C. Ning and Y. Shi)}

\author[a,b]{Chao Ning}\ead{chao.ning@sjtu.edu.cn},    
\author[a,b]{Han Wang}\ead{h.wang@sjtu.edu.cn},               
\author[a,b]{Longyan Li}\ead{lilongyan@sjtu.edu.cn},  
\author[d]{Yang Shi}\ead{yshi@uvic.ca}  

\address[a]{School of Automation and Intelligent Sensing, Shanghai Jiao Tong University, Shanghai, China}  
\address[b]{Key Laboratory of System Control and Information Processing, Ministry of Education of China, Shanghai, China}             
\address[d]{Department of Mechanical Engineering, University of Victoria, Victoria, Canada}        


\begin{abstract}                          
This paper develops a novel COllaborative-Online-Learning (COOL)-enabled motion control framework for multi-robot systems to avoid collision amid randomly moving obstacles whose motion distributions are partially observable through decentralized data streams. To address the notable challenge of data acquisition due to occlusion, a COOL approach based on the Dirichlet process mixture model is proposed to efficiently extract motion distribution information by exchanging among robots selected learning structures. By leveraging the fine-grained local-moment information learned through COOL, a data-stream-driven ambiguity set for obstacle motion is constructed. We then introduce a novel ambiguity set propagation method, which theoretically admits the derivation of the ambiguity sets for obstacle positions over the entire prediction horizon by utilizing obstacle current positions and the ambiguity set for obstacle motion. Additionally, we develop a compression scheme with its safety guarantee to automatically adjust the complexity and granularity of the ambiguity set by aggregating basic ambiguity sets that are close in a measure space, thereby striking an attractive trade-off between control performance and computation time. Then the probabilistic collision-free trajectories are generated through distributionally robust optimization problems. The distributionally robust obstacle avoidance constraints based on the compressed ambiguity set are equivalently reformulated by deriving separating hyperplanes through tractable semi-definite programming. Finally, we establish the probabilistic collision avoidance guarantee and the long-term tracking performance guarantee for the proposed framework. The numerical simulations are used to demonstrate the efficacy and superiority of the proposed approach compared with state-of-the-art methods.
\end{abstract}

\end{frontmatter}

\section{Introduction}  \label{section1}
Exhibiting superior task execution efficiency, robustness, flexibility, and fault tolerance compared with single-robot systems, multi-robot systems are extensively applied across various domains including military, agriculture, and transportation \cite{arai2002advances,rizk2019cooperative}. Safety is one of the most fundamental challenges for the motion control of robotic systems \cite{wang2017safety}.  In uncertain and dynamic environments the unexpected motion of obstacles undeniably poses a risk to robot collision-free navigation and brings disastrous consequences. A crucial aspect of collision avoidance involves acquiring accurate motion probability distributions of obstacles. However, achieving perfect accuracy in obtaining the distributions is typically impractical, necessitating estimation from available data \cite{A20}. Therefore, the key to obstacle avoidance for multi-robot systems lies in the efficient utilization of decentralized data from different robots. The objective of this work is to develop a collaborative-online-learning-enabled distributionally robust motion control framework for multi-robot systems. This framework aims to learn motion distributions in a collaborative manner, to exhibit robustness against distribution estimation errors, and to be suitable for real-time requirements. 

A plethora of methods have been proposed for collision avoidance motion control in robotic systems. We categorize and discuss these methods separately based on their application environments and target systems. Based on application environments, the methods can be classified into those designed for deterministic environments and those for uncertain environments. Collision avoidance in deterministic environments is a widely studied problem, with numerous well-established algorithms having been proposed, such as the penalty methods \cite{A1}, control barrier function methods \cite{dai2023safe}, virtual holonomic constraints methods \cite{A19}, and among many others. Although the assumption of a deterministic environment simplifies the collision avoidance problem, it can incur safety issues in reality because prediction errors for obstacle positions are inevitable. Optimistically assuming perfect knowledge of future obstacle states may lead to unsafe trajectories when actual obstacle motions deviate from their predicted paths. Uncertain environments are more aligned with real-world scenarios, and our research is focused on this setting. In recent years, a large body of work on safe learning-based control has also been developed to address environmental uncertainty. For instance, \cite{fisac2018general} proposed a general framework that quantifies uncertainty from learning and ensures safety by augmenting the control space with safety margins. \cite{yuan2024lightweight} developed a lightweight distributed Gaussian Process (GP) regression method for online multi-robot learning with limited bandwidth and computational power. \cite{yuan2025distributed} proposed a distributed safe learning and planning framework for multi-robot systems, which leverages GP to model external disturbance dynamics and derive the sufficient conditions to guarantee the safety. \cite{9683423} proposed an online GP-based formation control law that learns unknown dynamics during runtime and guarantees bounded formation errors with high probability. \cite{7798979} developed a GP-based safe learning framework to estimate regions of attraction in uncertain nonlinear systems, enabling safe exploration with high-probability stability guarantees. For collision avoidance motion control in uncertain environments, chance-constrained methods are the most popular, ascribing to their capability to directly curb the chance of collision \cite{hakobyan2019risk, sharma2020risk}. However, the safety probability guarantee of all the aforementioned chance-constrained methods heavily relies on the precise knowledge of distribution that is rarely available. To enhance the robustness against distribution errors, a novel control paradigm based on Distributionally Robust Optimization (DRO), called Distributionally Robust Model Predictive Control (DRMPC) was firstly proposed in \cite{li2021distributionally} and \cite{li2023distributionally}. Furthermore, a learning based data-driven MPC by using Conditional Value-at-Risk (CVaR) was proposed in \cite{ning2021online}. Building upon DRMPC, a novel motion control framework, referred to as Distributionally Robust Motion Control (DRMC) has attracted widespread attention, where the motion control is accomplished by solving a DRO problem at each control time. For the DRO problem, chance constraints should hold for all distributions within a set, called the ambiguity set \cite{van2015distributionally}. \cite{hakobyan2021wasserstein} utilized CVaR to measure safety risk and constructed a distance-based ambiguity set comprising all distributions within a Wasserstein ball centered at the empirical distribution. Building upon the work above, GP regression was employed to learn the motion distribution from real-time data and the estimated unimodal Gaussian distribution serves as the center of the Wasserstein ball \cite{hakobyan2022distributionally}. Note that all of the methods above require solving a non-convex optimization problem at each control time, which is computationally intensive and unsuitable for real-time operation. Fortunately, the latest study \cite{safaoui2023distributionally} integrated the separating hyperplane technique with DRMC to convexify collision avoidance constraints, which enhances computational efficiency at the expense of control performance. The above literature review indicates that existing DRMC methods either entail substantial computational burdens or sacrifice a significant amount of control performance to enhance computational efficiency. DRMC that ensures safety while flexibly balancing computation time and control performance is currently lacking in academia. Moreover, all the aforementioned DRMC methods inherently impose restrictive assumptions, namely that the nominal trajectories of obstacles are the priori knowledge which is challenging to satisfy in practice.

Based on the target systems, the methods can be categorized into those designed for single-robot systems and those for multi-robot systems. The majority of methods for collision avoidance motion control in uncertain environments focus on single-robot systems. Among them, some relied solely on historical data to estimate distributions and construct ambiguity sets without utilizing real-time data \cite{hakobyan2021wasserstein, safaoui2023distributionally}, while others assumed that data streams are readily accessible \cite{hakobyan2022distributionally}. However, in practice, the occlusions between obstacles and robots make online data acquisition of single-robot systems challenging. For multi-robot systems, the simplest method is to adopt a central coordinator to decide every robot’s motion \cite{honig2018trajectory, soria2021predictive}. However, the centralized method tends to be intractable and unrealistic for large-scale multi-robot systems due to the unbearable computation time of the central estimator. Fortunately, the literature presents a variety of methods for motion control in decentralized multi-robot systems, including sequential methods \cite{richards2004decentralized}, synchronous methods \cite{luis2019trajectory}, and asynchronous methods \cite{csenbacslar2022asynchronous}. However, the majority of collision avoidance motion control methods for multi-robot systems are only applicable to deterministic environments. A research gap currently exists in motion control for multi-robot systems in uncertain environments, particularly in the area of DRMC.

Motivated by these considerations, this paper proposes a novel COllaborative-Online-Learning-enabled DRMC (COOL-DRMC) approach for decentralized multi-robot systems to guide all the robots in tracking their reference trajectories while containing the chance of collisions. We consider a realistic setting where (i) the obstacle motion distribution can be partially inferred from data streams, (ii) the occlusion of obstacles or other robots impedes online data acquisition of obstacle positions, and (iii) the nominal positions of obstacles in the future step are unavailable. At each control time, a novel collaborative-online-learning method based on the Dirichlet Process Mixture Model (DPMM) \cite{gomes2008incremental} is leveraged to decipher the structural properties of the motion distribution. By leveraging the local-moment information of the motion distributions, an ambiguity set for obstacle motion is constructed. Without making the restrictive assumption that nominal trajectories of obstacles are known \cite{safaoui2023distributionally}, we put forward a novel ambiguity set propagation method, enabling the theoretical derivation of the ambiguity sets for obstacle positions throughout the entire prediction horizon. Regarding the ambiguity set, we develop a provably safe compression scheme to adjust its complexity and granularity. By controlling the degree of compression, the proposed scheme enjoys an attractive trade-off between control performance and computation time. The distributionally robust robot-obstacle collision avoidance constraints are reformulated as affine constraints through separating hyperplanes, which can be computed through tractable Semi-Definite Programming (SDP). Furthermore, the inter-robot separating hyperplanes are generated by a safe spatial allocation protocol. Finally, we theoretically establish the probabilistic collision avoidance guarantee and the long-term tracking performance guarantee for COOL-DRMC. The main contributions of this paper are given below.

Many studies based on DRO focus on addressing the uncertainties of systems. For example, \cite{li2023distributionally} proposed a DRMPC method with output feedback to address uncertainties with unbounded support sets in systems. In contrast, our method focuses on environmental uncertainties and emphasizes online updating of the ambiguity set through online learning. Additionally, \cite{kandel2023safe} proposed an end-to-end framework for safe learning-based control, which online learns the system model and constructs an ambiguity set for the learning error. Specifically, this paper constructs a Wasserstein-based ambiguity set using real-time errors as sample points, whereas our method employs a non-parametric Bayesian method for decentralized real-time data and constructs a local-moment-based ambiguity set.

The remainder of this paper is structured below. The system modeling and problem statement are presented in Section \ref{section2}. Section \ref{section3} introduces the COllaborative Online Learning (COOL) method. Section \ref{section4} describes the separating hyperplane computation and control algorithm. Numerical simulations are presented in Section \ref{section5}. Conclusions are given in Section \ref{section6}.

\textit{Notation}: $\mathbb{N}_+$ and $\mathbb{N}_{++}$ denote the sets of non-negative integers and positive integers, respectively. For a positive integer $N$, $\mathbb{N}_+^N$ and $\mathbb{N}_{++}^N$ denote the sets $\{ 0,1,...,N \}$ and $\{ 1,2,...,N \}$, respectively. For a real number $a$, $(a)^+ = \max(a,0)$. $0_d$ is a vector of zeros in $\mathbb{R}^d$. For matrices $A$ and $B$, $A\bullet B$ represents their Frobenius inner product. The sum of a vector $a$ and a set $\mathcal{A}$ is defined by $a \oplus \mathcal{A} := \{y+a |a \in \mathcal{A}\}$. For a convex compact set $\mathcal{C}$ and a unit vector $h$, the support function and the support point of $\mathcal{C}$ in the direction $h$ are $\mathcal{S}_{\mathcal{C}}(h) = \max_{y \in \mathcal{C}}y^Th$ and $\mathcal{P}_{\mathcal{C}}(h) = \arg\max_{y \in \mathcal{C}}y^Th$, respectively.

\section{System modeling and problem statement}  \label{section2}
Consider a group of $N$ robots navigating a cluttered environment with $L$ static or random moving obstacles. The obstacles can prevent the robots from obtaining environmental information by blocking the observation of robots. The objective is to control them in a collaborative way to follow their reference trajectories $y^r_i$, $i \in \mathbb{N}_{++}^{N}$ obtained from the upper-level planner \cite{karaman2011sampling}, while curbing collision probabilities below a prescribed tolerance.

\subsection{Modeling of robots and obstacles}     \label{section2.1}
In this work, we consider robot $i \in \mathbb{N}_{++}^{N}$ with convex and compact rigid bodies, denoted by $\mathcal{R}_i \subset \mathbb{R}^d$, where $0_d \in \mathcal{R}_i$ and $d \in \{2,3\}$. The dynamic of robot $i$ is modeled by the following time-invariant discrete-time system.
\begin{equation} \label{e2.1.1}
\begin{split}
&x_i(t+1)=f(x_i(t),u_i(t))\\
&y_i(t)=h(x_i(t))
\end{split}
\end{equation}
where $x_i(t) \in \mathcal{X} \subseteq \mathbb{R}^{n_x}$, $u_i(t) \in \mathcal{U} \subseteq \mathbb{R}^{n_u}$ and $y_i(t) \in \mathbb{R}^d$ are the state, input and output, respectively. $f: \mathbb{R}^{n_x} \times \mathbb{R}^{n_u} \rightarrow \mathbb{R}^{n_x}$ and $h: \mathbb{R}^{n_x} \times \mathbb{R}^{n_u} \rightarrow \mathbb{R}^{d}$ represent the system dynamics and the output mapping, respectively. The system output $y_i(t)$ is defined as the Cartesian coordinates of the robot Center of Mass (CoM). The admissible state set $\mathcal{X}$ and input set $\mathcal{U}$ are convex and compact polytopes. Let $k|t$ denote the $k$-step ahead of time $t$. The prediction horizon is taken as $K$, and then $k \in \mathbb{N}_{+}^{K-1}$. The region expected to be occupied by robot $i$ at time $t+k$ is denoted by $y_i(k|t) \oplus \mathcal{R}_i$.

For the modeling of obstacles, we consider obstacle $\ell \in \mathbb{N}_{++}^{L}$ with convex and compact rigid bodies, denoted by $\mathcal{O}_{\ell} \subset \mathbb{R}^d$, where $0_d \in \mathcal{O}_{\ell}$. The motion of the obstacle between two time steps is modeled by translation in (\ref{e2.1.3}).
\begin{equation} \label{e2.1.3}
\begin{split}
\varphi_{\ell}(t+1)=\varphi_{\ell}(t)+\omega_{\ell,t}
\end{split}
\end{equation}
where $\varphi_{\ell}(t) \in \mathbb{R}^d$ denotes the CoM of the obstacle $\ell$ at time $t$. $\omega_{\ell,t} \in \mathbb{W}_{\omega}$ is a random translation vector, whose distribution is unknown but can be partially inferred from data. $\mathbb{W}_{\omega}$ is the convex support set. Given the observation of the position of obstacle $\ell$ at time $t$, the predicted positions over the prediction horizon are given by
\begin{equation} \label{e2.1.4}
\begin{split}
\varphi_{\ell}(k|t)&=\varphi_{\ell}(t)+\omega_{\ell,t}^{(1)}+...+\omega_{\ell,t}^{(k)}\\
&=\varphi_{\ell}(t)+\omega_{\ell,t,k}
\end{split}
\end{equation}
where $\omega_{\ell,t}^{(1)},...,\omega_{\ell,t}^{(k)}$ are the random motions between two adjacent time steps, which are independent and identically distributed (i.i.d.) and follow the same distribution conditioned on the information available at current time $t$. $\varphi_{\ell}(k|t)$ and $\omega_{\ell,t,k}$ are the predicted position at time $t+k$ and the random motion between $k$ time steps, which are all random variables.
\begin{rem} \label{rem2.1.0}
The motion of the obstacle can be modeled as the linear dynamics $\varphi_{\ell}(t+1) = A_{\ell}\varphi_{\ell}(t) + B_{\ell}\omega_{\ell,t}$. In this way, the predicted positions over the prediction horizon are given by $\varphi_{\ell}(k|t)= A_{\ell}^{k} \varphi_{\ell}(t) + A_{\ell}^{k-1}B_{\ell} \omega_{\ell,t}^{(1)} +...+ B_{\ell}\omega_{\ell,t}^{(k)}$. Our method can still handle such obstacle dynamics but results in highly complex expressions. For ease of exposition, we only consider translational uncertain motion in this paper.
\end{rem}
Obstacles in practical scenarios are typically complex and non-convex. However, obstacles can be decomposed into multiple simple obstacles and over-approximated as their convex hulls. Additionally, obstacles are not always observable for robots due to the occlusion of other obstacles or robots, especially when there are large static obstacles in the environment. Therefore, we define observable obstacles for robot $i$ at time $t$ as follows, whose indices are denoted as $\mathbb{O}_{i,t} \subseteq \mathbb{N}_{++}^{L}$.
\begin{defn}    \label{defn2.1.3}
(Observable obstacle for robot $i$ at time $t$) The obstacle $\ell$ is called the observable obstacle for robot $i$ at time $t$ if and only if (iff) there is no other obstacle or robot in the way between part of the obstacle $\ell$ and the robot $i$. Mathematically, there exists $\alpha \in \mathcal{O}_{\ell,t}$ such that $\{\theta \alpha + (1-\theta)y_i(t) | \theta \in (0,1)\} \cap \{[\bigcup_{j=1,j\neq N}^{I}(y(t) \oplus \mathcal{R}_{j})] \cup [\bigcup_{s=1,s\neq \ell}^{L}(\varphi_{s}(t) \oplus \mathcal{O}_{s})]\} = \emptyset$.
\end{defn}

\subsection{Collision avoidance constraints}    \label{section2.3}
We aim to design a controller that drives robot $i$ to follow its reference trajectory while ensuring system safety. Specifically, the safety of a multi-robot system means limiting inter-robot and robot-obstacle collisions below the acceptable threshold. Therefore, we formalize the inter-robot collision avoidance and robot-obstacle collision avoidance constraints in this subsection.
\begin{assum}   \label{assum2.3.1}
We assume that a robot can communicate without delay with other robots.
\end{assum}
For the inter-robot collision avoidance constraints, we require that all robots avoid collision among themselves. At time $t$, robots $i$ and $j$ will not collide at $k$-step ahead of $t$, iff the following constraint is satisfied.
\begin{equation} \label{e2.3.1.1}
\begin{split}
(y_i(k|t) \oplus \mathcal{R}_{i}) \cap (y_j(k|t) \oplus \mathcal{R}_{j}) = \emptyset
\end{split}
\end{equation}
Note that for the two convex and compact sets $y_i(k|t) \oplus \mathcal{R}_i$ and $y_j(k|t) \oplus \mathcal{R}_j$, if they do not intersect, there always exists a separating hyperplane $\mathcal{H}_{i,j,t,k}^{rob}:=\left\{y \in \mathbb{R}^d \big| y^T h_{i,j,t,k}^{rob} + g_{i,j,t,k}^{rob} = 0\right\}$ to separate them. By employing $\mathcal{H}_{i,j,t,k}^{rob}$ and the support function \cite{boyd2004convex}, constraint (\ref{e2.3.1.1}) can be convexified into the following two affine constraints.
\begin{align}
y_i(k|t)^T h_{i,j,t,k}^{rob} + \mathcal{S}_{\mathcal{R}_i}(h_{i,j,t,k}^{rob}) + g_{i,j,t,k}^{rob} + \tfrac{r}{2} \leq 0   \label{e2.3.1.5}\\
\hspace{-5pt} -y_j(k|t)^T h_{i,j,t,k}^{rob} + \mathcal{S}_{\mathcal{R}_j}(-h_{i,j,t,k}^{rob}) - g_{i,j,t,k}^{rob} + \tfrac{r}{2} \leq 0    \label{e2.3.1.6}
\end{align}
where $r>0$ is the safe margin. If constraints (\ref{e2.3.1.5}) and (\ref{e2.3.1.6}) are simultaneously satisfied, constraint (\ref{e2.3.1.1}) is satisfied.

For the robot-obstacle collision avoidance constraints, we require that the robots avoid the obstacles in the environment. Similarly, if the following constraint is satisfied, robot $i$ can avoid collision with obstacle $\ell$ at $k$-step ahead of $t$.
\begin{equation} \label{e2.3.2.1}
\begin{split}
(y_i(k|t) \oplus \mathcal{R}_{i}) \cap (\varphi_{\ell}(k|t) \oplus \mathcal{O}_{\ell}) = \emptyset
\end{split}
\end{equation}
Note that $\varphi_{\ell}(k|t) \oplus \mathcal{O}_{\ell}$ is a convex and compact set, and thus constraint (\ref{e2.3.2.1}) can also be reformulated as the following deterministic affine constraints using the separating hyperplane $\mathcal{H}_{i,\ell,t,k}^{obs}:=\left\{y \in \mathbb{R}^d \big| y^T h_{i,\ell,t,k}^{obs} + g_{i,\ell,t,k}^{obs} = 0\right\}$ and the support function.
\begin{align}
y_i(k|t)^T h_{i,\ell,t,k}^{obs} + \mathcal{S}_{\mathcal{R}_i}(h_{i,\ell,t,k}^{obs}) + g_{i,\ell,t,k}^{obs} + r \leq 0   \label{e2.3.2.2}\\
-\varphi_{\ell}(k|t)^T h_{i,\ell,t,k}^{obs} + \mathcal{S}_{\mathcal{O}_{\ell}}(-h_{i,\ell,t,k}^{obs}) - g_{i,\ell,t,k}^{obs} \leq 0    \label{e2.3.2.3}
\end{align}
If constraints (\ref{e2.3.2.2}) and (\ref{e2.3.2.3}) are satisfied concurrently, constraint (\ref{e2.3.2.1}) is also satisfied. We note that if the position of the obstacle $\ell$ at the predicted time $k|t$ is known, $h_{i,\ell,t,k}^{obs}$ can be directly computed as follows.
\begin{equation}   \label{e_new1}
\begin{split}
    h_{i,\ell,t,k}^{obs} = \frac{\varphi_{\ell}(k|t) - y_i(k|t)}{\Vert \varphi_\ell(k|t) - y_i(k|t) \Vert_2}
\end{split}
\end{equation}
However, due to the uncertainty of $\varphi_{\ell}(k|t)$, $h_{i,\ell,t,k}^{obs}$ cannot be directly computed using (\ref{e_new1}). Instead, $\varphi_{\ell}(k|t)$ or its distribution needs to be learned or estimated from data. Therefore, we denote the left-hand side of (\ref{e2.3.2.3}) as $\mathcal{L}_{\mathcal{H}_{i,\ell,t,k}^{obs}}(\varphi_{\ell}(k|t)\oplus\mathcal{O}_{\ell})$ and consider the following Distributionally Robust CVaR (DR-CVaR) constraint.
\begin{equation} \label{e2.3.2.5}
\begin{split}
\mathop{\sup}\limits_{\upsilon \in \widetilde{\mathbb{D}}_{\ell,t,k}} {\rm CVaR}_{\alpha^{u}}^{\upsilon}(\mathcal{L}_{\mathcal{H}_{i,\ell,t,k}^{obs}}(\varphi_{\ell}(k|t)\oplus\mathcal{O}_{\ell})) \leq 0
\end{split}
\end{equation}
where $\alpha^{u}$ is the used-defined confidence level and $\widetilde{\mathbb{D}}_{\ell,t,k}$ is the ambiguity set for $\varphi_{\ell}(k|t)$. Additionally, ${\rm CVaR}_{\alpha}^{\upsilon}(X)=\min_{z\in\mathbb{R}}\mathbb{E}^{\upsilon}\left[ z+(X-z)^+/(1-\alpha) \right]$ measures the conditional expectation of the loss within the $(1-\alpha)$ worst-case quantile. If $\widetilde{\mathbb{D}}_{\ell,t,k}$ encapsulates the underlying distribution of $\varphi_{\ell}(k|t)$ and constraints (\ref{e2.3.2.2}) and (\ref{e2.3.2.5}) are simultaneously satisfied, constraint (\ref{e2.3.2.1}) is satisfied with a probability of at least $\alpha^{u}$.

It is evident that when constraints (\ref{e2.3.2.2}) and (\ref{e2.3.2.5}) are simultaneously satisfied, the distributionally robust collision avoidance requirement between agent $i$ and obstacle $\ell$ is guaranteed. However, since obstacles are not controllable, constraint (\ref{e2.3.2.5}) cannot be directly enforced during the motion control problem. Instead, constraint (\ref{e2.3.2.5}) is taken into account in the computation of the separating hyperplane $\mathcal{H}_{i,\ell,t,k}^{obs}$, while constraint (\ref{e2.3.2.2}) is explicitly enforced in the motion control problem of agent $i$. Through this decoupled treatment, the distributionally robust agent--obstacle collision avoidance requirement is effectively ensured.

\subsection{Distributionally robust optimization for motion control}
Given the separating hyperplanes $\mathcal{H}_{i,j,k,t}^{rob}$ and $\mathcal{H}_{i,\ell,k,t}^{obs}$, the optimal control strategy is designed through the following optimization problem.
\begin{subequations}        \label{e2.4.2.1}
\begin{alignat}{2}
&\hspace{-15pt} \min_{\bm{u_i(t)}} J_{i}(\bm{u_{i}(t)}) = \sum_{k=0}^{K-1} l_{i}(x_{i}(k|t), u_{i}(k|t)) \notag\\
        &\hspace{70pt}+ q_{i}(x_{i}(K|t))    \label{e2.4.2.1a}\\
s.t.&\ x_i(k+1|t)=f(x_i(k|t),u_i(k|t))     \label{e2.4.2.1b}\\
    &y_i(k|t)=h(x_i(k|t))     \label{e2.4.2.1c}\\
    &x_i(0|t)=x_i(t)        \label{e2.4.2.1d}\\
    &x_i(k|t) \in \mathcal{X}    \label{e2.4.2.1e}\\
    &u_i(k|t) \in \mathcal{U}    \label{e2.4.2.1f}\\
    &\hspace{-7pt} y_i(k|t)^T h_{i,j,t,k}^{rob}+\mathcal{S}_{\mathcal{R}_i}(h_{i,j,t,k}^{rob})+g_{i,j,t,k}^{rob}+\tfrac{r}{2} \leq 0    \label{e2.4.2.1g}\\
    &\hspace{-7pt} y_i(k|t)^T h_{i,\ell,t,k}^{obs} + \mathcal{S}_{\mathcal{R}_i}(h_{i,\ell,t,k}^{obs}) + g_{i,\ell,t,k}^{obs} + r \leq 0    \label{e2.4.2.1h}\\
    &x_{i}(K|t) \in \mathcal{X}_{i}^{f}    \label{e2.4.2.1i}
\end{alignat}
\end{subequations}
where $\bm{u_i(t)}=(u_i(0|t),...,u_i(K-1|t))$. Constraints (\ref{e2.4.2.1b}) and (\ref{e2.4.2.1f}) should hold for $k\in\mathbb{N}_{+}^{K-1}$. Constraints (\ref{e2.4.2.1c}) and (\ref{e2.4.2.1e}) should be imposed for $k\in\mathbb{N}_{+}^{K}$. Constraint (\ref{e2.4.2.1g}) should be satisfied for $k\in\mathbb{N}_{++}^{K-1}$ and ${\forall} j \neq i$. Constraint (\ref{e2.4.2.1h}) should be satisfied for $k\in\mathbb{N}_{++}^{K-1}$ and $\ell \in \mathbb{N}_{++}^{L}$. Here, a typical choice for the stage cost function $l_{i}: \mathbb{R}^{n_x} \times \mathbb{R}^{n_u} \times \mathbb{R} \rightarrow \mathbb{R}$ and the terminal cost function $q_{i}: \mathbb{R}^{n_x} \times \mathbb{R} \rightarrow \mathbb{R}$ is as follows.
\begin{equation}    \label{e2.4.2.2}
\begin{split}
    &l_{i}(x_{i}(k|t),u_{i}(k|t)) = \\
    &\Vert x_{i}(k|t)-r_{i}^{x}(t+k) \Vert_{Q} + \Vert u_{i}(k|t)-r_{i}^{u}(t+k) \Vert_{R}  \\
    &q_{i}(x_{i}(K|t)) = \Vert x_{i}(K|t)-r_{i}^{x}(t+K) \Vert_{P}
\end{split}
\end{equation}
where $(r_{i}^{x},r_{i}^{u})$ is the reference trajectory generated from the upper-level planner, and $Q,R,P$ are positive semi-definite matrices. Constraints (\ref{e2.4.2.1b}) and (\ref{e2.4.2.1c}) calculate the predicted state and position, and (\ref{e2.4.2.1e}) and (\ref{e2.4.2.1f}) are the constraints of system state and input. Note that by respectively incorporating constraints (\ref{e2.3.1.5}) and (\ref{e2.3.1.6}) into the optimization problems of robot $i$ (as constraint (\ref{e2.4.2.1g})) and robot $j$, the inter-robot collision avoidance constraint with respect to (w.r.t) $i$ and $j$ is satisfied. Additionally, the DR-CVaR constraint (\ref{e2.3.2.5}) is taken into account during the computation of $\mathcal{H}_{i,\ell,t,k}^{obs}$. Consequently, the probabilistic obstacle avoidance constraint can also be satisfied by incorporating constraint (\ref{e2.3.2.2}) into (\ref{e2.4.2.1h}) as constraint (\ref{e2.4.2.1h}). Finally, constraint (\ref{e2.4.2.1i}) is a terminal constraint, which is crucial for establishing the long-term tracking performance guarantee.

Although the collision avoidance constraints are tractable affine constraints, several challenges arise in computing the separating hyperplane $\mathcal{H}_{i,j,k,t}^{rob}$ and $\mathcal{H}_{i,\ell,k,t}^{obs}$. The first challenge is to leverage the advantages of the multi-robot system to construct the ambiguity sets throughout the prediction horizon efficiently by integrating the observed data streams from all robots in the decentralized system which can accurately capture the underlying distribution information. Moreover, it is also a challenge for the robots in a decentralized system to independently compute the separating hyperplane $\mathcal{H}_{i,j,t,k}^{rob}$ while ensuring safety and rationality. Additionally, $\mathcal{H}_{i,\ell,t,k}^{obs}$ is crucial for both control performance and safety but nontrivial to compute due to the DR-CVaR constraint (\ref{e2.3.2.5}), which involves an infinite number of distributions in the ambiguity set.

\section{Preliminaries} \label{section03}
In this section, we provide the essential background necessary to support the proposed framework.
\subsection{Background on DPMM}   \label{section3.1}
In this paper, the DPMM is used to model the stochastic motion of obstacles by learning the probability distribution of their inter-sample displacements. In this subsection, our introduction to the DPMM is organized into two parts. The first part provides a detailed explanation of the classical DPMM theory. The second part presents the online variational inference algorithm that enables practical implementation in an online setting, along with its mathematical derivation.
\subsubsection{DPMM theory}  \label{section3.1.1}
As one of the most popular nonparametric Bayesian models, the DPMM can automatically decipher the structural property of the distribution. Specifically, the Dirichlet process is a distribution over distributions, and thus the draw from a Dirichlet process on the parameter space produces a discrete probability measure on the same parameter space. By considering the mixture model as this discrete probability measure, the DPMM is capable of deciphering the multimodal data distribution without prior information on the number of mixture components. To elaborate in detail, the sampling procedure of the DPMM is presented as follows.
\begin{equation}    \label{e3.1.a1}
\begin{split}
    &G \sim \mathrm{DP} (\alpha,H); \quad \phi \sim G; \quad \omega \sim F(\phi)
\end{split}
\end{equation}
where $\mathrm{DP}$ represents the Dirichlet process (DP) with the concentration parameter $\alpha$ and the base distribution $H$ on the space of parameter $\phi$. $G$ is a discrete probability measure on the same space. $F(\phi)$ is the parametric distribution family with $\phi$ being its parameters and the data point $\omega$ is sampled from the distribution $F(\phi)$. Here, the data point $\omega$ corresponds to the observed displacement of an obstacle between two consecutive time steps, capturing its random motion over one sampling interval. Drawing $G$ from $\mathrm{DP}(\alpha,H)$ admits the following stick-breaking procedure \cite{sethuraman1994constructive}.
\begin{equation}    \label{e3.1.a2}
\begin{split}
&v_{n} \sim \mathrm{Beta}(1,\alpha)  \qquad \qquad \phi_n \sim H \\
&\gamma_{n} = v_n \prod\nolimits_{l=1}^{n-1}(1-v_{l}) \quad G=\sum\nolimits_{n=1}^{+\infty}\gamma_{n}\delta_{\phi_{n}} 
\end{split}
\end{equation}
where $\mathrm{Beta}$ is the Beta distribution. $\delta_{\phi_{n}}$ denotes the Dirac delta function at $\phi_{n}$, and $\gamma_{n}$ represents the weight of the mass point $n$. For ease of presentation, we denote the collection of mixture weights as $\bm{\gamma} = \{ \gamma_1,\gamma_2,...\}$. According to (\ref{e3.1.a1}), we know that $\boldsymbol{\gamma}$ only depend on scale parameter $\alpha$. As a result, we aggregate the sample process of $\boldsymbol{\gamma}$ by $\boldsymbol{\gamma} \sim \mathrm{GEM}(\alpha)$, where $\mathrm{GEM}$ represents the stick-breaking distribution. Therefore, (\ref{e3.1.a2}) can be rewritten as follows.
\begin{equation}    \label{e3.1.a3}
\begin{split}
&\phi_n \sim H  \qquad \boldsymbol{\gamma} \sim \mathrm{GEM}(\alpha) \qquad G=\sum\nolimits_{n=1}^{+\infty}\gamma_{n}\delta_{\phi_{n}}
\end{split}
\end{equation}
Combining (\ref{e3.1.a1}) and (\ref{e3.1.a3}), the distribution of the data point $\omega$ is an infinite Bayesian mixture as follows.
\begin{equation}    \label{e3.1.a4}
\begin{split}
p(\omega) = \int F(\phi) dG = \sum_{n=1}^{\infty} \gamma_{n} F(\phi_{n})
\end{split}
\end{equation}
Given the weighted sum expression of discrete distribution $G$ as in (\ref{e3.1.a3}), drawing a sample $\phi_{n}$ from $G$ can be viewed as picking a component whose probability is $\gamma_{n}$. Let the latent variable $z$ denote the selected mass point. In this paper, each component corresponds to a distinct obstacle motion mode,  which represents one possible stochastic behavior of the obstacle over a single sampling interval. Since each component $\phi_{n}$ is chosen with probability $\gamma_{n}$, the variable $z$ follows a multinomial distribution over the weight vector $\boldsymbol{\gamma}$, denoted as follows.
\begin{equation}    \label{e3.1.a5}
\begin{split}
z \sim \mathrm{Mult}(\boldsymbol{\gamma})
\end{split}
\end{equation}
Finally, by substituting (\ref{e3.1.a3}) and (\ref{e3.1.a5}) into (\ref{e3.1.a1}), the complete generative process of the DPMM can be reformulated as follows.
\begin{equation} \label{e3.1.3}
\begin{array}{cc}
\bm{\gamma} \sim {\rm GEM}(\alpha)&\phi_n \sim H\\
z \sim {\rm Mult}(\bm{\gamma})&\qquad \omega \sim F(\phi_z)
\end{array}
\end{equation}
The reader is referred to \cite{sethuraman1994constructive} for a more comprehensive introduction to the DPMM.
\begin{rem}\label{rem3.1.1}
    Compared with other methods, the DPMM not only enables learning the number of mixture components in the underlying distribution but also provides posterior estimates for each component. On the one hand, by automatically inferring the number of mixture components from the data rather than predefining it, the DPMM can adapt to the complexity of the data structure. On the other hand, the posterior distribution is utilized to extract fine-grained local distribution features, which are subsequently used to construct the ambiguity set.
\end{rem}
\subsubsection{Online variational inference algorithm of DPMM}  \label{section3.1.2} 
The previous subsection has outlined the fundamental principles of the DPMM. In practice, however, the objective is not to sample from the generative process, but to infer the posterior distribution of random variable $\omega$ given a sequence of observations $\bm{\hat{\omega}}$. This is achieved via the online variational inference algorithm for the DPMM \cite{gomes2008incremental}, which updates the posterior distribution incrementally as new data arrives, while maintaining computational and memory efficiency. Each learning round consists of two phases: a model-building phase and a compression phase.

In the model-building phase, Variational Bayes (VB) is used to approximate the intractable posterior with a tractable proxy distribution by maximizing the variational free energy as follows.
\begin{equation}    \label{e3.2.a1}
\begin{split}
\mathcal{F}(\boldsymbol{\hat{\omega}}; q) = \int q(V, \Phi, Z) \log \frac{p(V, \Phi, Z, \boldsymbol{\hat{\omega}} | \alpha, H)}{q(V, \Phi, Z)}
\end{split}
\end{equation}
Here, $q(V, \Phi, Z)$ denotes the variational distribution over the stick-breaking weights $V$, component parameters $\Phi$, and assignment variables $Z$, while $p(\cdot)$ denotes the true joint distribution. This optimization is conducted under the clump constraints identified from previous rounds, which improves the efficiency and consistency of posterior updates.

In the compression phase, the algorithm reduces memory usage by grouping similar observations into clumps, small sets of statistically similar data points that are treated as a unit during inference. Instead of storing each data point individually, the algorithm summarizes each clump by its sufficient statistics, thereby significantly reducing memory overhead. To support clump-level inference, a modified variational free energy is introduced as follows.
\begin{align}
    \mathcal{F}_{C} = & -\sum_{k=1}^{K} KL(q(v_{k}) || p(v_{k}|\alpha))  \label{e3.2.a2}\\
&-\sum_{k=1}^{K} KL(q(\phi_{k}) || p(\phi_{k}|H)) + \frac{N}{T} \sum_{s}n_{s} \log(S_{sk}) \notag
\end{align}
where $KL(\cdot | \cdot)$ denotes the Kullback–Leibler divergence, $K$ is the number of mixture components, $N$ the total data volume, $T$ the number of processed observations, and $n_s$ the size of clump $s$. The term $S_{sk}$ measures the statistical affinity between clump $s$ and component $k$. Based on this clump-level objective (\ref{e3.2.a2}), the algorithm recursively partitions the data to identify optimal clump structures. It is important to note that clumps are algorithmic structures introduced for efficient online inference and information aggregation, and do not correspond to any physical quantity of obstacle motion. The full procedure is detailed in Algorithm \ref{algorithm2}.
\begin{algorithm}[t] 
    \caption{Online Variational Inference for DPMM Upon Receiving New Observations}
    \label{algorithm2}
        \KwIn{New observations; Clumps from the previous round}
        Model-building phase based on the new observations and the clump constraints from the previous round (\ref{e3.2.a1}) \;
        Each clump is first assigned to the component with the highest responsibility\;
        \While{Respecting the predefined memory budget}{
            \For{$k=1$ to $K$}{
                \If{evaluated(k) = FALSE}{
                    Split partition $k$ into two partitions \;
                    $\Delta BFE(k) = $ change in Equations (\ref{e3.2.a2}) \;
                    $evaluated(k) = TRUE$ \;
                }
            }
            Split partition $\arg\max_{k}\Delta BFE(k)$ \;
            $K = K + 1$ \;
        }
        Retain clumps into next round \;
        Discard summarized data points \;
        \KwOut{Clumps in this round; The proxy distribution $q(\cdot)$}
\end{algorithm}

\subsection{Distance Metrics between Ambiguity Sets}  \label{section03.2} 
Consider the moment-based ambiguity set as follows.
\begin{small}
\begin{equation}    \label{e3.3.a1}
\begin{split}
    \mathbb{D}(\hat{\mu}, \hat{\Sigma}) = \left\{  p
        \left|\begin{array}{l}
              p \in \mathcal{M}(\mathbb{W}_{\omega}) \\
             \left( \mathbb{E}^{ p}[x] - \hat{\mu}  \right)^{T} \hat{\Sigma}^{-1} \left( \mathbb{E}^{ p}[x] - \hat{\mu}  \right) \leq \beta  \\
             \mathbb{E}^{ p}\left[ \left(x-\hat{\mu} \right)\left(x-\hat{\mu} \right)^{T} \right] \preceq \varepsilon \hat{\Sigma} 
        \end{array}\right. \right\}
\end{split}
\end{equation}
\end{small}
where $\hat{\mu}$ and $\hat{\Sigma}$ denote the estimated first and second moments, respectively. Directly measuring the distance between two such ambiguity sets, say $\mathbb{D}_1$ and $\mathbb{D}_2$, is generally intractable due to the infinite-dimensional nature of the set of probability distributions they represent. Fortunately, the distance between two individual distributions can be effectively quantified using metrics such as the Wasserstein distance \cite{panaretos2019statistical}. As a result, the distance between two moment-based ambiguity sets $\mathbb{D}_1 = \mathbb{D}(\hat{\mu}_1, \hat{\Sigma}_1)$ and $\mathbb{D}_2 = \mathbb{D}(\hat{\mu}_2, \hat{\Sigma}_2)$ can be measured by computing the Wasserstein distance between the corresponding nominal Gaussian distributions $\mathcal{N}(\hat{\mu}_1, \hat{\Sigma}_1)$ and $\mathcal{N}(\hat{\mu}_2, \hat{\Sigma}_2)$ as follows.
\begin{align}       
&\mathcal{D}(\mathbb{D}^{1}(\mathbb{W}, \hat{\mu}^{1}, \hat{\Sigma}^{1}), \mathbb{D}^{2}(\mathbb{W}, \hat{\mu}^{2}, \hat{\Sigma}^{2}))=      \label{e3.3.a2}\\
&\Vert \mu_1-\mu_2 \Vert_2^2+tr(\Sigma_1)+tr(\Sigma_2)-2tr(( \Sigma_1^{\frac{1}{2}}\Sigma_2\Sigma_1^{\frac{1}{2}})^{\frac{1}{2}})        \notag
\end{align}

\section{Collaborative online learning for the position distributions of obstacles}  \label{section3}
\begin{figure}[thbp]
\begin{center}
\includegraphics[width=8.4cm]{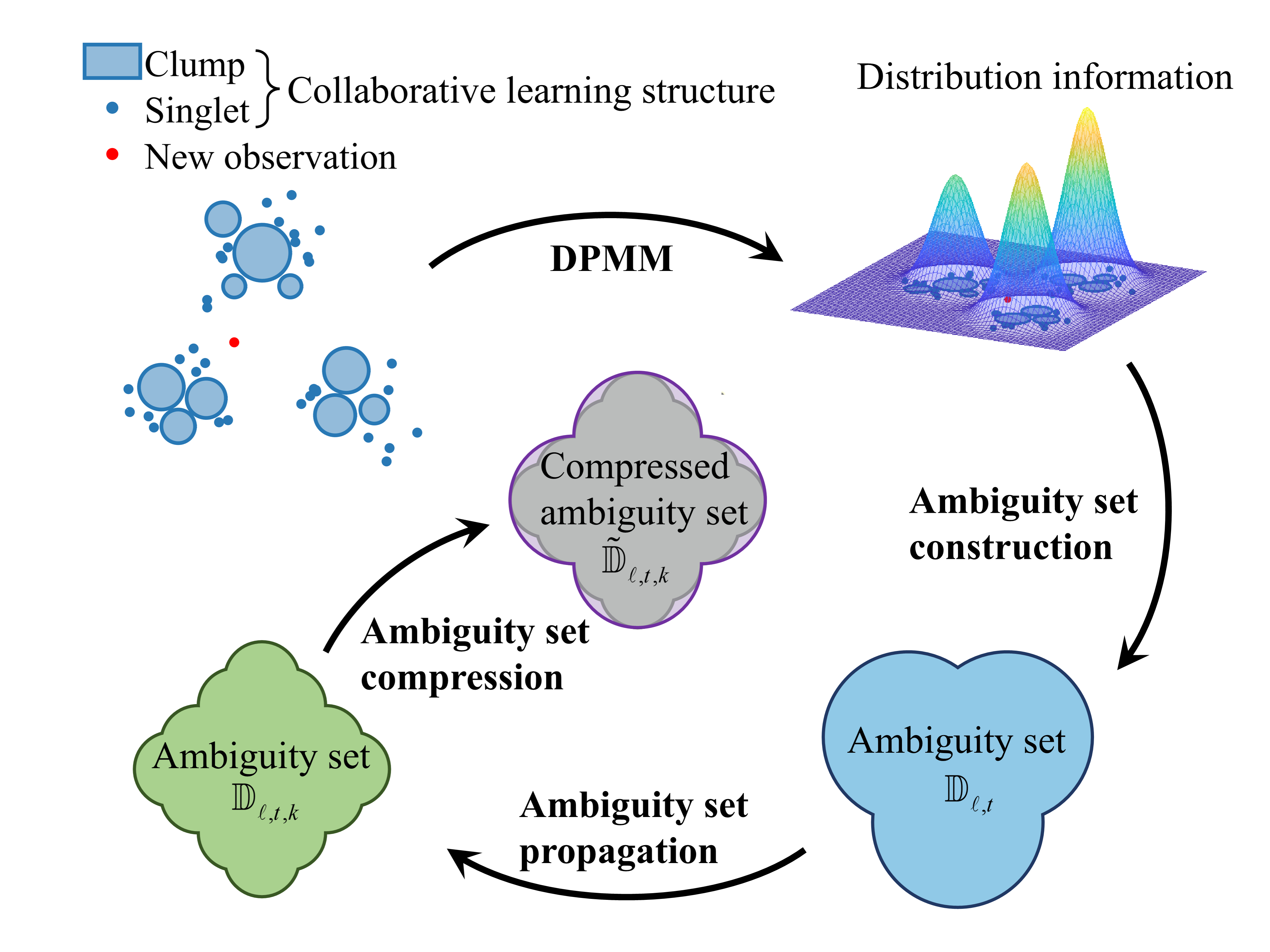}    
\caption{Construction of a compressed ambiguity set}  
\label{fig3.1}                                 
\end{center}                                 
\vspace{-0.3cm}
\end{figure}

To learn the complex uncertainty distribution of $\varphi_{\ell}(k|t)$ and decipher its structural properties, the data stream of the obstacle motion should be extensively exploited. However, the complex environment renders it difficult to obtain real-time data of random motion due to the occlusion between robots and obstacles. To harness the multi-perspective characteristics of multi-robot systems for addressing this challenge, we propose a COOL method to effectively extract distribution information from the decentralized data by exchanging the selected learning structure. We first introduce the COOL method. Subsequently, we construct the ambiguity set for obstacle motion and propagate it to obtain the ambiguity sets for the obstacle positions throughout the prediction horizon. To inhibit the growing complexity of ambiguity sets along the prediction horizon, we propose an flexible ambiguity set compression method and establish its safety guarantee. The process of constructing the compressed ambiguity sets is illustrated in Fig. \ref{fig3.1}.

\subsection{Collaborative online learning method}    \label{section3.1.3}
In this subsection, we introduce a COOL method for multi-robot systems by expanding the online variational inference algorithm introduced in Subsection \ref{section3.1.2}. To begin with, we define the learning structure using the clumps and singlets defined in Section \ref{section3.1.2} as follows.
\begin{defn}    \label{defn3.1.1}
(Learning structure). Suppose that there are $N_c$ clumps summarized by sufficient statistics and $N_s$ singlets, the learning structure is defined as follows.
\begin{equation} \label{e3.1.4}
\begin{split}
S^{tr} = \Big\{ \left( \mu_i, \Sigma_i, N_i \right) \Big| i \in \mathbb{N}_{++}^{N_c} \Big\} \bigcup \Big\{ \hat{\omega}^{(j)} \Big| j \in \mathbb{S} \Big\}
\end{split}
\end{equation}
where $\mu_i$, $\Sigma_i$, and $N_i$ represent the mean, covariance matrix, and the number of the original data in the clump $i$, respectively. $\mathbb{S}$ denotes the indices of $N_s$ singlets. 
\end{defn}
The learning structure encapsulates all the information necessary for the online variational inference, which provides a foundation for the continuous operation of the online variational algorithm across multiple robots. The number of the total learning data $N^{dat}$ is calculated as $N^{dat}=N_s+\sum_{i=1}^{N_c}N_i$. We denote the learning structure and the learning data number of robot $i$ regarding obstacle $\ell$ at time $t$ as $S^{tr}_{i,\ell,t}$ and $N^{dat}_{i,\ell,t}$, respectively. COOL is achieved by choosing and exchanging the selected learning structure among robots. Note that the selected learning structure can be flexibly selected based on a priori knowledge of the obstacle motion. In this work, we choose the selected learning structure as the learning structure with the largest number of learning data represented as $S^{tr}_{j,\ell,t}$, where $j = \mathop{\arg\max}_{j\in \mathbb{N}_{++}^{N}} N^{dat}_{j,\ell,t}$.

The implementation steps of the proposed method are as follows. If robot $i$ observes $\hat{\omega}_{\ell,t}$, it will request $S^{tr}_{k,\ell,t-1}$ from all other robots $k \neq i$. To reduce the amount of data transmitted, other robots only return $N^{dat}_{k,\ell,t-1}$. Then robot $i$ will re-establish communication with the robot having the selected learning structure (e.g. robot $j$) and receive $S^{tr}_{j,\ell,t-1}$ to replace $S^{tr}_{i,\ell,t-1}$. Finally, robot $i$ performs the online variational inference algorithm using $\hat{\omega}_{\ell,t}$ to obtain the posterior distribution and the updated learning structure $S^{tr}_{i,\ell,t}$. It is clear that each learning for $\omega_{\ell,t}$ is based on the current selected learning structure of the system and further improves this learning structure.

Note that we choose $F(\cdot)$ as the Gaussian distribution family, and thus at time $t$, the learning result for the obstacle motion $\omega_{\ell,t}$ includes the number of the mixture components $m_{\ell,t}$ and for each component $i$: (1) the number of learning data $N_{\ell,t}^{i}$, (2) the estimated weight $\hat{\gamma}_{\ell,t}^{i}$, and (3) the estimated mean $\hat{\mu}_{\ell,t}^{i}$ and covariance $\hat{\Sigma}_{\ell,t}^{i}$ which are essentially computed by (\ref{e3.1.7}).
\begin{equation} \label{e3.1.7}
\begin{split}
&\hat{\mu}_{\ell,t}^{i} = \frac{1}{N_{\ell,t}^{i}}\sum\nolimits_{j=1}^{N_{\ell,t}^{i}} \hat{\omega}_{\ell,t}^{i(j)} \\
&\hat{\Sigma}_{\ell,t}^{i} = \frac{1}{N_{\ell,t}^{i}}\sum\nolimits_{j=1}^{N_{\ell,t}^{i}}(\hat{\omega}_{\ell,t}^{i(j)} - \hat{\mu}_{\ell,t}^{i}) (\hat{\omega}_{\ell,t}^{i(j)} - \hat{\mu}_{\ell,t}^{i})^T
\end{split}
\end{equation}
where $\hat{\omega}_{\ell,t}^{i(j)}$ is the learning data from component $i$.
\begin{rem} \label{rem3.1.3}
     Note that a key challenge in COOL is how to efficiently extract distribution information from decentralized data streams. Due to the complexity of the positions of robots and obstacles, the observed data from different robots may overlap but are not entirely identical. Therefore, we adopt the online variational inference algorithm that uses the learning structure to store the extracted data information. By exchanging the selected learning structure among robots, the distribution information in decentralized data streams can be fully exploited while avoiding redundant computations.
\end{rem}

\subsection{Ambiguity set construction and propagation}     \label{section3.2}
Based on the results of COOL, we develop an online-updated ambiguity set for $\omega_{\ell,t}$. This ambiguity set is designed to characterize the uncertainty in the probability distribution of the obstacle motion vector $\omega_{\ell,t}$ learned from data, and serves as the foundation for constructing distributionally robust collision avoidance constraints. The definition of this ambiguity set is first stated explicitly, followed by a detailed explanation.
\begin{small}
\begin{subequations}    \label{e3.2.1}
    \begin{alignat}{2}
        &\mathbb{D}_{\ell,t} := \left\{ \sum_{i=1}^{m_{\ell,t}} \gamma_{\ell,t}  p^{i}
        \left|\begin{array}{l}
             \boldsymbol{\gamma_{\ell,t}} \in \Delta_{\ell,t}(\boldsymbol{\hat{\gamma}_{\ell,t}})  \\
              p^{i} \in \mathbb{D}^{i} (\mathbb{W}_{\omega}, \hat{\mu}^{i} , \hat{\Sigma}^{i} )
        \end{array}\right. \right\}     \label{e3.2.1a}\\
        &\Delta_{\ell,t} := \left\{ \boldsymbol{\gamma_{\ell,t}} \left| \sum_{i=1}^{m_{\ell,t}}\gamma_{\ell,t}^{i} =1, \sum_{i=1}^{m_{\ell,t}} \vert \gamma_{\ell,t}^{i} -\hat{\gamma}_{\ell,t}^{i}  \vert \leq \theta_{\ell,t} \right. \right\}   \label{e3.2.1b} \\
        &\mathbb{D}^{i}  :=         \label{e3.2.1c} \\
        &\left\{  p
        \left|\begin{array}{l}
              p \in \mathcal{M}(\mathbb{W}_{\omega}) \\
             \left( \mathbb{E}^{p}[x] - \hat{\mu}_{\ell,t}^{i}  \right)^{T} (\hat{\Sigma}_{\ell,t}^{i} )^{-1} \left( \mathbb{E}^{p}[x] - \hat{\mu}_{\ell,t}^{i}  \right) \leq \beta_{\ell,t}^{i}  \\
             \mathbb{E}^{p}\left[ \left(x-\hat{\mu}_{\ell,t}^{i} \right)\left(x-\hat{\mu}_{\ell,t}^{i} \right)^{T} \right] \preceq \varepsilon_{\ell,t}^{i} \hat{\Sigma}_{\ell,t}^{i} 
        \end{array}\right. \right\} \notag
    \end{alignat}
\end{subequations}  
\end{small}
where $\boldsymbol{\gamma_{\ell,t}} = (\gamma_{\ell,t}^{1},...,\gamma_{\ell,t}^{m_{\ell,t}})$ and $\boldsymbol{\hat{\gamma}_{\ell,t}} = (\hat{\gamma}_{\ell,t}^{1},...,\hat{\gamma}_{\ell,t}^{m_{\ell,t}})$. $\mathbb{D}_{\ell,t}$ in (\ref{e3.2.1a}) is the ambiguity set devised as a weighted Minkowski sum of $m_{\ell,t}$ basic ambiguity sets $\mathbb{D} $. Intuitively, $\mathbb{D}_{\ell,t}$ incorporates two layers of robustness: the first accounts for the uncertainty in the estimated weights, modeled by the set $\Delta_{\ell,t}$, and the second addresses the uncertainty in the estimated moments, captured by the structure of the basic ambiguity set $\mathbb{D}_{\ell,t}^{i}$. Specifically, $\Delta_{\ell,t}$ in (\ref{e3.2.1b}) is the confidence region of the weight $\boldsymbol{\gamma_{\ell,t}}$, where $\theta_{\ell,t}$ is a parameter indicating the total variation distance. For the definition of basic ambiguity set $\mathbb{D}_{\ell,t}^{i}$ in (\ref{e3.2.1c}), the second constraint specifies that the first-order moment of any distribution in $\mathbb{D}_{\ell,t}^{i}$ belongs to a confidence region scaled by the parameter $\beta_{\ell,t}^{i}$. Similarly, the third constraint ensures that the second-order moment lies within a confidence region scaled by the parameter $\varepsilon_{\ell,t}^{i}$. $\hat{\mu}_{\ell,t}^{i}$ and $\hat{\Sigma}_{\ell,t}^{i}$ are estimated first-order and second-order moments, respectively. $\mathcal{M}(\mathbb{W}_{\omega})$ is the set of positive Borel measures on $\mathbb{W}_{\omega}$. From a physical perspective, the ambiguity set $\mathbb{D}_{\ell,t}$ represents a family of plausible probabilistic models for the obstacle displacement $\omega_{\ell,t}$, capturing the uncertainty in both the learned motion modes and their relative likelihoods.

There are several advantages of the ambiguity set. First, by leveraging the local moment information, the ambiguity set is well-suited for describing multimode structures. Moreover, each basic ambiguity set is devised using mean and covariance, which endows the resulting control method with enormous computational benefits. The ambiguity set is general enough to encompass a wide class of existing ambiguity sets. For example, it reduces to the data-driven ambiguity set \cite{delage2010distributionally}, when $m_{\ell,t} = 1$. To obtain the finite-sample guarantee of the ambiguity set, we make the following assumption.
\begin{assum}   \label{assum4.1.1}
    Assume that the underlying distribution $p_{\ell,t}^{*}$ is time-invariant with $m_{\ell,t}$ component and the date $\hat{\omega}_{\ell}^{i(j)}$ is from its $i$-th component.
\end{assum}
Assumption \ref{assum4.1.1} is reasonable, as powerful nonparametric learning tools such as DPMM can effectively identify the number of modes in a multimodal distribution and assign data points to their corresponding components when sufficient data is available. With the above assumption, we provide the finite-sample guarantee of the ambiguity set below.
\begin{thm} \label{thm3.2.1}
    (Finite-sample guarantee) Suppose Assumption \ref{assum4.1.1} holds, given the confidence level $\chi_{\ell}$ and $\alpha_{\ell}^{i}$, let $\hat{\alpha}_{\ell,t}^{i} = 1-(\alpha_{\ell}^{i})^{\frac{1}{2}}$, $\overline{R}_{\ell,t}^{i} = (1-(( \hat{R}_{\ell,t}^{i} )^{2} + 2) (2 + (2\ln ( 4/\hat{\alpha}_{\ell}^{i} ))^{\frac{1}{2}}) / (N_{\ell,t}^{i})^{\frac{1}{2}})^{-\frac{1}{2}} \hat{R}_{\ell,t}^{i}$ and $\hat{R}_{\ell,t}^{i}  = \sup_{\omega \in \mathbb{W}_{\omega}}\Vert (\hat{\Sigma}_{\ell,t}^{i} )^{-\frac{1}{2}}(\omega - \hat{\mu}_{\ell,t}^{i} ) \Vert_2$, while the parameters $\theta_{\ell,t}$, $\beta $ and $\varepsilon $ of the ambiguity set $\mathbb{D}_{\ell,t}$ are defined as follows (the subscripts of $\hat{\alpha}_{\ell}^{i}$, $\overline{R}_{\ell,t}^{i}$, $\hat{R}_{\ell,t}^{i}$, and part of $N_{\ell,t}^{i}$ are omitted)
    \begin{small}
    \begin{align}
        &\theta_{\ell,t} = 2\sqrt{\frac{m_{\ell,t}\ln2 - \ln (1 - \chi_{\ell})}{2\left( \sum_{i=1}^{m_{\ell,t}}N_{\ell,t}^{i}  - 1\right)}}     \label{e3.2.2}\\
        &\beta  = \frac{\frac{\overline{R}^{2}}{N} \left( 2+\sqrt{2\ln(1/\hat{\alpha})} \right)^2}{1- \frac{\overline{R}^{2}}{\sqrt{N}} \left( \sqrt{1-\frac{N}{\overline{R}^4}} + \sqrt{\ln(\frac{1}{\hat{\alpha}})} \right) - \frac{\overline{R}^{2}}{N}\left( 2+ \sqrt{2\ln(\frac{1}{\hat{\alpha}})} \right)^2}      \notag\\
        &\varepsilon  = \frac{1+ \frac{\overline{R}^{2}}{N} \left( 2+\sqrt{2\ln(1/\hat{\alpha})} \right)^2}{1- \frac{\overline{R}^{2}}{\sqrt{N}} \left( \sqrt{1-\frac{N}{\overline{R}^4}} + \sqrt{\ln(\frac{1}{\hat{\alpha}})} \right) - \frac{\overline{R}^{2}}{N}\left( 2+ \sqrt{2\ln(\frac{1}{\hat{\alpha}})} \right)^2}     \notag
    \end{align}
    \end{small}
    If $N $ satisfies the following condition
    \begin{small}
    \begin{align}
        N  > \max\left\{ \begin{array}{l}
             \left(\hat{R}^2+2\right)^2 \left(2+\sqrt{2\ln(1/\hat{\alpha})}\right),  \\
             \left( 8+32\sqrt{32\ln(4/\hat{\alpha})} \right)^2 \big/ \left( \sqrt{\hat{R}+4} - \hat{R} \right)^4 
        \end{array} \right\}    \notag
    \end{align}
    \end{small}
    then we have the following guarantee.
    \begin{small}
    \begin{equation}    \label{e3.2.3}
    \begin{split}
        \mathbb{P}\{ p_{\ell,t}^* \in \mathbb{D}_{\ell,t} \} \geq \chi_{\ell} \prod\nolimits_{i=1}^{m_{\ell,t}} \alpha_{\ell}^{i} = \alpha_{\ell}
    \end{split}
    \end{equation}        
    \end{small}
\end{thm}
Its proof is provided in Appendix \ref{appendix A}. Theorem \ref{thm3.2.1} (Finite-sample guarantee) states that the probability that the unknown underlying distribution $p_{\ell,t}^{*}$ lies within the ambiguity set $\mathbb{D}_{\ell,t}$ is no less than a prescribed confidence level $\alpha_{\ell}$. Notably, Theorem \ref{thm3.2.1} is distribution-agnostic in the sense that it does not require the underlying distribution $p_{\ell,t}^{*}$ to follow any specific form, thereby accommodating non-Gaussian distributions as well. However, constraining the collision probability across the prediction horizon requires the ambiguity set for $\varphi_{\ell}(k|t)$. According to (\ref{e2.1.4}), $\varphi_{\ell}(k|t)$ is the sum of $\varphi_{\ell}(t)$ and $k$ i.i.d. random variables obeying the distribution in (\ref{e3.2.1}) with high probability. We then propagate $\mathbb{D}_{\ell,t}$ over the prediction horizon to obtain the ambiguity set for $\varphi_{\ell}(k|t)$ through the following Lemma, whose proof is given in Appendix \ref{appendix B}.
\begin{lem}         \label{lem3.2.1}
    (Ambiguity set propagation) Assume that the underlying distribution $p_{\ell,t}^{*}$ lies in the ambiguity set $\mathbb{D}_{\ell,t}$ in (\ref{e3.2.1}), the distribution $p_{\varphi_{\ell}(k|t)}^{*} = \delta_{\varphi_{\ell}(t)} * [p_{\ell,t}^{*}]^{k}$ lies in the following ambiguity set.
    \begin{small}
    \begin{align}   
        &\mathbb{D}_{\ell,t,k} := \left\{ \sum_{j=1}^{m_{\ell,t,k}} \gamma_{\ell,t,k}^{j}  p^{j}
        \left|\begin{array}{l}
             \boldsymbol{\gamma_{\ell,t,k}} \in \Delta_{\ell,t,k}(\boldsymbol{\hat{\gamma}_{\ell,t,k}})  \notag\\
              p^{j} \in \mathbb{D}_{\ell,t,k}^{j}(\mathbb{W}_{\ell,t+k}, \hat{\mu}_{\ell,t,k}^{j}, \hat{\Sigma}_{\ell,t,k}^{j})
        \end{array}\right. \right\} \\
        &\Delta_{\ell,t,k} := \left\{ \boldsymbol{\gamma_{\ell,t,k}} \left|\begin{array}{l}
             \sum_{j=1}^{m_{\ell,t,k}}\gamma_{\ell,t,k}^{j}=1  \\
             \sum_{j=1}^{m_{\ell,t,k}} \vert \gamma_{\ell,t,k}^{j}-\hat{\gamma}_{\ell,t,k}^{j} \vert \leq \theta_{\ell,t,k} 
        \end{array} \right. \right\}    \notag\\
        &\mathbb{D}_{\ell,t,k}^{j} :=   \label{e3.2.4}\\
        &\left\{  p
        \left|\begin{array}{l}
              p \in \mathcal{M}(\mathbb{W}_{\ell,t+k})   \\
             \left( \mathbb{E}^{ p}[x] - \hat{\mu}_{\ell,t,k}^{j} \right)^{T} (\hat{\Sigma}_{\ell,t,k}^{j})^{-1} \left( \mathbb{E}^{ p}[x] - \hat{\mu}_{\ell,t,k}^{j} \right) \leq \beta_{\ell,t,k}^{j} \\
             \mathbb{E}^{ p}\left[ \left(x-\hat{\mu}_{\ell,t,k}^{j}\right)\left(x-\hat{\mu}_{\ell,t,k}^{j}\right)^{T} \right] \preceq \varepsilon_{\ell,t,k}^{j}\hat{\Sigma}_{\ell,t,k}^{j}
        \end{array}\right. \right\} \notag
    \end{align}
    \end{small}
    where
    \begin{subequations}    \label{e3.2.5}
    \begin{alignat}{2}
        &m_{\ell,t,k} = \binom{m_{\ell,t}-1}{k+m_{\ell,t}-1} \label{e3.2.5a}\\ 
        &\hat{\gamma}_{\ell,t,k}^{j} = (k!) \big/ (k_1^j!k_2^j!...k_{m_{\ell,t}}^j!) \prod\nolimits_{i=1}^{m_{\ell,t}}(\hat{\gamma}_{\ell,t}^i)^{k_i^j}  \label{e3.2.5b}\\
        &\hat{\mu}_{\ell,t,k}^{j} = \varphi_{\ell}(t)+\sum\nolimits_{i=1}^{m_{\ell,t}}k_i^j\hat{\mu}_{\ell,t}^i  \label{e3.2.5c}\\
        &\hat{\Sigma}_{\ell,t,k}^{j} = \sum\nolimits_{i=1}^{m_{\ell,t}}k_i^j\hat{\Sigma}_{\ell,t}^i \label{e3.2.5d}\\
        &\theta_{\ell,t,k} = k\theta_{\ell,t}(1+2\theta_{\ell,t})^{k-1} \label{e3.2.5e}\\
        &\beta_{\ell,t,k}^{j} = \sum\nolimits_{i=1}^{m_{\ell,t}}k_i^j\beta   \label{e3.2.5f}\\
        &\varepsilon_{\ell,t,k}^{j} = \beta_{\ell,t,k}^{j} + \max_{i=1,...,m_{\ell,t}} \varepsilon_{\ell,t}^{i}   \label{e3.2.5g}\\
        &\sum\nolimits_{i=1}^{m_{\ell,t}}k_i^j = k, \quad k_i^j \in \mathbb{N}_{+}   \label{e3.2.5h}
    \end{alignat}
    \end{subequations}
    where $\mathbb{W}_{\ell,t}$ is the support set of the position of obstacle $\ell$ at time $t$, $\binom{\cdot}{\cdot}$ denotes the number of combinations, and $k_{1}^{j},...,k_{m_{\ell,t}}^{j}$ is a particular delineation of $k$ with index $j$.
\end{lem}
\begin{rem} \label{rem3.2.1}
    In (\ref{e3.2.5}), $m_{\ell,t,k}$ essentially represents the number of non-negative integer-solution groups to $\sum_{i=1}^{m_{\ell,t}}k_i=k$, which is calculated through (\ref{e3.2.5a}). And $k_{1}^{j},...,k_{m_{\ell,t}}^{j}$ denote the $j$-th group of solutions.
\end{rem}

\subsection{Ambiguity set compression}     \label{section3.3}
Since the reformulation of the DR-CVaR constraint (\ref{e2.3.2.5}) only requires an ambiguity set over the obstacle position, the ambiguity set (\ref{e3.2.4}) can be directly employed to construct (\ref{e2.3.2.5}). However, as indicated by (\ref{e3.2.5a}), the number of basic ambiguity sets $m_{\ell,t,k}$ grows rapidly with $k$, which leads to a high computational burden in practice. To address this, we introduce a compression scheme that limits the total number of basic ambiguity sets to a user-defined maximum $M_{\ell,t,k}$ by aggregating similar components.

To mitigate the impact of compression on control performance, it is important to ensure that only similar ambiguity sets are merged. As discussed in Subsection \ref{section03.2}, the distance between two basic ambiguity sets can be computed using appropriate metrics (e.g., Wasserstein distance). This naturally leads to the use of clustering techniques to group the $m_{\ell,t,k}$ ambiguity sets in (\ref{e3.2.4}) into $M_{\ell,t,k}$ clusters. We do not impose any specific requirement on the choice of clustering algorithm.

For each group $n\in\mathbb{N}_{++}^{M_{\ell,t,k}}$, we aim to compress all the basic ambiguity sets in this group into a single basic ambiguity set $\widetilde{\mathbb{D}}_{\ell,t,k}^{n} \left( \mathbb{W}_{\ell,t+k}, \hat{\widetilde{\mu}}_{\ell,t,k}^n, \hat{\widetilde{\Sigma}}_{\ell,t,k}^n \right)$. However, the uncertainty of the weight in (\ref{e3.2.4}) poses significant challenges to the compression process. Fortunately, we propose the following safe compression of ambiguity set.
\begin{lem}  \label{lem3.3.1}
    (Safe compression of ambiguity set) For the ambiguity set (\ref{e3.2.4}) and any given grouping $\mathbb{M}_{\ell,t,k}^{1},..., \mathbb{M}_{\ell,t,k}^{M_{\ell,t,k}}$, if the compressed ambiguity set $\widetilde{\mathbb{D}}_{\ell,t,k}$ is calculated by
    \begin{small}
    \begin{align}
        &\widetilde{\mathbb{D}}_{\ell,t,k} := \left\{ \sum_{n=1}^{M_{\ell,t,k}} \widetilde{\gamma}_{\ell,t,k}^{n} \widetilde{ p}_{\ell,t,k}^{n}
        \left|\begin{array}{l}
             \boldsymbol{\widetilde{\gamma}_{\ell,t,k}} \in \widetilde{\Delta}_{\ell,t,k}  \notag\\
             \widetilde{ p}_{\ell,t,k}^{n} \in \widetilde{\mathbb{D}}_{\ell,t,k}^{n}
        \end{array}\right. \right\}  \\
        &\widetilde{\Delta}_{\ell,t,k} := \left\{ \boldsymbol{\widetilde{\gamma}_{\ell,t,k}} \left|\begin{array}{l}
             \sum_{n=1}^{M_{\ell,t,k}} \widetilde{\gamma}_{\ell,t,k}^{n}=1  \\
             \sum_{n=1}^{M_{\ell,t,k}} \vert \widetilde{\gamma}_{\ell,t,k}^{n}-\hat{\widetilde{\gamma}}_{\ell,t,k}^{n} \vert \leq \widetilde{\theta}_{\ell,t,k} 
        \end{array} \right. \right\}        \notag\\
        &\widetilde{\mathbb{D}}_{\ell,t,k}^{n} :=           \label{e3.3.2.1}\\
        &\left\{  p
        \left|\begin{array}{l}
              p \in \mathcal{M}(\mathbb{W}_{\ell,t+k}) \\
             \left( \mathbb{E}^{ p}[x] - \hat{\widetilde{\mu}}_{\ell,t,k}^{n} \right)^{T} (\hat{\widetilde{\Sigma}}_{\ell,t,k}^{n})^{-1} \left( \mathbb{E}^{ p}[x] - \hat{\widetilde{\mu}}_{\ell,t,k}^{n} \right) \leq \widetilde{\beta}_{\ell,t,k}^{n} \\
             \mathbb{E}^{ p}\left[ \left(x-\hat{\widetilde{\mu}}_{\ell,t,k}^{n}\right)\left(x-\hat{\widetilde{\mu}}_{\ell,t,k}^{n}\right)^{T} \right] \preceq \hat{\widetilde{\Phi}}_{\ell,t,k}^{n}
        \end{array}\right. \right\}     \notag
    \end{align}
    \end{small}
    where (the subscripts $\ell,t,k$ are omitted)
    \begin{small}
    \begin{subequations}    \label{e3.3.2.2}
    \begin{alignat}{2}
        &\hat{\widetilde{\gamma}}^{n} = \sum\nolimits_{j \in \mathbb{M}^{n}} \hat{\gamma}^{j}, \quad \widetilde{\theta} = \theta  \label{e3.3.2.2a}\\
        &\hat{\widetilde{\mu}}^{n} = \left(\sum\nolimits_{j \in \mathbb{M}^{n}} \hat{\gamma}^{j} \hat{\mu}^{j}\right) \big/ \left(\sum\nolimits_{j \in \mathbb{M}^{n}} \hat{\gamma}^{j}\right) \label{e3.3.2.2b}\\ 
        &\widetilde{\beta}^{n} = \sum\nolimits_{j \in \mathbb{M}^{n}} \left\{ \overline{\gamma}^{n} \beta^{j} + \breve{\gamma}^{n} \right\} \label{e3.3.2.2c}\\
        &\hat{\widetilde{\Sigma}}^{n} = \sum\nolimits_{j \in \mathbb{M}^{n}} \left\{ \overline{\gamma}^{n} \hat{\Sigma}^{j} + \breve{\gamma}^{n} \left(\hat{\mu}^{j}\right) \left(\hat{\mu}^{j}\right)^T \right\} \label{e3.3.2.2d}\\
        &\hat{\widetilde{\Phi}}^{n} = 4\hat{\beta}^{n} \hat{\widetilde{\Sigma}}^{n} +  \label{e3.3.2.2e}\\
        &\sum_{j \in \mathbb{M}^{n}} \overline{\gamma}^{n} \left\{ \left( \varepsilon^{j} + 3\beta^{j} \right) \hat{\Sigma}^{j} + 3\left( \hat{\mu}^{j} - \hat{\widetilde{\mu}}^{n} \right) \left( \hat{\mu}^{j} - \hat{\widetilde{\mu}}^{n} \right)^T \right\}  \notag\\
        &\overline{\gamma}^{n} = \max_{\boldsymbol{\gamma} \in \Delta} \frac{\gamma^{j}}{\sum_{j \in \mathbb{M}^{n}} \gamma^{j}} \label{e3.3.2.2f}\\ 
        &\breve{\gamma}^{n} = \max_{\boldsymbol{\gamma} \in \Delta} \left| \frac{\gamma^{j}}{\sum_{j \in \mathbb{M}^{n}} \gamma^{j}} - \frac{\hat{\gamma}^{j}}{\sum_{j \in \mathbb{M}^{n}} \hat{\gamma}^{j}} \right| \label{e3.3.2.2g}
    \end{alignat}
    \end{subequations}
    \end{small}
    then we have $\mathbb{D}_{\ell,t,k} \subseteq \widetilde{\mathbb{D}}_{\ell,t,k}$.
\end{lem}
Its proof is provided in Appendix \ref{appendix C}. Lemma \ref{lem3.3.1} enables flexible adjustment of $M_{\ell,t,k}$. It is important to note that a larger value of $M_{\ell,t}$ leads to more detailed local moment information, which we refer to as ``granularity''. Similarly, the ambiguity set $\widetilde{\mathbb{D}}_{\ell,t,k}$ is composed of a greater number of basic ambiguity sets, which we refer to as its ``complexity''. Therefore, Lemma \ref{lem3.3.1} effectively enables flexible adjustment of the complexity and granularity of $\widetilde{\mathbb{D}}_{\ell,t,k}$. Additionally, Lemma \ref{lem3.3.1} illustrates that $\mathbb{D}_{\ell,t,k}$ is definitely contained within $\widetilde{\mathbb{D}}_{\ell,t,k}$, which ensures the safety of constructing the DR-CVaR constraints (\ref{e2.3.2.5}) using $\widetilde{\mathbb{D}}_{\ell,t,k}$.
\begin{rem} \label{rem3.4.1}
    Note that a key challenge of the proposed framework lies in providing the probabilistic collision avoidance guarantee under uncertain distributions. To address this challenge, we propose an ambiguity set based on the confidence region and provide the finite-sample guarantee for this ambiguity set using the technique of data-driven DRO. Additionally, we propose a general compression method in (\ref{e3.3.2.2}), which can flexibly adjust the complexity of the ambiguity set with a safe guarantee. Together, the probabilistic collision avoidance guarantee is finally obtained in Theorem \ref{thm4.4.1} (see Section \ref{section4.3}).
\end{rem}

\section{The COOL-DRMC framework for multi-robot motion control}  \label{section4}
\begin{figure}[htbp]
\begin{center}
\includegraphics[width=8.4cm]{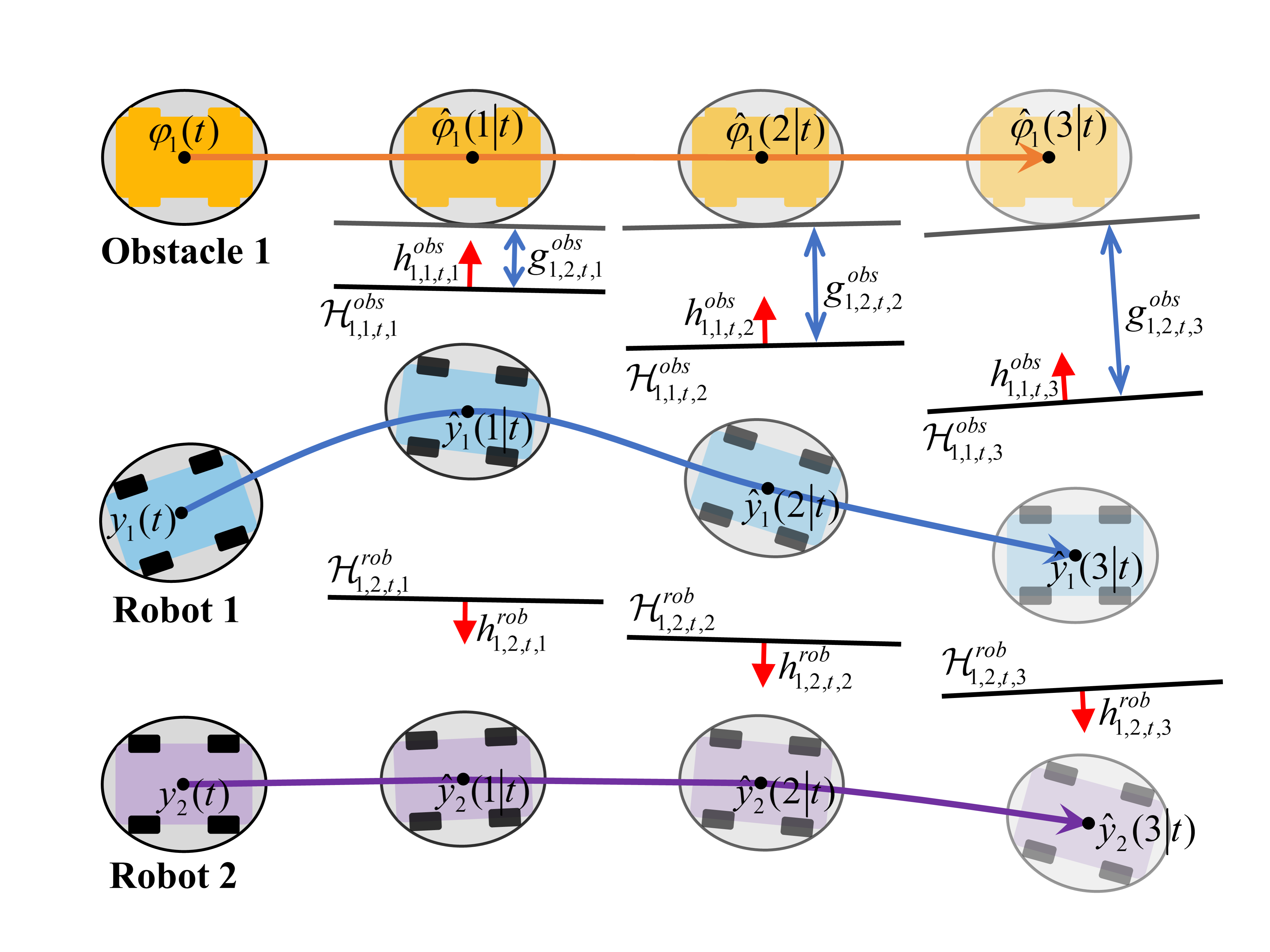}    
\caption{Separating hyperplanes between robot 1 and both obstacle 1 and robot 2 from the view of robot 1}  
\label{fig4.1}                                 
\end{center}                                 
\end{figure}

In this section, we first introduce a spatial allocation protocol for decentralized systems to enable robots to independently compute the separating hyperplane $\mathcal{H}_{i,j,t,k}^{rob}$ while ensuring safety. Additionally, we utilize the DRO technique to transform the computation of the separating hyperplane $\mathcal{H}_{i,\ell,t,k}^{obs}$ into a SDP problem based on $\widetilde{\mathbb{D}}_{\ell,t,k}$. Subsequently, we summarize the COOL-DRMC algorithm and present its performance analysis. An example of the separating hyperplanes between robot 1 and both obstacle 1 and robot 2 from the view of robot 1 at time $t$ is illustrated in Fig. \ref{fig4.1}

\subsection{Inter-robot separating hyperplane computation}  \label{section4.1}
We consider the spatial allocation protocol between robots $i$ and $j$ at prediction time $k|t$ aimed at enabling robots $i$ and $j$ to independently compute the separating hyperplanes $\mathcal{H}_{i,j,k,t}^{rob}$ and $\mathcal{H}_{j,i,k,t}^{rob}$. By imposing the affine constraint (\ref{e2.4.2.1g}) in each robot perspective, constraints (\ref{e2.3.1.5}) and (\ref{e2.3.1.6}) can be naturally satisfied. We begin introducing the process of obtaining $\mathcal{H}_{i,j,k,t}^{rob}$ with the definition of committed trajectory at time $t$.
\begin{defn}    \label{defn4.1.1}
(Committed trajectory). Given an optimal state $\bm{x_i^*(t-1)}=(x_i^*(0|t-1),...,x_i^*(K|t-1))$ by solving the optimization problem (\ref{e2.4.2.1}) at time $t-1$, the committed trajectory at time $t$ is defined by
\begin{equation}    \label{e4.1.1}
\bm{\hat{y}_i(t)} = (\hat{y}_i(0|t),...,\hat{y}_i(K|t))
\end{equation}
where $\hat{y}_i(k|t)= h(x_i^*(k+1|t-1))$ for $k \in \mathbb{N}_{+}^{K-1}$, and $\hat{y}_i(K|t)= h(f(x_i^*(K|t-1), \kappa_i^f(x_i^*(K|t-1))))$ is obtained through a terminal control law $\kappa_i^f(\cdot)$. The committed trajectory at $t=0$ is set as the reference trajectory $h(r_i^x)$.
\end{defn}
For robot $i$, it is clear that the committed trajectory $\bm{\hat{y}_i(t)}$ is obtained before solving (\ref{e2.4.2.1}). Moreover, based on Assumption \ref{assum2.3.1}, robot $i$ can also obtain $\bm{\hat{y}_j(t)}$ by communicating without delay with robot $j$. Then robot $i$ calculates $\mathcal{H}_{i,j,k,t}^{rob}$ independently based on $\bm{\hat{y}_i(t)}$ and $\bm{\hat{y}_j(t)}$. Specifically, the normal vector of this separation hyperplane is calculated as follows.
\begin{equation} \label{e4.1.2}
\begin{split}
h_{i,j,t,k}^{rob} = \frac{\hat{y}_j(k|t)-\hat{y}_i(k|t)}{\Vert \hat{y}_j(k|t)-\hat{y}_i(k|t) \Vert_2}
\end{split}
\end{equation}
The separating hyperplane is defined as the plane perpendicular to $h_{i,j,t,k}^{rob}$ and passing through the point $(\mathcal{P}_{\hat{y}_i(k|t) \oplus \mathcal{R}_i}(h_{i,j,t,k}^{rob})+\mathcal{P}_{\hat{y}_j(k|t) \oplus \mathcal{R}_j}(-h_{i,j,t,k}^{rob}))/2$. As a result, the bias $g_{i,j,t,k}^{rob}$ can be obtained below.
\begin{align}
  & g_{i,j,t,k}^{rob}        \label{e4.1.3}\\
= & \left( \mathcal{S}_{\hat{y}_j(k|t) \oplus \mathcal{R}_j}(-h_{i,j,t,k}^{rob}) - \mathcal{S}_{\hat{y}_i(k|t) \oplus \mathcal{R}_i}(h_{i,j,t,k}^{rob})  \right) \big/ 2    \notag
\end{align}
We note that even when the two convex sets $\hat{y}_{i}(k|t) \oplus \mathcal{R}_{i}$ and $\hat{y}_{j}(k|t) \oplus \mathcal{R}_{j}$ are disjoint, the separating hyperplane obtained via the specific construction in (\ref{e4.1.2}) and (\ref{e4.1.3}) may fail to separate them. In this case, we directly use the separating hyperplane from the previous time step, i.e., $\mathcal{H}_{i,j,t,k}^{rob} = \mathcal{H}_{i,j,t-1,k+1}^{rob}$. In particular, for the case $k=K-1$, since $\mathcal{H}_{i,j,t-1,K}^{rob}$ does not exist, we compute $\mathcal{H}_{i,j,t,K-1}^{rob}$ as the separating hyperplane between $\hat{y}_{i}(K-1|t) \oplus \mathcal{R}_{i}$ and $\hat{y}_{j}(K-1|t) \oplus \mathcal{R}_{j}$. In summary, robots $i$ and $j$ are always separated into two half-spaces.

\subsection{Robot-obstacle separating hyperplane computation}  \label{section4.2}
We utilize the DRO technique to equivalently transform the computation of $\mathcal{H}_{i,\ell,t,k}^{obs}$ satisfied the DR-CVaR constraint (\ref{e2.3.2.5}) into an SDP problem. Based on $\bm{\hat{y}_{i}(t)}$, the normal vector $h_{i,\ell,t,k}^{obs}$ is first calculated below.
\begin{equation} \label{e4.2.1}
\begin{split}
h_{i,\ell,t,k}^{obs} = \frac{\hat{\varphi}_{\ell}(k|t)-\hat{y}_i(k|t)}{\Vert \hat{\varphi}_{\ell}(k|t)-\hat{y}_i(k|t) \Vert_2}
\end{split}
\end{equation}
where $\hat{\varphi}_{\ell}(k|t)$ is the mean CoM of obstacle $\ell$ based on the compressed ambiguity set (\ref{e3.3.2.1}) and calculated as follows.
\begin{equation} \label{e4.2.2}
\begin{split}
\hat{\varphi}_{\ell}(k|t) = \sum\nolimits_{n=1}^{M_{\ell,t,k}} \hat{\widetilde{\gamma}}_{\ell,t,k}^n \hat{\widetilde{\mu}}_{\ell,t,k}^n
\end{split}
\end{equation}
Note that with a given normal vector $h_{\ell,t,k}$, the decrease in $g_{\ell,t,k}$ results in $\mathcal{H}_{i,\ell,t,k}^{obs}$ moving further away from $\hat{y}(k|t)$, consequently leading to a less conservative trajectory. Therefore, we minimize $g_{i,\ell,t,k}^{obs}$ subject to constraint (\ref{e2.3.2.5}). However, (\ref{e2.3.2.5}) requires that the CVaR condition be satisfied for the worst-case distribution within the ambiguity set $\widetilde{\mathbb{D}}_{\ell,t,k}$ and $\widetilde{\mathbb{D}}_{\ell,t,k}$ involves an infinite number of distributions; thus (\ref{e2.3.2.5}) is infinite-dimensional because it involves an infinite number of distributions in $\widetilde{\mathbb{D}}_{\ell,t,k}$. Fortunately, the optimization problem can be equivalently reformulated as an SDP problem, as derived in the following theorem, where the superscript $obs$ and the subscripts $i$, $\ell$, $t$, and $k$ are suppressed.
\begin{thm}     \label{thm4.2.1}
Suppose that $\mathbb{W}:=\{\varphi \in \mathbb{R}^d \big| E\varphi \leq f\}$ is the support set of the obstacle position $\varphi$, the optimization problem of minimizing $g$ subject to (\ref{e2.3.2.5}) admits the finite-dimensional reformulation as follows.
\vspace{-0.3cm}
\begin{alignat}{2}
    &\min_{g,z,\vartheta,\zeta,r_1^n,r_2^n,s^n,\Lambda^n,\xi^n,\tau^n,\Omega^n,\eta^n,\lambda^n} \quad g    \label{e4.2.3}\\
    &s.t.\ z + \tfrac{1}{1-\alpha^{u}}\left( \hat{\widetilde{\theta}} \zeta + \sum\nolimits_{n=1}^{M} \hat{\widetilde{\gamma}}^{n} (r_1^n-r_2^n) + \vartheta \right) \leq 0      \notag\\
    & \zeta - r_1^n - r_2^n = 0         \notag\\
    & r_1^n - r_2^n + \vartheta \geq \tau^{n}\hat{\widetilde{\beta}}^{n} - 2(\hat{\widetilde{\mu}}^{n})^T\xi^{n} - (\hat{\widetilde{\mu}}^{n})^T \Omega^{n} \hat{\widetilde{\mu}}^{n}     \notag\\
    & \quad\quad\quad\quad\quad\quad + \hat{\widetilde{\Sigma}}^{n} \bullet \Lambda^{n} + \hat{\widetilde{\Phi}}^{n} \bullet \Omega^{n} + s^{n} \notag\\
    & \left(\begin{array}{cc}
      \Omega^{n}   & \dfrac{1}{2} \left(\begin{array}{c}
                                            E^T\eta^{n} - 2\xi^{n}  \\
                                            +h -2\Omega^{n}\hat{\widetilde{\mu}}^{n}
                                      \end{array} \right) \\
      \dfrac{1}{2} \left(\begin{array}{c}
                               E^T\eta^{n} - 2\xi^{n}  \\
                               +h -2\Omega^{n}\hat{\widetilde{\mu}}^{n}
                          \end{array} \right)^T   & \left.\begin{array}{l}
                                                              g-\mathcal{S}_{\mathcal{O}}(-h)+  \\
                                                              z+s^{n}-f^T\eta^n 
                                                         \end{array}\right.
    \end{array}\right) \succeq 0        \notag\\
    & \left(\begin{array}{cc}
      \Omega^{n}   & \dfrac{1}{2}\left(\begin{array}{c}
                                            E^T\lambda^{n} - 2\xi^{n}  \\
                                            -2\Omega^{n}\hat{\widetilde{\mu}}^{n}
                                      \end{array} \right) \\
      \dfrac{1}{2}\left(\begin{array}{c}
                               E^T\lambda^{n} - 2\xi^{n}  \\
                               -2\Omega^{n}\hat{\widetilde{\mu}}^{n}
                          \end{array} \right)^T   & s^{n}-f^T\lambda^n 
    \end{array}\right) \succeq 0        \notag\\
    & \left(\begin{array}{cc}
            \Lambda^{n} & \xi^{n}    \\\
            (\xi^{n})^T & \tau^{n}
        \end{array}\right) \succeq 0, \quad \Omega^{n} \succeq 0        \notag\\
    & \zeta \geq 0, \quad r_1^n \geq 0, \quad r_2^n \geq 0, \quad \eta^n \geq 0, \quad \lambda^n \geq 0     \notag
\end{alignat}
\end{thm}
\vspace{-0.6cm}
Its proof is provided in Appendix \ref{appendix D}. If the hyperplane defined by (\ref{e4.2.1}) and (\ref{e4.2.3}) does not split $\hat{y}_{i}(k|t) \oplus \mathcal{R}_{i}$ and the risk region w.r.t. $\varphi_{\ell}(k|t)$ (i.e., $\hat{y}_{i}(k|t)$ does not satisfy constraint (\ref{e2.4.2.1h}) based on $\mathcal{H}_{i,\ell,t,k}^{obs}$), we use the separating hyperplane from the previous time step $\mathcal{H}_{i,\ell,t,k}^{obs} = \mathcal{H}_{i,\ell,t-1,k+1}^{obs}$. In particular, for the case $k=K-1$, we compute $\mathcal{H}_{i,\ell,t,K-1}$ as the separating hyperplane between $\hat{y}_{i}(K-1|t)\oplus\mathcal{R}_{i}$ and $\mathbb{W}_{\ell,t+K-1}\oplus\mathcal{O}_{\ell}$. For the case $\ell \notin \mathbb{O}_{i,t}$, we use the compressed ambiguity set at the most recent time $t-l$ when $\ell \in \mathbb{O}_{i,t-l}$. If $\widetilde{\mathbb{D}}_{\ell,t-l,k+l}$ does not exist, we directly compute $\mathcal{H}_{i,\ell,t,k}$ as the separating hyperplane between $\hat{y}_{i}(k|t)\oplus\mathcal{R}_{i}$ and $\mathbb{W}_{\ell,t+k}\oplus\mathcal{O}_{\ell}$.
\vspace{-0.2cm}
\begin{rem} \label{rem4.2.1}
    Note that a major challenge lies in how to reformulate the infinite-dimensional DR-CVaR constraint into a finite-dimensional constraint that is computationally tractable. To address this challenge, we employ DRO techniques to first express the DR-CVaR constraint as an optimization problem, then derive its dual problem, and finally transform it into Linear Matrix Inequality (LMI) constraints using the Schur complement.
\end{rem}

\subsection{COOL-DRMC with a novel update scheme}    \label{section4.3}
Based on the results in Sections~\ref{section4.1} and~\ref{section4.2}, separating hyperplanes can be computed to enforce collision-avoidance constraints and ensure the safety of the multi-robot system. However, the affine constraints induced by these hyperplanes may render the trajectory optimization problem~(\ref{e2.4.2.1}) infeasible. On the one hand, even when the committed convex regions are geometrically separated, the hyperplanes constructed by the adopted specific procedure may fail to separate them into non-intersecting half-spaces, causing previously committed trajectories to violate the resulting affine constraints. On the other hand, the ambiguity sets $\mathbb{D}_{\ell,t}$ are updated online, and the resulting changes in uncertainty may lead to newly computed affine constraints that are incompatible with existing committed trajectories. Therefore, a safe scheme for updating separating hyperplanes is required. Before introducing this scheme, we first define the terminal set $\mathcal{X}_{i}^{f}$ as follows.
\begin{defn}[Terminal set]  \label{defn4.3.1}
The terminal set $\mathcal{X}_{i}^{f}$ of robot $i$ is defined as a robust positively invariant set over $\mathcal{X}$. Specifically, for any $x_{i}(t) \in \mathcal{X}_{i}^{f}$, there exists a terminal control law $u_{i}(t) = \kappa_{i}^{f}(x_{i}(t)) \in \mathcal{U}$ entailing that $x_{i}(t+1) = f(x_{i}(t), u_{i}(t)) \in \mathcal{X}_{i}^{f}$. Moreover, for all $j \neq i$, $\big( h(\mathcal{X}_{i}^{f}) \oplus \mathcal{R}_{i} \big) \cap \big( h(\mathcal{X}_{j}^{f}) \oplus \mathcal{R}_{j} \big) = \emptyset$. And for all $\ell \in \mathbb{N}_{+}^{L}$, $\big( h(\mathcal{X}_{i}^{f}) \oplus \mathcal{R}_{i} \big) \cap \big( \mathbb{W}_{\ell,\varphi} \oplus \mathcal{O}_{\ell} \big) = \emptyset$.
\end{defn}
For robotic systems, such a robust terminal set can be defined by enforcing the velocity to be zero and the position to lie in the safe region, i.e., the complement of the union of the support sets of all obstacles. This safe region is nonempty in virtually all motion planning scenarios (at least in all case studies in Section \ref{section6}). Moreover, for linear robotic systems, less conservative terminal sets can be computed using tools such as Multi-Parametric Toolbox (MPT) \cite{herceg2013multi}, which allows for the construction of the maximal robust control invariant set inside the safe region, thereby reducing conservatism compared to a zero-velocity constraint. Readers can refer to Section IV.B of \cite{batkovic2022safe} for further explanations and practical examples. Next, we present the safe update schemes for the inter-robot separating hyperplanes and the robot–obstacle separating hyperplanes, respectively.

For the inter-robot case, at time $t$ we consider the update of the separating hyperplane $\mathcal{H}_{i,j,t,k}^{rob}$ between robots $i$ and $j$ at the predicted time step $k|t$. Robot $i$ checks whether the committed occupied regions $\hat{y}_{i}(k|t) \oplus \mathcal{R}_{i}$ and $\hat{y}_{j}(k|t) \oplus \mathcal{R}_{j}$ lie strictly on opposite sides of the newly computed separating hyperplane. If so, the newly computed hyperplane is adopted. Otherwise, for $k \in \mathbb{N}_{++}^{K-1}$, the separating hyperplane from the previous time step is retained, i.e., $\mathcal{H}_{i,j,t,k}^{rob} = \mathcal{H}_{i,j,t-1,k+1}^{rob}$. For $k = K$, since $\mathcal{H}_{i,j,t-1,K+1}^{rob}$ does not exist, a separating hyperplane is constructed based on the terminal occupied regions $h(\mathcal{X}_{i}^{f}) \oplus \mathcal{R}_{i}$ and $h(\mathcal{X}_{j}^{f}) \oplus \mathcal{R}_{j}$, which is guaranteed to exist by Definition~\ref{defn4.3.1}. Robot $j$ follows the same update rule, ensuring synchronous updates of the separating hyperplane. As a result, the inter-robot separating hyperplane remains feasible with respect to the committed trajectories.

For the robot-obstacle case, we consider the update of the separating hyperplane $\mathcal{H}_{i,\ell,t,k}^{obs}$ between robot $i$ and obstacle $\ell$ at the predicted time step $k|t$. For observable obstacles $\ell \in \mathbb{O}_{i,t}$, robot $i$ verifies whether the committed occupied region $\hat{y}_{i}(k|t) \oplus \mathcal{R}_{i}$ satisfies the affine constraint~(\ref{e2.4.2.1h}) induced by the newly computed separating hyperplane. If the constraint is satisfied, the new hyperplane is adopted. Otherwise, for $k \in \mathbb{N}_{++}^{K-1}$, the hyperplane from the previous time step is retained, i.e., $\mathcal{H}_{i,\ell,t,k}^{obs} = \mathcal{H}_{i,\ell,t-1,k+1}^{obs}$. For $k = K$, since $\mathcal{H}_{i,\ell,t-1,K+1}^{obs}$ does not exist, a separating hyperplane is constructed based on the terminal occupied region $h(\mathcal{X}_{i}^{f}) \oplus \mathcal{R}_{i}$ and the overall occupied region of obstacle $\ell$, $\mathbb{W}_{\ell,\varphi} \oplus \mathcal{O}_{\ell}$, whose existence is ensured by Definition~\ref{defn4.3.1}. For unobservable obstacles $\ell \notin \mathbb{O}_{i,t}$, the same update rule as in the case where the newly computed separating hyperplane violates the committed trajectory constraints is applied. Under these safe update rules, the robot-obstacle separating hyperplanes remains feasible with respect to the committed trajectories at all times.

\begin{algorithm}[htbp]
    \caption{COOL-DRMC for robot $i$ at time $t$}
     \label{algorithm1}
    \KwIn{$x_i(t)$, $\bm{\hat{y}_i(t)}$, $\{ N^{dat}_{i,\ell,t-1}, S^{tr}_{i,\ell,t-1} \} {\forall}\ell \in \mathbb{N}_{++}^{L}$}
    $\{ \varphi_{\ell}(t), \hat{\omega}_{\ell,t} \} {\forall}\ell \in \mathbb{O}_{i,t} \leftarrow$ Observation data \;
    \For{$j=1:N$ and $j\neq i$}{	
        $\{ \bm{\hat{y}_j(t)}, N^{dat}_{j,\ell,t-1} \} \leftarrow$ Robot $j$\;
        \For{$\ell \in \mathbb{O}_{i,t}$}{
            \If{$N^{dat}_{j,\ell,t-1} > N^{dat}_{i,\ell,t-1}$}{
                $S^{tr}_{j,\ell,t-1} \leftarrow$ Robot $j$\;
                $\{ N^{dat}_{i,\ell,t-1}, S^{tr}_{i,\ell,t-1} \} \leftarrow \{ N^{dat}_{j,\ell,t-1}, S^{tr}_{j,\ell,t-1} \}$\;
            }
        }
        \For{$k=1:K$}{
            $\mathcal{H}_{i,j,t,k}^{rob} \leftarrow$ Inter-robot separating hyperplane calculation (Section \ref{section4.1}) and check whether it should be updated (Section \ref{section4.3})\;
        }
    }
    \For{$\ell \in \mathbb{N}_{++}^{L}$}{
        \If{$\ell \in \mathbb{O}_{i,t}$}{
            $\{N^{dat}_{i,\ell,t}, S^{tr}_{i,\ell,t} \} \leftarrow$ Online variational inference for DPMM with $S^{tr}_{i,\ell,t-1}$ and $\hat{\omega}_{\ell,t}$  (\ref{e3.2.1})\;
            $\mathbb{D}_{\ell,t} \leftarrow$ Ambiguity set construction (\ref{e3.2.1})\;
            \For{$k=1:K$}{
                $\mathbb{D}_{\ell,t,k} \leftarrow$ Ambiguity set propagation (\ref{e3.2.4})\;
                \If{$m_{\ell,t,k} > M_{\ell,t,k}$}{
                    $\widetilde{\mathbb{D}}_{\ell,t,k} \leftarrow$Ambiguity set compression\;
                }
            }
        }
        \For{$k=1:K$}{
            $\mathcal{H}_{i,\ell,t,k}^{obs} \leftarrow$ Robot-obstacle separating hyperplane calculation (Section \ref{section4.2}) and check whether it should be updated (Section \ref{section4.3})\;
        }
    }
    $\bm{u^*_i(t)} \leftarrow$ Optimal decision of (\ref{e2.4.2.1})\;
    \KwOut{First step of optimal decision $u^*_i(0|t)$}
\end{algorithm}

With the above safe update rules in place, the COOL-DRMC algorithm for robot $i$ at time $t$ is delineated in Algorithm \ref{algorithm1}. The robot $i$ first communicates with other robots to execute COOL, and then calculates the separating hyperplanes $\mathcal{H}_{i,j,t,k}^{rob}$ and $\mathcal{H}_{i,\ell,t,k}^{obs}$, and finally solves the optimization problem (\ref{e2.4.2.1}). The first step of $\bm{u^*_i(t)}$ is utilized as control input $u^*_i(0|t)$. Additionally, we provide the performance analysis of COOL-DRMC. Specifically, we theoretically establish the probabilistic collision avoidance guarantee and the long-term tracking performance guarantee of COOL-DRMC. The probabilistic collision avoidance guarantee is as follows.
\begin{thm}   \label{thm4.4.1}
    Let \ref{e2.4.2.1}) be feasible at time $t$, then the control trajectory $y_{i}(t+1) = y_{i}(1|t)$ is collision avoidance with at least $\prod_{\ell=1}^{L} (\alpha_{\ell} \alpha^{u})$. Furthermore, let (\ref{e2.4.2.1}) be feasible at each time $t = 0,...,T-1$, then the control trajectory $y_{i}(1),...,y_{i}(T)$ is collision-free with at least $\prod_{\ell=1}^{L}(\alpha_{\ell} \alpha^{u})^{T}$.
\end{thm}
Its proof is provided in Appendix \ref{appendix E}. Theorem \ref{thm4.4.1} provides a specific collision avoidance probability for control trajectories. We note that $\alpha_{\ell}$ and $\alpha^{u}$ are freely adjustable. Therefore, for any given $\alpha\in [0, 1)$, we can adjust $\alpha_{\ell}$ and $\alpha^{u}$ such that the collision avoidance probability exceeds $\alpha$. It is worth noting that a given confidence level $\alpha$ can be decomposed into different combinations of $\alpha_{\ell}$ and $\alpha^u$. In practice, the choice of this decomposition can significantly affect the control performance. Empirically and theoretically, allocating $\alpha$ in a balanced manner, avoiding either $\alpha_{\ell}$ or $\alpha^u$ approaching 1, typically leads to improved performance. A detailed theoretical and numerical analysis of this sensitivity is presented in Section \ref{section5.6}. To obtain a long-term performance guarantee of COOL-DRMC, we make the following assumptions.
\begin{assum}   \label{assum4.4.1}
    The system model $f$ is continuous, and the stage cost $l_{i}$ and the terminal cost $q_{i}$ are continuous at the origin and satisfy $l_{i}(r_{i}^{x}(t),r_{i}^{u}(t)) = 0$ and $q_{i}(r_{i}^{x}(t)) = 0$. In addition, $l_{i}(x_{i}(t),u_{i}(t)) \geq g(\Vert x_{i}(t)-r_{i}^{x}(t) \Vert)$ for all feasible $x_{i}(t)$, $u_{i}(t)$, where $g$ is a $\mathcal{K}_{\infty}$-function.
\end{assum}
This is a standard assumption commonly used to prove the stability of the receding horizon controller. The typical choice for $l_{i}$ and $q_{i}$ in (\ref{e2.4.2.2}) satisfies this assumption.
\begin{assum}   \label{assum4.4.2}
    For any $x_{i}(t) \in \mathcal{X}_{i}^{f}$, under the terminal control law $u_{i}(t) = \kappa_{i}^{f}(x_{i}(t))$, the resulting next state $x_{i}(t+1) = f(x_{i}(t), u_{i}(t))$ satisfies$q_i(x_{i}(t+1)) - q_i(x_{i}(t)) \leq -l_i(x_{i}(t), u_{i}(t))$.
\end{assum}
Assumption \ref{assum4.4.2} essentially requires the terminal control policy to ensure a strict decrease in the cost-to-go function. While this assumption guarantees long-term tracking stability, it can be overly restrictive in motion planning tasks involving complex environments. This is because randomly moving obstacles may frequently intrude into the reference path, making long-term tracking infeasible in essence. Fortunately, this assumption is only necessary for ensuring the long-term tracking performance. If we are instead concerned with recursive feasibility, this condition can be safely relaxed. With the above two assumptions, we provide the long-term tracking performance guarantee of COOL-DRMC as follows, which is proven in Appendix \ref{appendix F}.
\begin{thm}   \label{thm4.4.2}
    Suppose Assumptions \ref{assum4.4.1} and \ref{assum4.4.2} hold and the optimization problem (\ref{e2.4.2.1}) is feasible at $t=0$, then the following long-term tracking performance guarantee is satisfied.
    \begin{equation}    \label{e4.4.1}
    \begin{split}
        \lim_{T\rightarrow\infty}\frac{1}{T}\sum\nolimits_{t=0}^{T}\Vert x_{i}(t)-r_{i}^{x}(t) \Vert = 0
    \end{split}
    \end{equation}
\end{thm}

\section{Simulation results}    \label{section5}
We conduct simulations to demonstrate the effectiveness and superiority of COOL-DRMC. We consider three scenarios: ($i$) two quadrotors navigating in a 3D environment with two dynamic obstacles whose motions are subject to multimodal distributions, ($ii$) two nonlinear vehicles navigating in an intersection-like scenario with two dynamic obstacles whose motions are subject to multimodal distributions, and ($iii$) three car-like robots navigating in a 2D environment with several static obstacles and dynamic obstacles whose motions are subject to time-varying distributions. Furthermore, we empirically discuss the impact of ambiguity set compression on control performance and computation time, the scalability of the proposed method in large-scale multi-robot teams, and the sensitivity of the proposed method to its key parameters. All the simulations are implemented in MATLAB R2023b. The YALMIP toolbox is utilized for mathematical modeling. The IPOPT 3.12.9 and MOSEK 10.1 solver are adopted to solve the optimization problem (\ref{e2.4.2.1}) and (\ref{e4.2.3}), respectively. The computational experiments are performed on a personal computer with 2.10 GHz Inter Core i7-13700 CPU and 32 GB RAM.

\subsection{Simulation 1 with multimodal distributions}     \label{section5.1}
\begin{figure}[bp]
\begin{center}
\includegraphics[width=8.4cm]{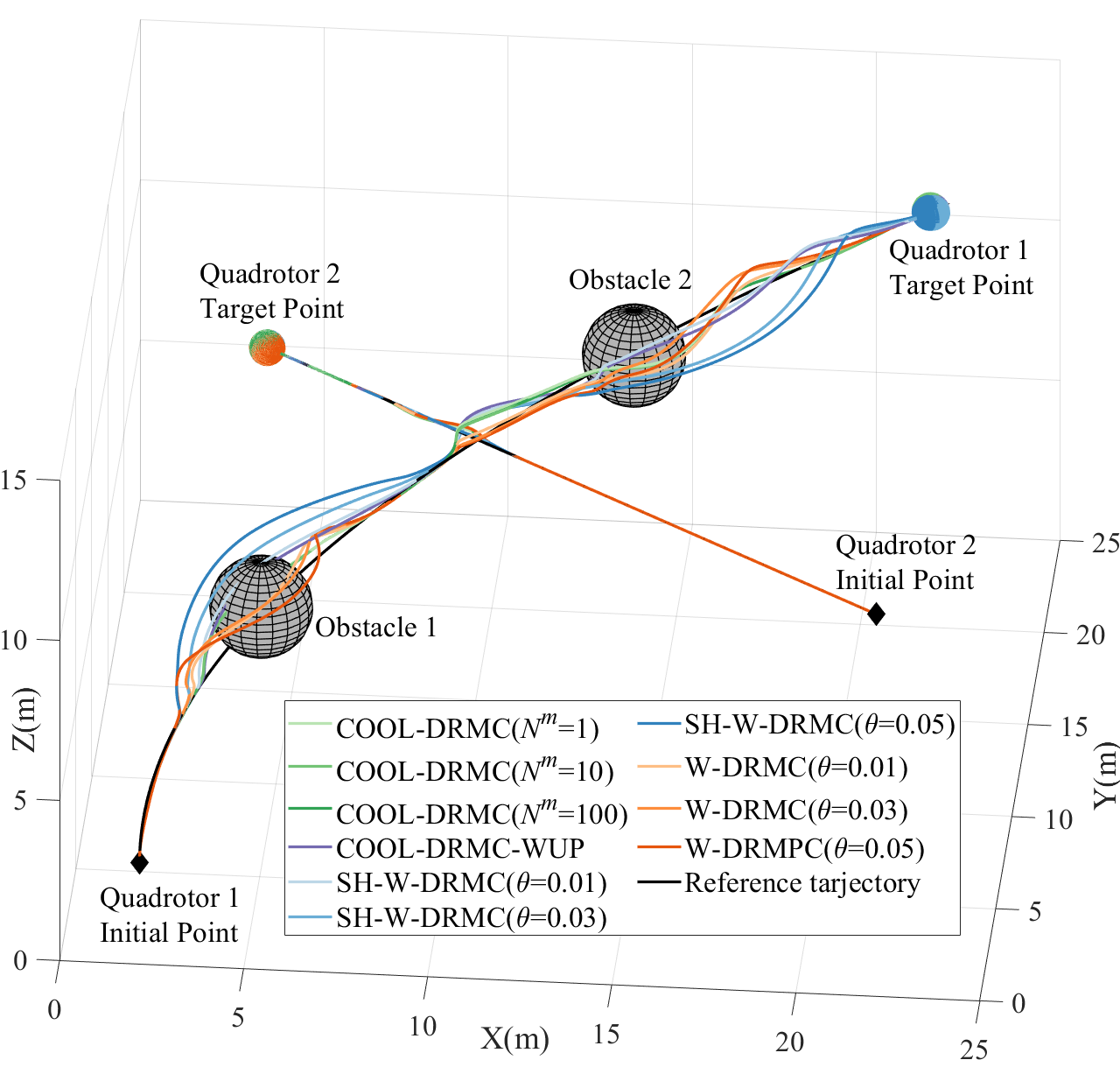}    
\caption{Trajectories of the quadrotors with different methods in simulation 1 (The positions of the obstacles and the quadrotors correspond to the positions at $T = 100s$)}.  
\label{fig5.1.1}                                 
\end{center}                                 
\end{figure}

Consider two quadrotors navigating in a 3D environment with the following dynamics \cite{mistler2001exact}.
\begin{equation} \label{e5.1.1}
\begin{split}
\begin{array}{lll}
     \ddot{x}=g\theta&\quad\ddot{y}=-g\phi&\quad\ddot{z}=\frac{1}{m_Q}u_1  \\
     \ddot{\phi}=\frac{l_Q}{I_{xx}}u_2&\quad\ddot{\theta}=\frac{l_Q}{I_{yy}}u_3&\quad\ddot{\psi}=\frac{l_Q}{I_{zz}}u_4 
\end{array}
\end{split}
\end{equation}
where $g=9.81kg/m^2$ is the gravitational acceleration, $m_Q=0.68kg$ is the mass of the quadrotor, $l_Q=0.23m$ is the distance between the CoM of the quadrotor and the rotor, and $I_{xx}=0.0075kg\cdot m^2$, $I_{yy}=0.0075kg\cdot m^2$ and $I_{zz}=0.013kg\cdot m^2$ represent the area moments of inertia about the principal axes in the body frame. The states are $\left( x, y, z, \dot{x}, \dot{y}, \dot{z}, \phi, \theta, \psi, \dot{\phi}, \dot{\theta}, \dot{\psi} \right) \in \mathbb{R}^{12}$, and the outputs are taken as the CoM of the quadrotor $\left( x,y,z \right)$. The objective is to control the quadrotors to track their reference trajectories while navigating around two randomly moving obstacles (grey) as shown in Fig. \ref{fig5.1.1}. The initial states of Quadrotor 1 and 2 are $x_1(0) = (2,2,2,0,0,0,0,0,0,0,0,0)^{T}$ and $x_2(0) = (22,2,11,0,0,0,0,0,0,0,0,0)^{T}$, respectively. The motions of the Obstacle 1 and 2 are sampled from the Gaussian mixture distributions with three components, denoted as $p_{1,t} = \sum_{i=1}^{3}\gamma_{1,t}^i\mathcal{N}(\mu_{1,t}^i,\Sigma_{1,t}^i)$ and $p_{2,t} = \sum_{i=1}^{3}\gamma_{2,t}^i\mathcal{N}(\mu_{2,t}^i,\Sigma_{2,t}^i)$ respectively. The parameters of the distribution are provided in Table \ref{tab5.1.a1}.
\begin{table}[h]
    \centering
    \caption{The parameters of motion distributions of Obstacle 1 and 2 in simulation 1 ($p_{1,t}$ and $p_{2,t}$).}
    \label{tab5.1.a1}
    \begin{tabular*}{\hsize}{@{}@{\extracolsep{\fill}}llccc@{}}
        \hline
        $\ell$ & $i$ & $\gamma $ & $\mu $ & $\Sigma $ \\
        \hline
        \multirow{3}{*}{$1$}    & $1$ & $0.3$ & $(0,0,0.1)^T$ & $diag(0.01,0.01,0.01)$ \\
                                & $2$ & $0.4$ & $(0,0,0.2)^T$ & $diag(0.01,0.01,0.01)$ \\
                                & $3$ & $0.3$ & $(0,0,0.3)^T$ & $diag(0.01,0.01,0.01)$ \\
        \hline
        \multirow{3}{*}{$2$}    & $1$ & $0.3$ & $(0,-0.2,0)^T$ & $diag(0.01,0.09,0.09)$ \\
                                & $2$ & $0.4$ & $(0.2,-0.2,0)^T$ & $diag(0.01,0.01,0.01)$ \\
                                & $3$ & $0.3$ & $(0.2,0,0)^T$ & $diag(0.01,0.01,0.01)$ \\
        \hline
    \end{tabular*}
\end{table}

The total simulation time is set to $T=100s$, and the system (\ref{e5.1.1}) is discretized using the sampling time $T_s=0.1s$. The prediction horizon is set to $K=10$. Increasing $K$ improves trajectory tracking accuracy but also increases computational complexity while decreasing $K$ can enhance the response speed of the controller. In our case study, $K=10$ is appropriate in our case study. The state and input weighting matrices are $Q=diag(1,1,1,0,...,0)$ and $R=0.01diag(1,1,1,1)$, respectively. We set $\alpha=0.95$. The relative increase in $Q$ and $R$ will respectively cause the controller to emphasize tracking accuracy and the smoothness of the control trajectory. In our case study, we choose a relatively large $Q$ to achieve a more precise control trajectory. The robotic motion control methods analyzed in this subsection are listed below.
\begin{itemize}
\item W-DRMC: Wasserstein Distribution Robust Motion Control with radius $\theta=0.01,0.03,0.05$ \cite{hakobyan2021wasserstein}.
\item SH-W-DRMC: Safe Halfspace based W-DRMC with radius $\theta=0.01,0.03,0.05$ \cite{safaoui2023distributionally}.
\item COOL-DRMC: The proposed method with the maximum number of basic ambiguity set $M_{\ell,t,k}=1,10,100$.
\item COOL-DRMC-WUP: The proposed method Without Uncertainty Propagation. The separating hyperplanes $\mathcal{H}_{\ell,t,k}^{obs}$ over the entire prediction horizon are computed based on the ambiguity set of $\varphi_{\ell}(1|t)$.
\end{itemize}
The trajectories generated by different motion control methods are shown in Fig. \ref{fig5.1.1}. We observe that, compared with other methods, COOL-DRMC exhibits a closer proximity to the reference trajectory when maneuvering around obstacles.
\begin{table}[t]
    \centering
    \caption{Average cost, tracking error, and computation time of quadrotor 1 with different methods.}
    \label{tab5.1.1}
    \resizebox{8.4cm}{2.3cm}{
    \begin{tabular}{llccc}
        \hline
        \multicolumn{2}{l}{Method} & Average Cost & Tracking Error & Computation Time  \\
        \hline
        \multirow{3}{*}{W-DRMC}    & $\theta=0.01$ & 8.6658 & 0.619 & 3.0535 \\
                                    & $\theta=0.03$ & 12.1944 & 0.712 & 3.4039 \\
                                    & $\theta=0.05$ & 17.1208 & 0.858 & 3.6547 \\
        \hline
        \multirow{3}{*}{SH-W-DRMC} & $\theta=0.01$ & 27.6748 & 0.801 & 0.0069 \\
                                    & $\theta=0.03$ & 44.2133 & 1.063 & 0.0070 \\
                                    & $\theta=0.05$ & 64.6716 & 1.313 & 0.0071 \\
        \hline
        \multicolumn{2}{l}{COOL-DRMC-WUP} & 22.0459 & 0.700 & 0.0328  \\
        \hline
        \multirow{3}{*}{COOL-DRMC} & $M_{\ell,t,k} = 1$ & 7.6468 & 0.517 & 0.0601 \\
                                   & $M_{\ell,t,k} = 10$ & 7.5875 & 0.515 & 0.1286 \\
                                   & $M_{\ell,t,k} = 100$ & 7.4860 & 0.513 & 0.8267 \\
        \hline
    \end{tabular}}
\end{table}

Table \ref{tab5.1.1} shows the average cost, tracking error, and computation time of different methods, where all results are averaged over 500 randomized trials to ensure statistical robustness. It can be observed that the average cost of the proposed method is significantly smaller than those of all other methods. Similarly, the proposed method also achieves a lower tracking error than competing methods, demonstrating improved trajectory tracking performance. The COOL-DRMC method with $M_{\ell,t,k}=10$ reduces the cost and tracking error by an average of 12.44\% and 35.7\% compared with the W-DRMC with $\theta=0.01$ by using the fined-grained distribution information in ambiguity sets. Furthermore, compared with W-DRMC, COOL-DRMC reduces the average computation time by one to two orders of magnitude depending on the value of $M_{\ell,t,k}$ through the utilization of separating hyperplanes. Note that while COOL-DRMC-WUP and SH-W-DRMC exhibit shorter computation time compared to COOL-DRMC, the quadrotors controlled by COOL-DRMC-WUP and SH-W-DRMC are confined within the same halfspace throughout the entire prediction horizon due to the absence of uncertainty propagation, which results in control lag and a cost increase of over 188\%.

Additionally, based on the ambiguity set compression, COOL-DRMC can balance the control performance and computation time by adjusting $M_{\ell,t,k}$. Specifically, when $M_{\ell,t,k}$ is large ($M_{\ell,t,k}=100$), the ambiguity set is scarcely compressed, thus boosting control performance but incurring a substantial computational burden. Conversely, as $M_{\ell,t,k}$ decreases ($M_{\ell,t,k}=1$), the ambiguity set is gradually compressed, leading to a drastic reduction in computational burden at the expense of slightly diminished control performance.

\subsection{Simulation 2 using nonlinear dynamics}         \label{section5.4}
Consider two vehicles navigating in a 2D environment with the following nonlinear dynamics \cite{pepy2006path}.
\begin{equation}    \label{e5.4.1}
\begin{aligned}
    \left[\begin{array}{c}
    p_{x,i}(t+1) \\
    p_{y,i}(t+1) \\
    \theta_{i}(t+1) \\
    v_{i}(t+1)
\end{array}\right] = \left[\begin{array}{c}
    p_{x,i}(t) + T_s \cdot v_{i}(t) \cos\theta_{i}(t) \\
    p_{y,i}(t) + T_s \cdot v_{i}(t) \sin\theta_{i}(t) \\
    \theta_{i}(t) + T_s \cdot \tfrac{v_{i}(t)}{l} \tan\phi_{i}(t) \\
    v_{i}(t) + T_s \cdot a_{i}(t)
\end{array}\right]
\end{aligned}
\end{equation}
\begin{table}[t]
\centering
\caption{The parameters of motion distributions of Obstacle 1 and 2 in simulation 2 ($p_{1,t}$ and $p_{2,t}$).}
\label{tab5.4.1}
\begin{tabular*}{\hsize}{@{}@{\extracolsep{\fill}}llccc@{}}
    \hline
    $\ell$ & $i$ & $\gamma $ & $\mu $ & $\Sigma $ \\
    \hline
    \multirow{3}{*}{$1$}    & $1$ & $0.25$ & $(0,0.4)^T$ & $diag(0.05^2,0.05^2)$ \\
                            & $2$ & $0.5$  & $(-0.2,-0.2)^T$ & $diag(0.05^2,0.05^2)$ \\
                            & $3$ & $0.25$ & $(0.4,0)^T$ & $diag(0.05^2,0.05^2)$ \\
    \hline
    \multirow{3}{*}{$2$}    & $1$ & $0.25$ & $(0,-0.4)^T$ & $diag(0.05^2,0.05^2)$ \\
                            & $2$ & $0.5$  & $(-0.2,0.2)^T$ & $diag(0.05^2,0.05^2)$ \\
                            & $3$ & $0.25$ & $(0.4,0)^T$ & $diag(0.05^2,0.05^2)$ \\
    \hline
\end{tabular*}
\end{table}
\hspace{-4pt}where $x_{i}(t) = (p_{x,i}(t), p_{y,i}(t), \theta_{i}(t), v_{i}(t))$ are the position, orientation, and velocity of the vehicle, respectively. $l := 0.25m$ is the length, and $T_s = 0.1$ is the sampling time. The system inputs are the steering angle $\phi_{i}(t) \in [-\pi/6,\pi/6]$ and the acceleration $a_{i}(t) \in [-5,5]$. The control horizon is set to $K=10$. The objective is to control the vehicles to track their reference trajectories through an intersection-like scenario while navigating around two randomly moving obstacles. The initial states of Vehicles 1 and 2 are $x_1(0) = (6,0.5,\pi/2,0)^{T}$ and $x_2(0) = (2,4,0,0)^{T}$, respectively. The motions of the Obstacles 1 and 2 are sampled from the Gaussian mixture distributions with three components, denoted as $p_{1,t} = \sum_{i=1}^{3}\gamma_{1,t}^i\mathcal{N}(\mu_{1,t}^i,\Sigma_{1,t}^i)$ and $p_{2,t} = \sum_{i=1}^{3}\gamma_{2,t}^i\mathcal{N}(\mu_{2,t}^i,\Sigma_{2,t}^i)$ respectively. The parameters of the distribution are provided in Table \ref{tab5.4.1}. The total simulation time is set to $T=100s$.
\begin{figure}[b]
\begin{center}
    \includegraphics[width=8.4cm]{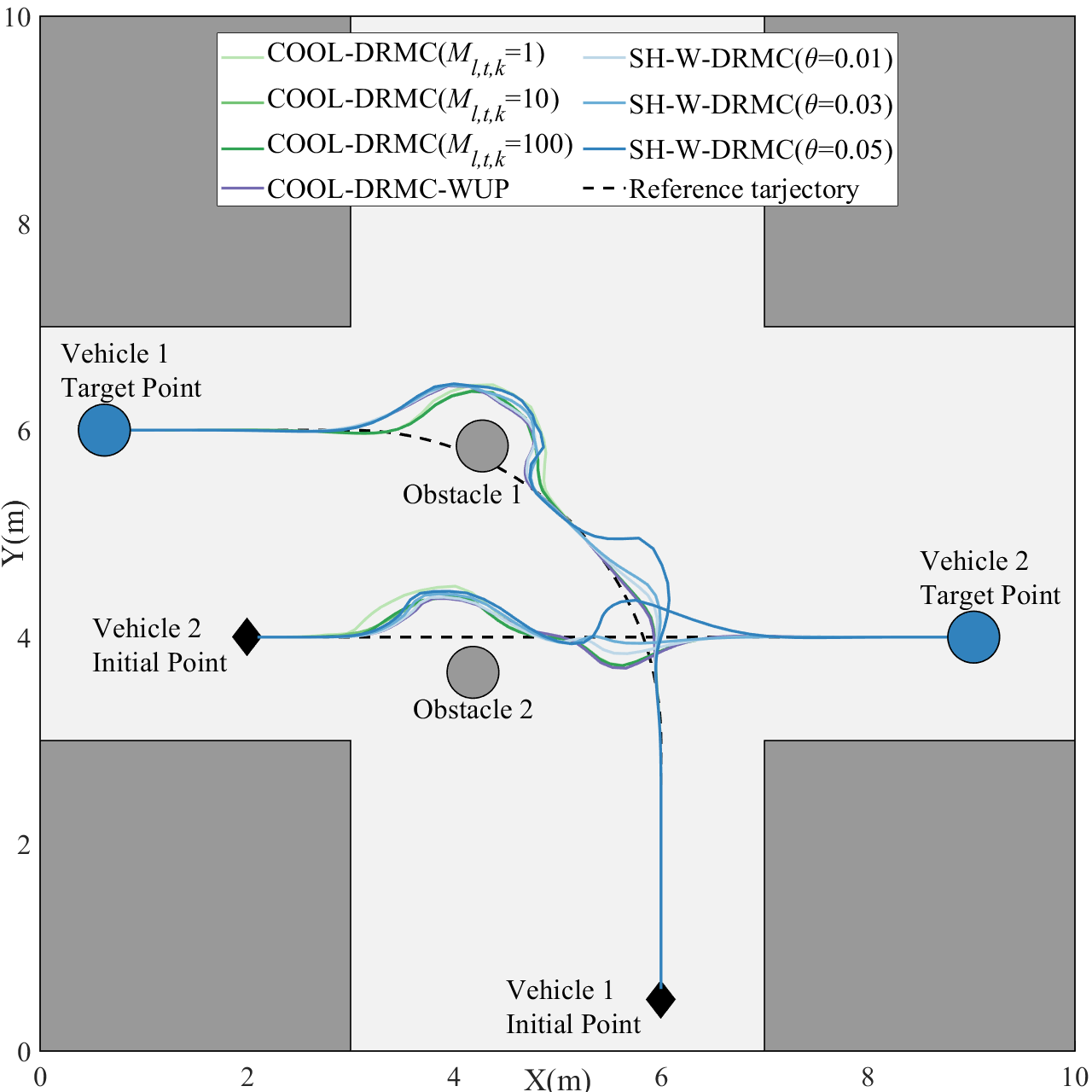}    
    \caption{Trajectories of the vehicles with different methods (The positions of the obstacles and the vehicles correspond to the positions at $T = 100s$).}  
    \label{fig5.4.1}                                 
\end{center}                                
\end{figure}
The state and input weighting matrices are $Q=diag(10,10,0,0)$ and $R=0.01diag(1,0.1)$, respectively. We set $\alpha=0.95$. The robotic motion control methods analyzed in this subsection, except for W-DRMC, are the same as those in simulation 1.
 
The trajectories generated by different motion control methods are shown in Fig. \ref{fig5.4.1}. The average cost and computation time of different methods are shown in Table \ref{tab5.4.2}, where all results are averaged over 500 randomized trials to ensure statistical robustness. It can be observed that COOL-DRMC is still capable of avoiding collisions with moving obstacles even for nonlinear vehicle dynamics. Additionally, the COOL-DRMC method with $M_{\ell,t,k}=10$ reduces the cost and tracking error by an average of 63.66\% and 31.45\% compared with the W-DRMC with $\theta=0.001$. Similar to the results of Simulation (i), COOL-DRMC can flexibly balance the trade-off between average cost (tracking error) and computation time by adjusting $M_{\ell,t,k}$. It is important to note that in this simulation, we do not compare the proposed method with W-DRMC. This is because W-DRMC does not employ the separating hyperplane technique to convexify the distributionally robust collision avoidance constraint, which, in combination with the nonlinear model, imposes a significant computational burden on the control optimization problem. Specifically, the average computation time of W-DRMC exceeds $100$ seconds at each sampling time, rendering it entirely unsuitable for real-time robotic control. In summary, the simulation results demonstrate that the proposed method can ensure safety and achieve excellent control performance even for nonlinear models.
\begin{table}[t]
    \centering
    \caption{Average cost, tracking error, and computation time of vehicle team with different methods.}
    \label{tab5.4.2}
    \resizebox{8.4cm}{1.7cm}{
    \begin{tabular}{llccc}
        \hline
        \multicolumn{2}{l}{Method} & Average Cost & Tracking Error & Computation Time(s)  \\
        \hline
        \multirow{3}{*}{SH-W-DRMC} & $\theta=0.001$ & 10.5347 & 0.248 & 0.0875 \\
                                    & $\theta=0.003$ & 10.5064 & 0.259 & 0.0730 \\
                                    & $\theta=0.005$ & 13.5858 & 0.331 & 0.0758 \\
        \hline
        \multicolumn{2}{l}{COOL-DRMC-WUP} & 10.0661 & 0.237 & 0.1164  \\
        \hline
        \multirow{3}{*}{COOL-DRMC} & $M_{\ell,t,k} = 1$ & 5.6015 & 0.197 & 0.1667 \\
                                   & $M_{\ell,t,k} = 10$ & 3.8281 & 0.170 & 0.3493 \\
                                   & $M_{\ell,t,k} = 100$ & 3.6395 & 0.169 & 1.4477 \\
        \hline
    \end{tabular}}
\end{table}

\subsection{Simulation 3 with time-varying distributions}   \label{section5.2}
\begin{figure*}[b]
    \centering
    \subfloat[$t=0$]{
        \label{fig5.2.1a}
        \includegraphics[width=8cm]{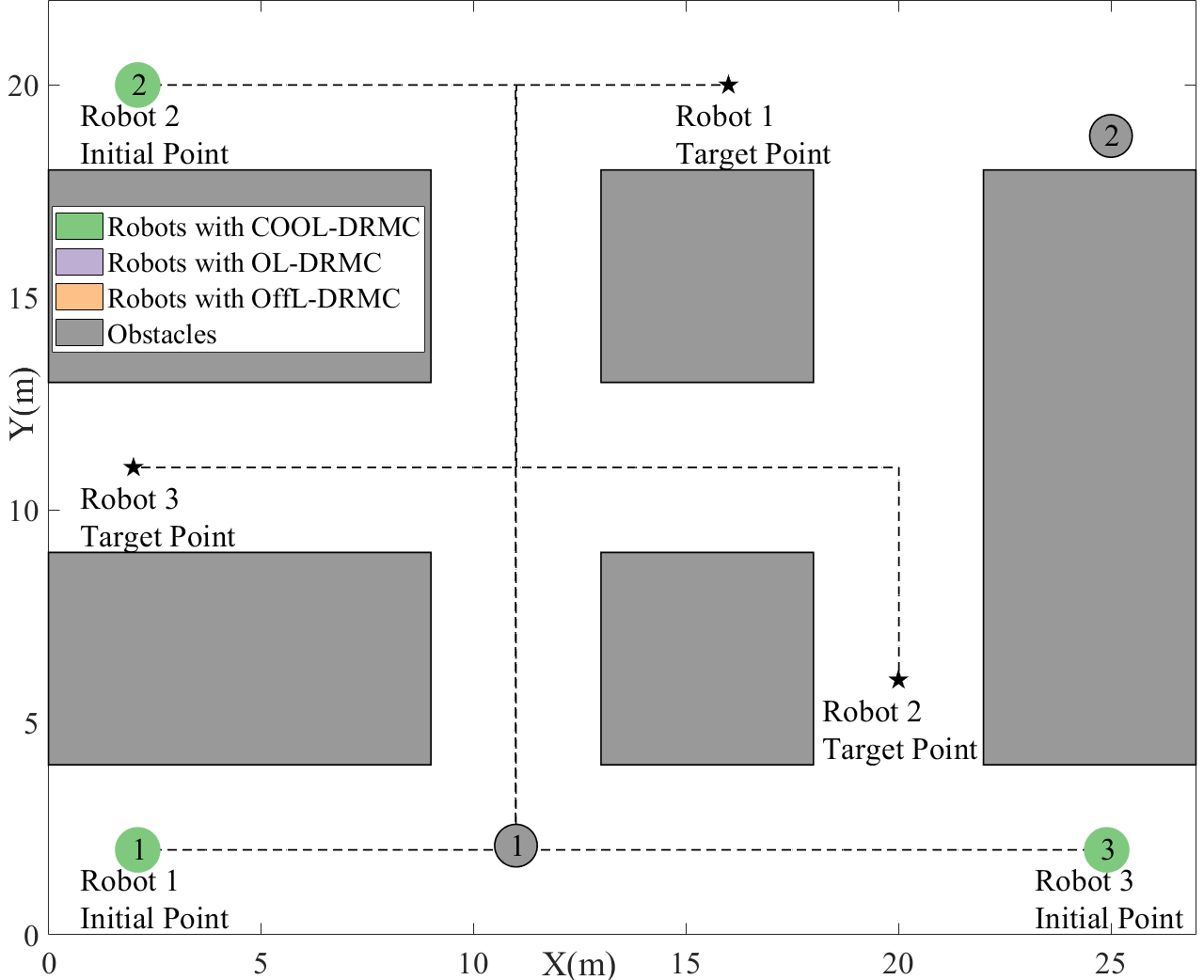}
    }
    \subfloat[$t=88$]{
        \label{fig5.2.1b}
        \includegraphics[width=8cm]{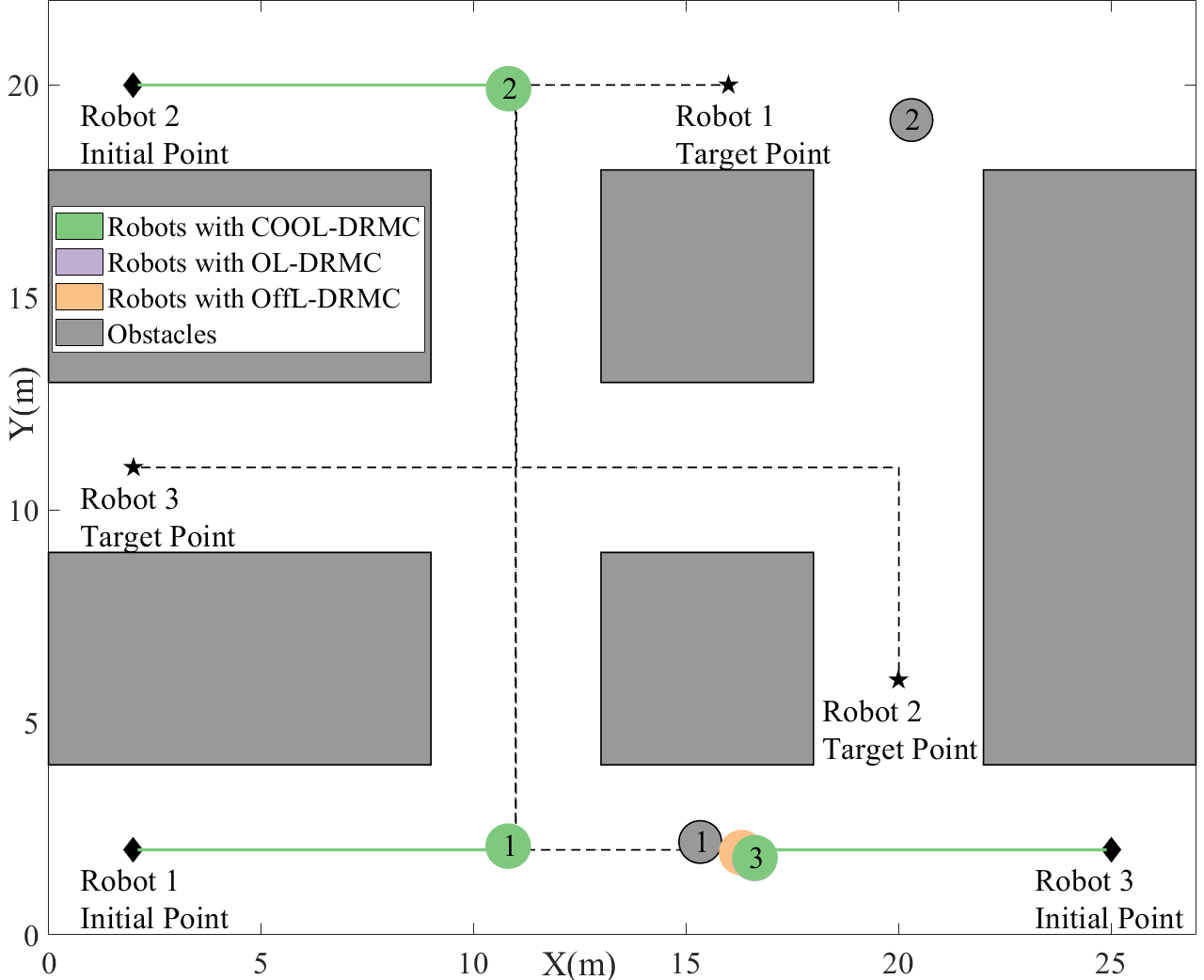}
    }\\
    \subfloat[$t=259$]{
        \label{fig5.2.1c}
        \includegraphics[width=8cm]{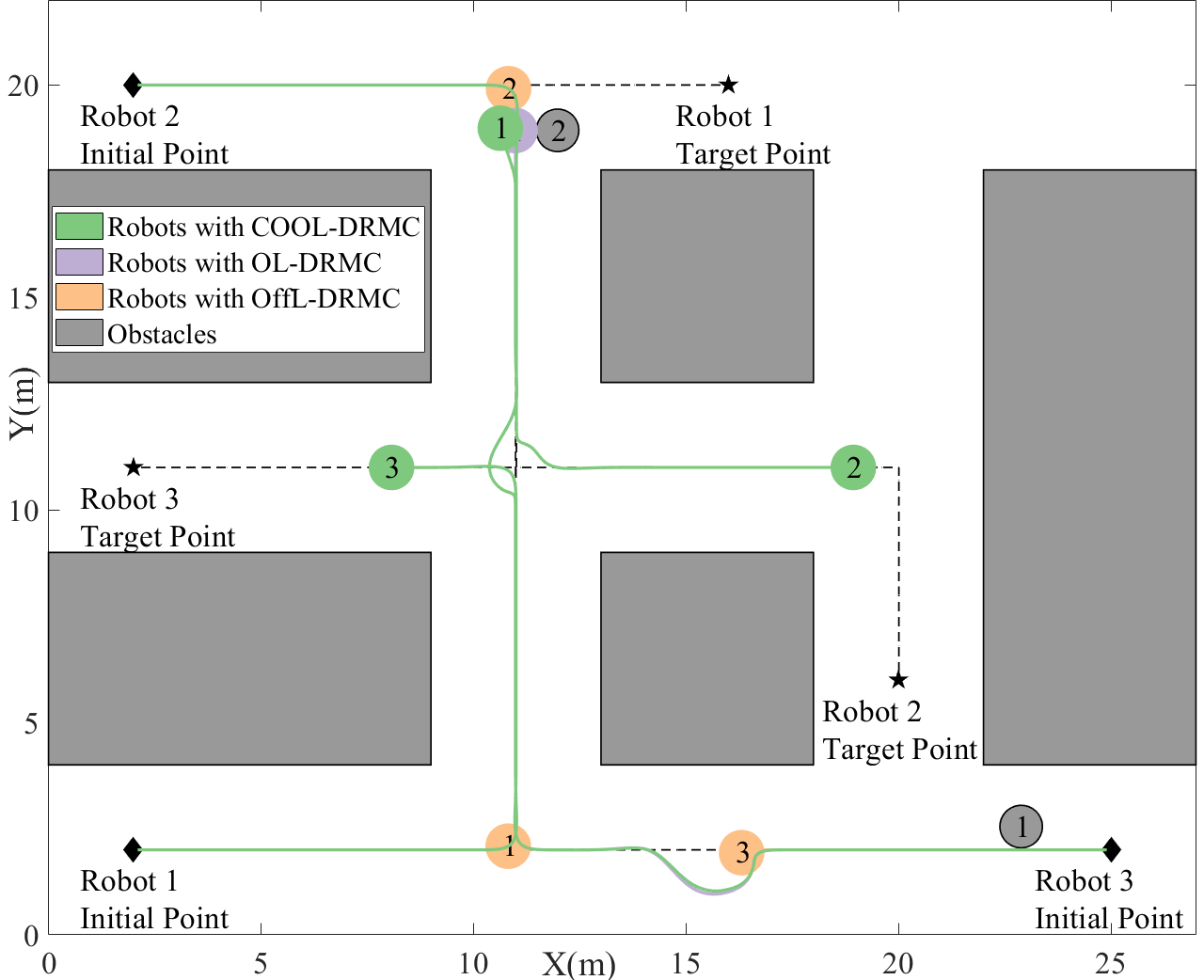}
    }
    \subfloat[$t=320$]{
        \label{fig5.2.1d}
        \includegraphics[width=8cm]{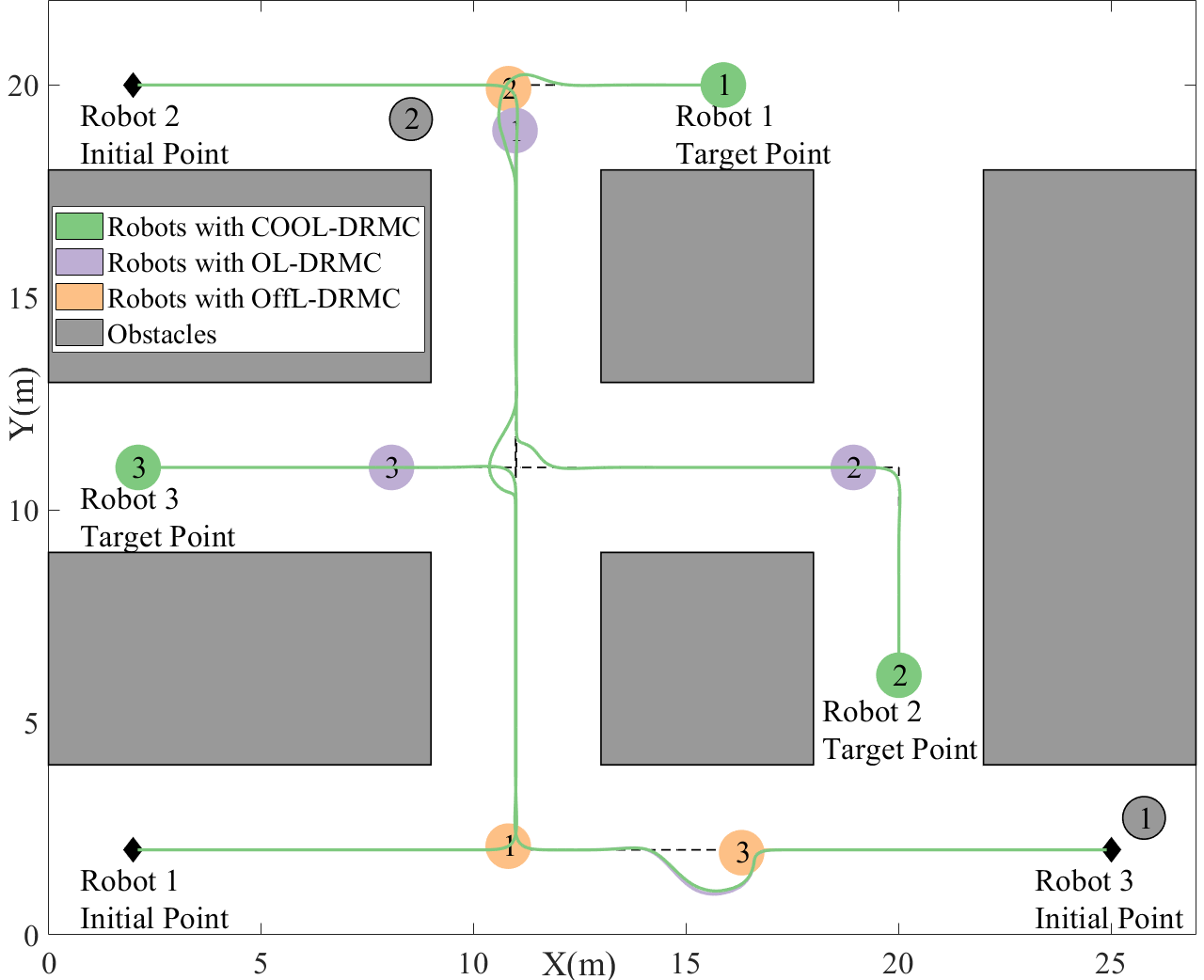}
    }
    \caption{Trajectories of robots with different control methods in simulation 2. (Rectangles represent static obstacles and circles represent the robots and obstacles. Numbers on the circles denote the indices of robots or obstacles and the different colors of the circles indicate the robots controlled by different methods as shown in the legend.)}  
    \label{fig5.2.1}                             
\end{figure*}
Consider three car-like robots navigating in a 2D environment with the following double integrator dynamics.
\begin{equation} \label{e5.2.1}
\begin{split}
x(t+1)=\left[\begin{array}{cccc}
     1&0&T_s&0  \\
     0&1&0&T_s  \\
     0&0&1&0  \\
     0&0&0&1
\end{array}\right]x(t)+\left[\begin{array}{cc}
    \frac{T_s^2}{2} & 0 \\
    0 & \frac{T_s^2}{2} \\
    T_s & 0 \\
    0 & T_s
\end{array}\right]u(t)
\end{split}
\end{equation}
where $\left( p^x,p^y,v^x,v^y \right)$ is the state of a robot, consisting of its CoM and velocity vector. The input $u=\left( a^x, a^y \right)$ is chosen as the acceleration vector. Similarly, the sampling time $T_s$ is selected as $0.1s$.

The simulation configuration is illustrated in Fig. \ref{fig5.2.1a}. The initial states of Robot 1, 2 and 3 are $x_1(0) = (2,2,0,0)^{T}$, $x_2(0) = (2,20,0,0)^{T}$ and $x_3(0) = (25,2,0,0)^{T}$, respectively. The circular robots with radius $r_r=0.5$ aim to track the given reference trajectories in a cluttered 2D environment with some static and dynamic obstacles where the static obstacles are modeled as rectangles representing topographic barriers and the dynamic obstacles are approximated by the circles with radius $r_o=0.5$ representing human or other robots that cannot be communicated with. The motions in each direction of the Obstacle 1 and 2 are sampled from time-varying Gaussian distributions. Specifically, the means of the motion distribution of Obstacle 1 and 2 are $(0.5,0)^T$ and $(-0.52,0)^T$ respectively, and the standard deviation of all distributions used to generate historical data is $0.01$, while the standard deviation increases to $0.3$ for the real-time data generation during the runtime of robotic motion control. The control input for robots is limited to lie in $\mathcal{U}:=\left\{ u\in \mathbb{R}^2 \big| \Vert u \Vert_{\infty} \leq 4 \right\}$ and its state is restricted to $\mathcal{X}:=\left\{ x\in \mathbb{R}^4 \big| (0,0,-2,-2) \leq x \leq (22,27,2,2) \right\}$. The reference trajectories of each robot are generated through the upper-level planner. The prediction horizon is set to $K=10$ and the state and input weighting matrices are chosen as $Q=diag(1,1,0,0)$ and $R=0.1diag(1,1,1,1)$. The maximum number of base ambiguity sets is selected as $M_{\ell,t,k}=10$ and we set $\alpha=0.95$. The simulation is independently repeated 500 times. The following robotic motion control methods are analyzed in this subsection to demonstrate the superiority of the proposed collaborative online learning.
\begin{itemize}
\item OffL-DRMC: Offline DRO-based Motion Control. Without online learning, the ambiguity sets are constructed only based on historical data.
\item OL-DRMC: Online-Learning-enabled DRO-based Motion Control. Without collaborative learning, each robot updates the ambiguity sets solely based on its own observed data.
\item COOL-DRMC: The proposed method. The ambiguity sets are updated through the collaborative online learning method.
\end{itemize}
Fig. \ref{fig5.2.1} shows the simulation results from one of the 500 independent simulations. For OffL-DRMC, the estimated distributions of random motions are based exclusively on historical data. Unfortunately, the motions of the obstacles become violent because of the increased standard deviation and the deviation between estimated and real-time distributions leads to the collision between obstacle 1 and robot 3 controlled by the OffL-DRMC at time $t=88$ as shown in Fig. \ref{fig5.2.1b}. In contrast, COOL-DRMC and OL-DRMC can learn the changes in the distribution and keep a safe distance from the obstacle.

For obstacle 2, it cannot be observed by robot 1 until it is approached due to the occlusion of static obstacles. For the OL-DRMC, only a few observed data of obstacle 2 are available when the controller attempts to avoid obstacle 2. The lack of data leads to the inability of robot 1 to learn the changed distribution in time, which causes the collision between obstacle 2 and robot 1 at time $t=259$ as shown in Fig. \ref{fig5.2.1c}. Fortunately, before robot 1 approaches obstacle 2, robot 2 has learned the changed distribution using sufficient observed data streams thanks to its advantageous position. Utilizing the collaborative online learning, robot 1 with COOL-DRMC learns the distribution of obstacle 2 based on the learning structure received from robot 2 instead of relying on the initial learning structure when robot 1 initially observes obstacle 2. Therefore, COOL-DRMC controls robot 1 to avoid obstacle 2 safely. Finally, as shown in Fig. \ref{fig5.2.1d}, the robot team controlled by COOL-DRMC successfully reached their target points, outperforming the other two methods.

\begin{table}[b]
    \centering
    \caption{Collision-free rate of different methods with 500 simulations.}
    \label{tab5.2.1}
    \resizebox{8.4cm}{1.25cm}{
    \begin{tabular}{lccc}
        \hline
        \multirow{2}{*}{Method} & \multicolumn{3}{c}{Collision-free rate}\\
        \cmidrule(r){2-4}
                                & With obstacle 1 & With obstacle 2 & Total\\
        \hline
        OffL-DRMC & $64.6\%$ & $69.2\%$ & $46.4\%$ \\
        \hline
        OL-DRMC & $99.4\%$ & $69.2\%$ & $68.8\%$ \\
        \hline
        COOL-DRMC & $99.6\%$ & $98\%$ & $97.6\%$ \\
        \hline
    \end{tabular}}
\end{table}
Table \ref{tab5.2.1} shows the collision-free rate of different methods with 500 independent simulations. The collision-free rate is computed as the ratio of successful trials without collision over 500 randomized experiments. Without online learning, OffL-DRMC has a $35.4\%$ rate of collision with obstacle 1. In contrast, COOL-DRMC and OL-DRMC demonstrate a collision-free rate of over $99\%$ in navigating around obstacle 1. Although OL-DRMC avoids collisions with obstacle 1 with a high probability, robot 1 controlled by OL-DRMC has a $30.8\%$ rate of collision with obstacle 2 due to the absence of real-time data, as described above. Therefore, the collision-free rate throughout the total navigation of OL-DRMC is only $68.8\%$. Thanks to the proposed collaborative online learning, the collision-free rate throughout the total navigation of COOL-DRMC is $97.6\%$.

\subsection{Discussions on ambiguity set compression}   \label{section5.3}
In this section, we discuss the influence of the maximum number of base ambiguity sets $M_{\ell,t,k}$ on the control performance and computation time. We consider a robot with the dynamics as (\ref{e5.2.1}) circumvents a dynamic obstacle whose motions in each direction are sampled from a Gaussian mixture distribution with three components.

\begin{figure}[htbp]
\begin{center}
\includegraphics[width=8.4cm]{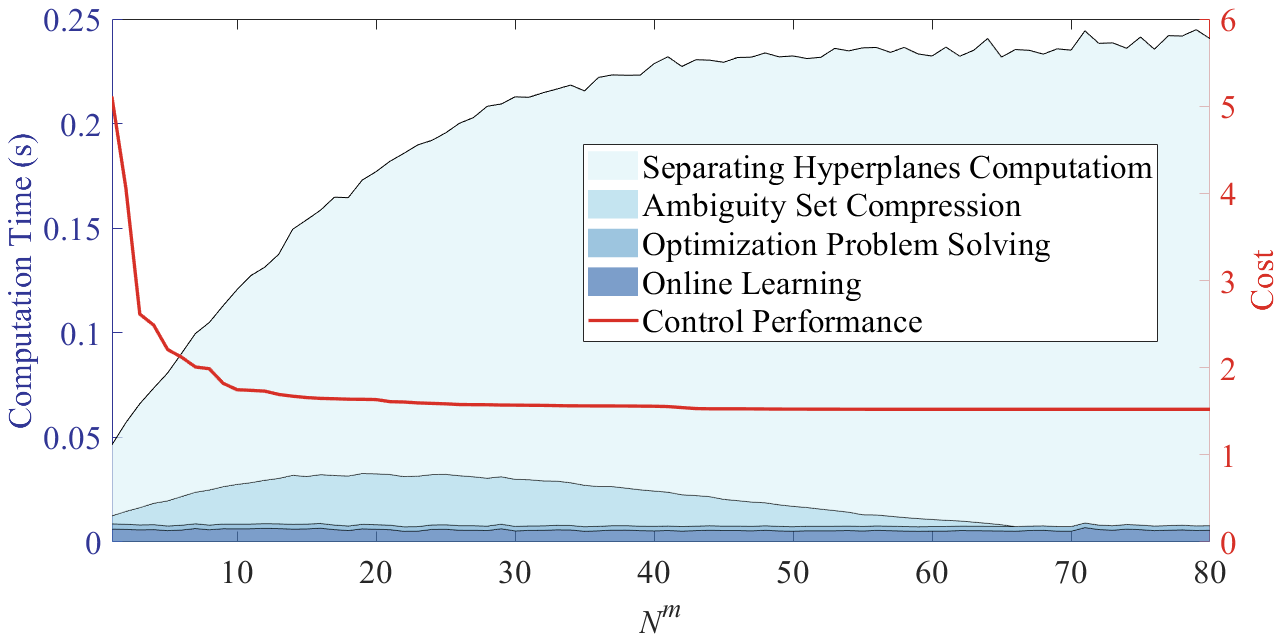}    
\caption{The average cost and the average computation time with $M_{\ell,t,k} = 1 \sim 80$.}  
\label{fig5.3.1}                                 
\end{center}                                 
\end{figure}

In Fig. \ref{fig5.3.1}, the average cost and computation time with $M_{\ell,t,k} = 1 \sim 80$ are shown. The average cost decreases rapidly with an increase of $M_{\ell,t,k}$ when $M_{\ell,t,k}<10$ and remains stable as $M_{\ell,t,k}>10$. Instead, the average computation time consistently ascends as $M_{\ell,t,k}$ increases until $M_{\ell,t,k}>66$, especially when $M_{\ell,t,k}<40$. Therefore, it is a wise choice to choose $M_{\ell,t,k}=10$ to balance the average cost and computation time. Additionally, the user can choose other values of $M_{\ell,t,k}$ according to the computational capability, the real-time requirement, and control performance.

The computation time is mostly spent on calculating the robot-obstacle separating hyperplanes over the entire prediction horizon, and therefore $M_{\ell,t,k}$ determines the maximum number of constraints in the optimal problem (11) and in turn determines the computation time. For the simulation configuration with prediction horizon $K=10$ and the number of components of distribution $m_{\ell,t}=3$, the maximum number of base ambiguity sets for the non-compressed ambiguity set is $m_{\ell,t,K}=\binom{m_{\ell,t}-1}{K+m_{\ell,t}-1}=66$. Thus the ambiguity set compression time reduces to $0$ and the computation time is no longer increased when $M_{\ell,t,k}>66$. Note that if the hardware computational capability is sufficient and better control performance is desired, $M_{\ell,t,k}$ can be increased. Conversely, to improve computational speed, $M_{\ell,t,k}$ can be decreased. Empirically, we show that $M_{\ell,t,k}=10$ provides a reasonable compromise in most cases.

\subsection{Scalability in large-scale multi-robot teams}   \label{section5.5}
To evaluate the scalability of the proposed COOL-DRMC framework in large-scale multi-robot systems, we conduct a comparative study between COOL-DRMC and a centralized baseline. As a distributed method, COOL-DRMC allows each robot to solve its own local optimization problem, which significantly alleviates the computational burden as the number of robots increases. This makes it particularly suitable for large-scale systems where centralized methods typically become intractable due to the curse of dimensionality. In this subsection, we consider a 3D navigation scenario where $n$ quadrotors are tasked with avoiding collisions with each other and two randomly moving dynamic obstacles. The dynamic obstacles follow multimodal stochastic trajectories as described in Scenario (i). We vary $n$ from 2 to 24 to evaluate the computation scalability of both methods. The centralized baseline solves a single joint optimization problem for all robots, including collision avoidance constraints among robots and constraints between robots and obstacles. In contrast, the decentralized COOL-DRMC solves separate optimization problems for each robot and uses a consensus mechanism based on committed trajectories to coordinate inter-robot collision avoidance via separating hyperplanes.
\begin{figure}[htbp]
\begin{center}
\includegraphics[width=8.4cm]{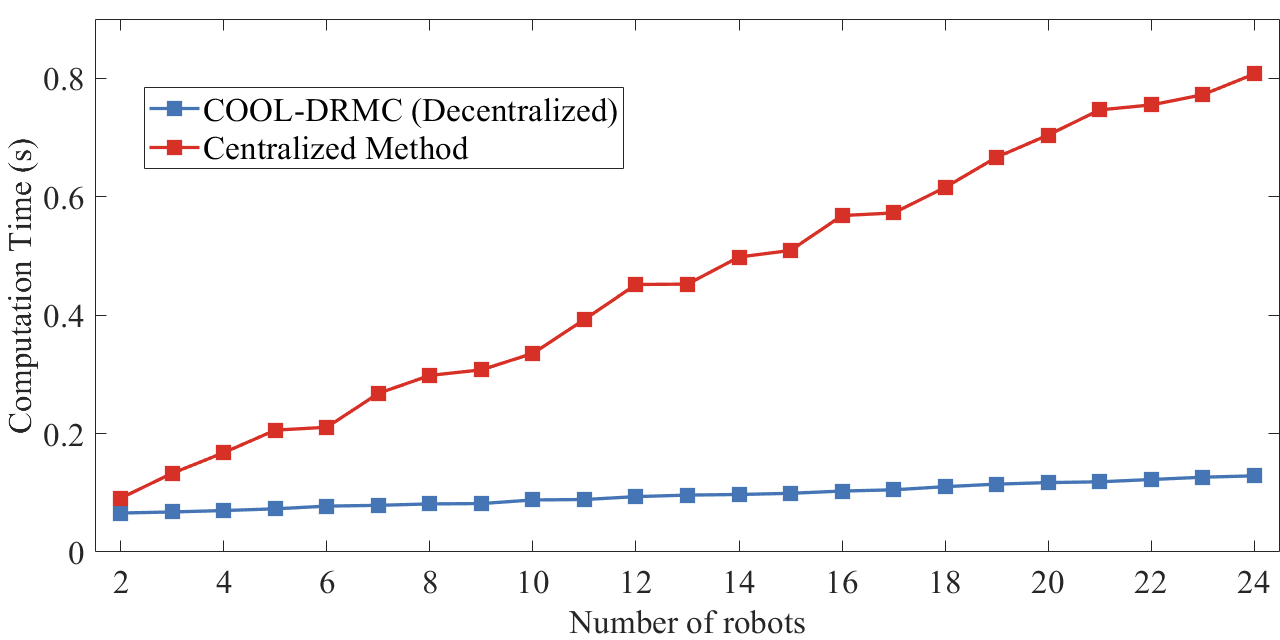}    
\caption{Average computation with different robot numbers using decentralized and centralized methods.}  
\label{fig5.5.1}                                 
\end{center}                                 
\end{figure}

Fig.\ref{fig5.5.1} compares the average computation time of the COOL-DRMC (decentralized) method with a centralized method, as the number of robots increases from 2 to 24. For COOL-DRMC, when additional robots are added, the extra burden for each robot arises from calculating the inter-robot separating hyperplanes, which relies on committed trajectories and involves lightweight geometric operations, resulting in a negligible increase in computation time. Therefore, the average computation time of COOL-DRMC grows only marginally with the team size. In contrast, the centralized method controls all robots jointly and must evaluate all robot-obstacle separating hyperplanes simultaneously. This results in a computational cost that grows almost proportionally with the number of robots, as shown in Fig. \ref{fig5.5.1}. These results demonstrate the scalability advantage of the proposed decentralized framework for large-scale multi-robot motion planning.

\subsection{Sensitivity to the key parameters}   \label{section5.6}
Given a fixed robustness level $\alpha$, different allocations to $\alpha_{\ell}$ (ambiguity set coverage) and $\alpha^{u}$ (DR-CVaR constraint confidence) can significantly impact control performance. In this subsection, we analyze the impact of different combinations of $\alpha_{\ell}$ and $\alpha^u$ on control performance from both theoretical and empirical perspectives.

From a qualitative analysis perspective, Theorem \ref{thm4.4.1} shows that the probability of robot $i$ avoiding collision with all obstacles $\ell=1,...,L$ at time $t$ is given by $\alpha = \prod_{\ell=1}^L \alpha_{\ell} (\alpha^u)^L$, indicating a trade-off between $\alpha_{\ell}$ and $\alpha^u$. First, according to Theorem 1, increasing $\alpha_{\ell}$ enlarges the ambiguity set $\mathbb{D}_{\ell,t,k}$ by amplifying its associated moment deviation parameters. As shown in (\ref{e3.2.2}), this effect becomes particularly severe as $\alpha_{\ell}$ approaches 1, leading to substantial conservativeness and degraded control performance. Similarly, increasing $\alpha^u$ tightens the distributionally robust constraint. For moment-based ambiguity sets, this may lead to overly pessimistic decisions, effectively degenerating to worst-case support-based constraints when $\alpha^u$ approaches 1.

To empirically discuss the sensitivity of different allocations of $\alpha$, we conduct a numerical case study using a double integrator robot tasked with avoiding a dynamic obstacle. To facilitate a wider range of $\alpha$ allocations and allow for a more comprehensive set of combinations, we fix the overall confidence level at $\alpha=0.8$. We then record the robot’s average cost and tracking error under various combinations of $\alpha_{u}$ and $\alpha_{\ell}$, as summarized in Table \ref{tab5.6.1}.
\begin{table}[htbp]
    \centering
    \caption{Average cost and tracking error under different combinations of $\alpha_{u}$ and $\alpha_{\ell}$.}
    \label{tab5.6.1}
    \resizebox{8.4cm}{1.2cm}{
    \begin{tabular}{lrcccccc}
        \specialrule{1.2pt}{1pt}{1pt}
            $\alpha_{u}$    & 0.810 & 0.820 & 0.850 & 0.895 & 0.942 & 0.976 & 0.988 \\
        \specialrule{0.7pt}{1pt}{1pt}
            $\alpha_{\ell}$ & 0.988 & 0.976 & 0.942 & 0.895 & 0.850 & 0.820 & 0.810 \\
        \specialrule{1.2pt}{1pt}{1pt}
            Average cost    & 11.97 & 8.60 & 6.77 & 6.19 & 7.57 & 12.98 & 34.22 \\
        \specialrule{0.7pt}{1pt}{1pt}
            Tracking error  & 0.307 & 0.270 & 0.233 & 0.222 & 0.250 & 0.349 & 0.661 \\
        \specialrule{1.2pt}{1pt}{1pt}
    \end{tabular}}
\end{table}
As shown in Table \ref{tab5.6.1}, the results are consistent with our theoretical analysis. Specifically, allocating the overall confidence level $\alpha$ more evenly between $\alpha_{u}$ and $\alpha_{\ell}$, so as to avoid either parameter approaching 1, achieves the lowest average cost and tracking error. In contrast, when either $\alpha_{u}$ or $\alpha_{\ell}$ becomes too close to 1, both the average cost and tracking error increase rapidly. These results confirm that balancing the confidence levels is crucial for achieving strong performance without excessive conservativeness.

\vspace{-0.2cm}
\section{Conclusion}    \label{section6}
\vspace{-0.2cm}
We proposed a COOL-DRMC framework for decentralized multi-robot systems to ensure safe navigation in dynamic and random environments. A novel COOL method based on the DPMM was introduced to integrate decentralized data streams and efficiently extract distribution information. An ambiguity set for motion was constructed by leveraging the local-moment information, and an ambiguity set propagation method was proposed to obtain the ambiguity set for positions over the entire prediction horizon. Additionally, we presented a compression scheme with safety guarantee for adjusting the complexity of this ambiguity set. The provably safe ambiguity set compression admitted the balance between control performance and computation time. The spatial allocation protocol ensures collision avoidance among robots. The robot-obstacle collision avoidance constraints were reformulated by deriving separating hyperplanes through tractable SDP problems. Finally, we provided the probabilistic collision avoidance guarantee and the long-term tracking performance guarantee for the proposed framework. The simulation results demonstrated the efficient learning capability, less conservative performance, lower computational overhead, flexibility, and safety of the proposed method.

\vspace{-0.2cm}
\begin{ack}                               
\vspace{-0.2cm}
This work was supported in part by the National Natural Science Foundation of China under Grants 62473256 and 62103264. The authors would like to thank the associate editor and all the reviewers for their valuable and constructive comments, which help improve this paper.  
\end{ack}

\appendix
\vspace{-0.2cm}
\section{Proof of Theorem \ref{thm3.2.1}}          \label{appendix A}
\vspace{-0.2cm}
\begin{pf}
We know that the underlying distribution $p_{\ell,t}^{*}$ has $m_{\ell,t}$ components, which is written as follows.
\begin{equation}    \label{eA1}
\begin{split}
    p_{\ell,t}^{*} = \sum\nolimits_{i=1}^{m_{\ell,t}} \gamma_{\ell,t}^{i*} p_{\ell,t}^{i*}
\end{split}
\end{equation}
Note that DPMM can correctly identify the component from which each learning data point is sampled. Therefore for $N $ learning data sampled from $p_{\ell,t}^{i*}$, we can have
\begin{equation}    \label{eA2}
\begin{split}
    \mathbb{P}\{ p_{\ell,t}^{i*} \in \mathbb{D}  \} \geq \alpha_{\ell}^{i}
\end{split}
\end{equation}
where the inequality derives from \textbf{Corollary 4 in \cite{delage2010distributionally}}. In addition, according to the \textbf{Bretagnolle-Huber-Carol} concentration inequality \cite{bellet2015robustness}, the underlying weight $\bm{\gamma_{\ell,t}^{*}}$ lies in set $\Delta_{\ell,t}$ parametrized by the selected value of $\theta_{\ell,t}$, with at least probability $\chi_{\ell}$.
\begin{equation}    \label{eA3}
\begin{split}
    \mathbb{P}\{ \bm{\gamma_{\ell,t}^{*}} \in \Delta_{\ell,t} \} \geq \chi_{\ell}
\end{split}
\end{equation}
where $\chi_{\ell}$ denotes the confidence level. Following the independence of random events, we have (\ref{e3.2.3}). This completes the proof. $\quad\square$
\end{pf}

\vspace{-0.2cm}
\section{Proof of Lemma \ref{lem3.2.1}}          \label{appendix B}
\vspace{-0.2cm}
\begin{pf}
For any $p_{\ell,t}$ lies in the ambiguity set $\mathbb{D}_{\ell,t}$ in (\ref{e3.2.1}), it can be written as follows.
\begin{equation}    \label{eB1}
\begin{split}
    p_{\ell,t} = \sum\nolimits_{i=1}^{m_{\ell,t}} \gamma_{\ell,t}^{i}  p_{\ell,t}^{i} 
\end{split}
\end{equation}
where $\bm{\gamma_{\ell,t}} \in \Delta_{\ell,t}$ and $p_{\ell,t}^{i} \in \mathbb{D}_{\ell,t}^{i} $. We note that $p_{\varphi_{\ell}(k|t)} = \delta_{\varphi_{\ell}(t)} * [p_{\ell,t}]^{k}$, which can be reformulated below.
\begin{alignat}{2}     
    &p_{\varphi_{\ell}(k|t)}    \label{eB2}\\
    \overset{(a)}{=}& \delta_{\varphi_{\ell}(t)} * \left[ \gamma_{\ell,t}^{1} p_{\ell,t}^{1} + ... + \gamma_{\ell,t}^{m_{\ell,t}} p_{\ell,t}^{m_{\ell,t}} \right]^{k} \notag\\
    \overset{(b)}{=}& \delta_{\varphi_{\ell}(t)}* \sum_{ \mathop{\sum}\limits_{i=1}^{m_{\ell,t}}k_i=k } \left\{ \frac{k!}{k_1!...k_{m_{\ell,t}}!} \mathop{\prod}\limits_{i=1}^{m_{\ell,t}}(\gamma_{\ell,t}^i)^{k_i} \mathop{\Xi}\limits_{i=1}^{m_{\ell,t}}\left[ p  \right]^{k_i} \right\} \notag\\
    \overset{(c)}{=}& \delta_{\varphi_{\ell}(t)}* \sum_{j=1}^{m_{\ell,t,k}} \left\{ \frac{k!}{k_1^j!...k_{m_{\ell,t}}^j!} \mathop{\prod}\limits_{i=1}^{m_{\ell,t}}(\gamma_{\ell,t}^i)^{k_i^j} \mathop{\Xi}\limits_{i=1}^{m_{\ell,t}}\left[ p  \right]^{k_i^j} \right\} \notag\\
    \overset{(d)}{=}& \sum\nolimits_{j=1}^{m_{\ell,t,k}} \gamma_{\ell,t,k}^{j} p_{\ell,t,k}^{j}  \notag
\end{alignat}
where $*$ denotes the convolution, $[\cdot]^{k}$ represents the convolution multiplier, and $\Xi$ denotes successive convolutions. Equation (b) is derived from the multinomial theorem \cite{kataria2016probabilistic}, where the summation symbol denotes all the non-negative integer-solution groups that satisfy $\sum_{i=1}^{m_{\ell,t}} k_i = k$. In equation (c), $m_{\ell,t,k}$ denotes the number of the groups calculated in (\ref{e3.2.5a}) and $k_{1}^{j}...k_{m_{\ell,t}}^{j}$ is the $j$-th group. The $\gamma_{\ell,t,k}^{j}$ and $p_{\ell,t,k}^{j}$ are as follows.
\begin{equation}        \label{eB3}
\begin{split}
    &\gamma_{\ell,t,k}^{j} = \frac{k!}{k_1^j!...k_{m_{\ell,t}}^j!} \left( (\gamma_{\ell,t}^{1})^{k_1^j} ... (\gamma_{\ell,t}^{m_{\ell,t}})^{k_{\ell,t}^j} \right) \\
    &p_{\ell,t,k}^{j} = \delta_{\varphi_{\ell}(t)}* [p_{\ell,t}^{1}]^{k_1^j} *...* [p_{\ell,t}^{m_{\ell,t}}]^{k_{\ell,t}^j} 
\end{split}
\end{equation}
It remains to prove that $\bm{\gamma_{\ell,t,k}}=(\gamma_{\ell,t,k}^{1},...,\gamma_{\ell,t,k}^{m_{\ell,t,k}}) \in \Delta_{\ell,t,k}$ and $p_{\ell,t,k}^{j}  \in \mathbb{D}_{\ell,t,k}^{j}$ in (\ref{e3.2.4}).

(1) Prove that $\bm{\gamma_{\ell,t,k}} \in \Delta_{\ell,t,k}$.\\
Based on the multinomial theorem, we know
\begin{equation}        \label{eB4}
\begin{split}
    \sum\nolimits_{j=1}^{m_{\ell,t,k}} \gamma_{\ell,t,k}^{j} = \left( \sum\nolimits_{i=1}^{m_{\ell,t}} \gamma  \right)^{k} = 1
\end{split}
\end{equation}
and thus the first condition in $\Delta_{\ell,t,k}$ is naturally satisfied, while the key lies in the second condition as follows.
\begin{alignat}{2}
   &\sum\nolimits_{j=1}^{m_{\ell,t,k}} \left| \gamma_{\ell,t,k}^{j}-\hat{\gamma}_{\ell,t,k}^{j} \right| \label{eB5}\\
   =& \sum\nolimits_{j=1}^{m_{\ell,t,k}} \frac{k!}{k_1^j!...k_{m_{\ell,t}}^j!} \left| \prod\nolimits_{i=1}^{m_{\ell,t}}(\gamma_{\ell,t}^i)^{k_i^j} - \prod\nolimits_{i=1}^{m_{\ell,t}}(\hat{\gamma}_{\ell,t}^i)^{k_i^j} \right| \notag
\end{alignat}
For the absolute value term in (\ref{eB5}), we have the following reformulation.
\begin{alignat}{2}
   &\left| \prod\nolimits_{i=1}^{m_{\ell,t}}(\gamma_{\ell,t}^i)^{k_i^j} - \prod\nolimits_{i=1}^{m_{\ell,t}}(\hat{\gamma}_{\ell,t}^i)^{k_i^j} \right| \label{eB6}\\
   =& \left| \sum_{l=1}^{m_{\ell,t}} \left( \prod_{i=1}^{l-1} (\hat{\gamma} )^{k_i^j} \left( (\gamma_{\ell,t}^{l})^{k_l^j} - (\hat{\gamma}_{\ell,t}^{l})^{k_l^j} \right) \prod_{i=l+1}^{m_{\ell,t}} (\gamma )^{k_i^j} \right) \right| \notag\\
   \overset{(a)}{\leq}& \sum_{l=1}^{m_{\ell,t}} \left| (\gamma_{\ell,t}^{l})^{k_l^j} - (\hat{\gamma}_{\ell,t}^{l})^{k_l^j} \right| \prod_{i=1}^{l-1} (\hat{\gamma} )^{k_i^j} \prod_{i=l+1}^{m_{\ell,t}} (\gamma )^{k_i^j} \notag\\
   \overset{(b)}{\leq}& \sum_{l=1}^{m_{\ell,t}} \left| (\gamma_{\ell,t}^{l})^{k_l^j} - (\hat{\gamma}_{\ell,t}^{l})^{k_l^j} \right| \left( \frac{\sum_{i=1}^{l-1}\hat{\gamma}  + \sum_{i=l+1}^{m_{\ell,t}}\gamma }{m_{\ell,t}-1} \right)^{k-k_l^j} \notag\\
   \overset{(c)}{\leq}& \sum_{l=1}^{m_{\ell,t}} \left| (\gamma_{\ell,t}^{l})^{k_l^j} - (\hat{\gamma}_{\ell,t}^{l})^{k_l^j} \right| \left( \frac{1-\hat{\gamma}_{\ell,t}^{l}+\theta_{\ell,t}}{m_{\ell,t}-1} \right)^{k-k_l^j}   \notag
\end{alignat}
where inequality (a) follows from the triangle inequality, inequality (b) results from the weighted arithmetic-geometric mean inequality, and inequality (c) is derived from the fact $\bm{\gamma_{\ell,t}} \in \Delta_{\ell,t}$ and the triangle inequality. Next, we further bound the absolute value term in (\ref{eB6}).
\begin{alignat}{2}
   &\left| (\gamma_{\ell,t}^{l})^{k_l^j} - (\hat{\gamma}_{\ell,t}^{l})^{k_l^j} \right| \label{eB7}\\
   =& \left| \gamma_{\ell,t}^{l} - \hat{\gamma}_{\ell,t}^{l} \right| \left( (\gamma_{\ell,t}^{l})^{k_l^j-1} + (\gamma_{\ell,t}^{l})^{k_l^j-2}(\hat{\gamma}_{\ell,t}^{l}) + (\hat{\gamma}_{\ell,t}^{l})^{k_l^j-1} \right)  \notag\\
   \leq& \left| \gamma_{\ell,t}^{l} - \hat{\gamma}_{\ell,t}^{l} \right| \cdot k_l^j (\hat{\gamma}_{\ell,t}^{l} + \theta_{\ell,t})^{k_l^j-1} \notag
\end{alignat}
where the inequality results from the fact $\bm{\gamma_{\ell,t}} \in \Delta_{\ell,t}$. Substituting (\ref{eB6}) and (\ref{eB7}) into (\ref{eB5}), we can obtain
\begin{alignat}{2}
   &\sum\nolimits_{j=1}^{m_{\ell,t,k}} \left| \gamma_{\ell,t,k}^{j}-\hat{\gamma}_{\ell,t,k}^{j} \right|      \label{eB8}\\
   \overset{(a)}{\leq}& \sum_{j=1}^{m_{\ell,t,k}} \frac{k!}{k_1^j!...k_{m_{\ell,t}}^j!} \sum_{l=1}^{m_{\ell,t}} \left\{\begin{array}{l}
        \left| \gamma_{\ell,t}^{l} - \hat{\gamma}_{\ell,t}^{l} \right| k_l^j (\hat{\gamma}_{\ell,t}^{l} + \theta_{\ell,t})^{k_l^j-1}  \\
        \left( \frac{1-\hat{\gamma}_{\ell,t}^{l}+\theta_{\ell,t}}{m_{\ell,t}-1} \right)^{k-k_l^j} 
   \end{array}\right\}          \notag\\
   \overset{(b)}{=}& \sum_{l=1}^{m_{\ell,t}} \left| \gamma_{\ell,t}^{l} - \hat{\gamma}_{\ell,t}^{l} \right| \sum_{j=1}^{m_{\ell,t,k}} \left\{\begin{array}{l}
        \frac{k!}{k_1^j!...k_{m_{\ell,t}}^j!} k_l^j (\hat{\gamma}_{\ell,t}^{l} + \theta_{\ell,t})^{k_l^j-1}  \\
        \left( \frac{1-\hat{\gamma}_{\ell,t}^{l}+\theta_{\ell,t}}{m_{\ell,t}-1} \right)^{k-k_l^j} 
   \end{array}\right\}   \notag\\
   \overset{(c)}{=}& \sum_{l=1}^{m_{\ell,t}} \left| \gamma_{\ell,t}^{l} - \hat{\gamma}_{\ell,t}^{l} \right| k \left( \hat{\gamma}_{\ell,t}^{l} + \theta_{\ell,t} + (m_{\ell,t}-1) \frac{1-\hat{\gamma}_{\ell,t}^{l}+\theta_{\ell,t}}{m_{\ell,t}-1} \right)^{k-1} \notag\\
   \overset{(d)}{=}& k\theta_{\ell,t}(1+2\theta_{\ell,t})^{k-1}    \notag
\end{alignat}
where equation (c) follows from the properties of multinomial expansion, where $\hat{\gamma}_{\ell,t}^{l}+\theta_{\ell,t}$ is treated as the $l$-th term of the $m_{\ell,t}$-term polynomial in equation (c), and expanding it based on multinomial theorem can recover equation (b). This finishes the proof that $\bm{\gamma_{\ell,t,k}} \in \Delta_{\ell,t,k}$.

(2) Prove that $p_{\ell,t,k}^{j}  \in \mathbb{D}_{\ell,t,k}^{j}$.\\
For ease of exposition, we use $\mu_{\ell,t}^{i}, \mu_{\ell,t,k}^{j}$ and $\Sigma_{\ell,t}^{i}, \Sigma_{\ell,t,k}^{j}$ represent the means and covariances of $p $ and $p_{\ell,t,k}^{j}$, respectively. As a result, the two conditions in $\mathbb{D}_{\ell,t}^{i}$ (\ref{e3.2.1}) can be stated as follows.
\begin{alignat}{2}
   & (\mu_{\ell,t}^{i} - \hat{\mu}_{\ell,t}^{i})^{T} (\hat{\Sigma}_{\ell,t}^{i})^{-1} (\mu_{\ell,t}^{i} - \hat{\mu}_{\ell,t}^{i}) \leq \beta_{\ell,t}^{i}  \label{eB9}\\
   & \Sigma_{\ell,t}^{i} + (\mu_{\ell,t}^{i} - \hat{\mu}_{\ell,t}^{i}) (\mu_{\ell,t}^{i} - \hat{\mu}_{\ell,t}^{i})^{T} \preceq \varepsilon_{\ell,t}^{i} \hat{\Sigma}_{\ell,t}^{i} \label{eB10}
\end{alignat}
Note that since $\hat{\Sigma}_{\ell,t}^{i}$ is positive definite, equation (\ref{eB9}) has the following three equivalent forms.
\begin{subequations}    \label{eB11}
\begin{alignat}{2}
   & (\mu_{\ell,t}^{i} - \hat{\mu}_{\ell,t}^{i})^{T} (\hat{\Sigma}_{\ell,t}^{i})^{-1} (\mu_{\ell,t}^{i} - \hat{\mu}_{\ell,t}^{i}) \leq \beta_{\ell,t}^{i}  \label{eB11a}\\
   \Longleftrightarrow& \left(\begin{array}{cc}
       \hat{\Sigma}_{\ell,t}^{i} & (\mu_{\ell,t}^{i} - \hat{\mu}_{\ell,t}^{i}) \\
       (\mu_{\ell,t}^{i} - \hat{\mu}_{\ell,t}^{i})^{T} & \beta_{\ell,t}^{i}
   \end{array}\right) \succeq 0   \label{eB11b}\\
   \Longleftrightarrow& (\mu_{\ell,t}^{i} - \hat{\mu}_{\ell,t}^{i}) (\mu_{\ell,t}^{i} - \hat{\mu}_{\ell,t}^{i})^{T} \preceq \beta_{\ell,t}^{i} \hat{\Sigma}_{\ell,t}^{i}    \label{eB11c}
\end{alignat}
\end{subequations}
Furthermore, based on (\ref{eB11c}), the sufficient and necessary conditions for (\ref{eB10}) can be derived as follows.
\begin{subequations}    \label{eB12}
\begin{alignat}{2}
   &\Sigma_{\ell,t}^{i} \preceq (\varepsilon_{\ell,t}^{i} - \beta_{\ell,t}^{i}) \hat{\Sigma}_{\ell,t}^{i}  \label{eB12a}\\
   \Longrightarrow  &\Sigma_{\ell,t}^{i} + (\mu_{\ell,t}^{i} - \hat{\mu}_{\ell,t}^{i}) (\mu_{\ell,t}^{i} - \hat{\mu}_{\ell,t}^{i})^{T} \preceq \varepsilon_{\ell,t}^{i} \hat{\Sigma}_{\ell,t}^{i}   \label{eB12b} \\
   \Longrightarrow  &\Sigma_{\ell,t}^{i} \preceq \varepsilon_{\ell,t}^{i} \hat{\Sigma}_{\ell,t}^{i}   \label{eB12c}
\end{alignat}
\end{subequations}
It is evident that (\ref{eB12a}) and (\ref{eB12c}) are, respectively, the sufficient and necessary conditions for (\ref{eB12b}). Clearly, the two conditions in $\mathbb{D}_{\ell,t,k}^{j}$ (\ref{e3.2.4}) have the corresponding forms above, which we will not elaborate further. We are now ready to begin the proof. Based on the relationship between $p_{\ell,t,k}^{j}$ and $p $ in (\ref{eB3}), we can obtain
\begin{equation}        \label{eB13}
\begin{split}
    \mu_{\ell,t,k}^{j} = \varphi_{\ell}(t) + \sum_{i=1}^{m_{\ell,t}} k_i^j \mu ,\  \Sigma_{\ell,t,k}^{j} = \sum_{i=1}^{m_{\ell,t}} k_i^j \Sigma 
\end{split}
\end{equation}
For the first condition in $\mathbb{D}_{\ell,t,k}^{j}$, we verify its form corresponding to (\ref{eB11b}) as follows.
\begin{alignat}{2}
   & \left(\begin{array}{cc}
       \hat{\Sigma}_{\ell,t,k}^{j} & (\mu_{\ell,t,k}^{j} - \hat{\mu}_{\ell,t,k}^{j}) \\
       (\mu_{\ell,t,k}^{j} - \hat{\mu}_{\ell,t,k}^{i})^{T} & \beta_{\ell,t,k}^{j}
   \end{array}\right)   \label{eB14}\\
   =& \sum\nolimits_{i=1}^{m_{\ell,t}} k_i^j \left(\begin{array}{cc}
       \hat{\Sigma}  & (\mu  - \hat{\mu} ) \\
       (\mu  - \hat{\mu} )^{T} & \beta 
   \end{array}\right) \succeq 0 \notag
\end{alignat}
Therefore $p_{\ell,t,k}^{j}$ satisfies the first condition in $\mathbb{D}_{\ell,t,k}^{j}$. For the second condition in $\mathbb{D}_{\ell,t,k}^{j}$, we verify the corresponding sufficient condition in (\ref{eB12a}) as follows.
\begin{small}
\begin{alignat}{2}
   & \Sigma_{\ell,t,k}^{j} = \sum\nolimits_{i=1}^{m_{\ell,t}} k_i^j \Sigma  \preceq \sum\nolimits_{i=1}^{m_{\ell,t}} k_i^j \varepsilon_{\ell,t}^{i} \hat{\Sigma}_{\ell,t}^{i} \label{eB15}\\
   \preceq &\max_{i=1...m_{\ell,t}} \varepsilon_{\ell,t}^{i} \sum\nolimits_{i=1}^{m_{\ell,t}} k_i^j \hat{\Sigma}_{\ell,t}^{i}  = (\varepsilon_{\ell,t,k}^{j} - \beta_{\ell,t,k}^{j}) \hat{\Sigma}_{\ell,t,k}^{j}    \notag
\end{alignat}    
\end{small}
\hspace{-3pt}Therefore $p_{\ell,t,k}^{j}$ satisfy the second condition in $\mathbb{D}_{\ell,t,k}^{j}$. This completes the proof that $p_{\ell,t,k}^{j} \in \mathbb{D}_{\ell,t,k}^{j}$.

Finally, the proof is complete. $\quad\square$ 
\end{pf}

\vspace{-0.2cm}
\section{Proof of Lemma \ref{lem3.3.1}}          \label{appendix C}
\vspace{-0.4cm}
\begin{pf}
In this proof, all the subscripts $\ell,t,k$ are omitted. Based on the ambiguity set (\ref{e3.2.4}), for all $ p \in \mathbb{D}$, we can obtain the following equation.
\begin{equation} \label{eC1}
\begin{split}
     p = \sum\nolimits_{j = 1}^{m} \gamma^j p^j = \sum\nolimits_{n=1}^{M} \widetilde{\gamma}^n \widetilde{p}^n
\end{split}
\end{equation}
where $p^j \in \mathbb{D}^j$, and $\widetilde{\gamma}^n, \widetilde{p}^n$ are as follows.
\begin{equation} \label{eC2}
\begin{split}
    &\widetilde{\gamma}^n = \sum\nolimits_{j \in \mathbb{M}^n}\gamma^j, \quad \widetilde{p}^n=\sum\nolimits_{j \in \mathbb{M}^n}(\gamma^j/\widetilde{\gamma}^n)p^j
\end{split}
\end{equation}
It remains to prove that $\bm{\widetilde{\gamma}} = (\widetilde{\gamma}^{1},...,\widetilde{\gamma}^{M}) \in \widetilde{\Delta}$ and $\widetilde{p}^n \in \widetilde{\mathbb{D}}^n$ in (\ref{e3.3.2.1}).

(1) Prove that $\bm{\widetilde{\gamma}} \in \widetilde{\Delta}$.\\
According to the definition of $\bm{\widetilde{\gamma}}$, we can obtain
\begin{small}
\begin{equation} \label{eC3}
\begin{split}
    &\sum\nolimits_{n=1}^{M} \widetilde{\gamma}^{n} = \sum\nolimits_{n=1}^{M} \sum\nolimits_{j \in \mathbb{M}^n}\gamma^j = \sum\nolimits_{j=1}^{m} \gamma^{j} = 1
\end{split}
\end{equation}    
\end{small}
Therefore, the first condition in $\widetilde{\Delta}$ is satisfied. For the second condition, we verify it as follows.
\begin{small}
\begin{alignat}{2}
    &\sum\nolimits_{n=1}^{M} \vert \widetilde{\gamma}^{n}-\hat{\widetilde{\gamma}}^{n} \vert = \sum\nolimits_{n=1}^{M} \left| \sum\nolimits_{j \in \mathbb{M}^n} (\gamma^{j} - \hat{\gamma}^{j}) \right|  \label{eC4}\\
    \leq & \sum\nolimits_{n=1}^{M} \sum\nolimits_{j \in \mathbb{M}^n} \vert \gamma^{j} - \hat{\gamma}^{j} \vert = \sum\nolimits_{j=1}^{m} \vert \gamma^{j} - \hat{\gamma}^{j} \vert  \leq \theta = \widetilde{\theta}     \notag
\end{alignat}    
\end{small}
\hspace{-3pt}where the first inequality follows from the triangle inequality. This completes the proof that $\bm{\widetilde{\gamma}} \in \widetilde{\Delta}$.

(2) Prove that $\widetilde{p}^n \in \widetilde{\mathbb{D}}^n$. \\
For simplicity of notation, we use $\mu^{j}, \widetilde{\mu}^{n}$ and $\Sigma^{j}, \widetilde{\Sigma}^{n}$ to represent the means and covariances of $p^{j}$ and $\widetilde{p}^{n}$. According the definition of $\widetilde{p}^{n}$ in (\ref{eC2}), we can obtain
\begin{alignat}{2}
    &\widetilde{\mu}^{n} = \sum\nolimits_{j \in \mathbb{M}_{n}} \frac{\gamma^{j}}{\widetilde{\gamma}^n} \mu^{j}   \label{eC5}\\
    &\widetilde{\Sigma}^{n} = \sum\nolimits_{j \in \mathbb{M}_{n}} \frac{\gamma^{j}}{\widetilde{\gamma}^n} \left\{ \Sigma^{j} + \left( \mu^{j} - \widetilde{\mu}^{n} \right) \left( \mu^{j} - \widetilde{\mu}^{n} \right)^T \right\}   \notag
\end{alignat}
We note that the two conditions in $\widetilde{\mathbb{D}}^n$ (\ref{e3.3.2.1}) also have the corresponding forms as in (\ref{eB11}) and (\ref{eB12}). Therefore, for the first condition in $\widetilde{\mathbb{D}}^n$, we verify its form corresponding to (\ref{eB11b}) below.
\begin{small}
\begin{alignat}{2}
    &\left(\begin{array}{cc}
        \hat{\widetilde{\Sigma}}^{n} & \left( \widetilde{\mu}^{n} - \hat{\widetilde{\mu}}^{n} \right) \\
        \left( \widetilde{\mu}^{n} - \hat{\widetilde{\mu}}^{n} \right)^T & \widetilde{\beta}^{n}
    \end{array}\right)          \label{eC6}\\
    \overset{(a)}{=}& \sum_{j \in \mathbb{M}_{n}} \left(\begin{array}{cc}
        \overline{\gamma}^{n} \hat{\Sigma}^{j} + \breve{\gamma}^{n} \left(\hat{\mu}^{j}\right) \left(\hat{\mu}^{j}\right)^T & * \\
        \left( \dfrac{\gamma^{j}}{\widetilde{\gamma}^n}\mu^{j} - \dfrac{\hat{\gamma}^{j}}{\hat{\widetilde{\gamma}}^n}\hat{\mu}^{j} \right)^T &  \overline{\gamma}^{n} \beta^{j} + \breve{\gamma}^{n}
    \end{array}\right) \notag
\end{alignat}    
\end{small}
\hspace{-3pt}where we only list the lower triangular elements of the matrix (a), since it is symmetric. The lower-left elements of matrix (a) can be further transformed as follows.
\begin{equation}    \label{eC7}
\begin{split}
    \frac{\gamma^{j}}{\widetilde{\gamma}^n}\mu^{j} - \frac{\hat{\gamma}^{j}}{\hat{\widetilde{\gamma}}^n}\hat{\mu}^{j}    =& \frac{\gamma^{j}}{\widetilde{\gamma}^n} \left( \mu^{j} - \hat{\mu}^{j} \right) + \left( \frac{\gamma^{j}}{\widetilde{\gamma}^n} - \frac{\hat{\gamma}^{j}}{\hat{\widetilde{\gamma}}^n} \right) \hat{\mu}^{j}
\end{split}
\end{equation}
Substituting (\ref{eC7}) into (\ref{eC6}), we can obtain
\begin{small}
\begin{alignat}{2}
    &\left(\begin{array}{cc}
        \hat{\widetilde{\Sigma}}^{n} & \left( \widetilde{\mu}^{n} - \hat{\widetilde{\mu}}^{n} \right) \\
        \left( \widetilde{\mu}^{n} - \hat{\widetilde{\mu}}^{n} \right)^T & \widetilde{\beta}^{n}
    \end{array}\right)          \label{eC8}\\
    =& \sum_{j \in \mathbb{M}_{n}} \left(\begin{array}{cc}
        \overline{\gamma}^{n} \hat{\Sigma}^{j} & * \\
        \dfrac{\gamma^{j}}{\widetilde{\gamma}^n} \left( \mu^{j} - \hat{\mu}^{j} \right)^T & \overline{\gamma}^{n} \beta^{j}
    \end{array}\right)+ \left(\begin{array}{cc}
        \breve{\gamma}^{n} \left(\hat{\mu}^{j}\right) \left(\hat{\mu}^{j}\right)^T & * \\
        \left( \dfrac{\gamma^{j}}{\widetilde{\gamma}^n} - \dfrac{\hat{\gamma}^{j}}{\hat{\widetilde{\gamma}}^n} \right) \left(\hat{\mu}^{j}\right)^T & \breve{\gamma}^{n}
    \end{array}\right) \notag\\
    \overset{(a)}{\succeq} & \sum_{j \in \mathbb{M}_{n}} \overline{\gamma}^{n} \left(\begin{array}{cc}
        \hat{\Sigma}^{j} & * \\
        \left( \mu^{j} - \hat{\mu}^{j} \right)^T & \beta^{j}
    \end{array}\right)+ \breve{\gamma}^{n} \left(\begin{array}{cc}
        \left(\hat{\mu}^{j}\right) \left(\hat{\mu}^{j}\right)^T & * \\
        \left(\hat{\mu}^{j}\right)^T & 1
    \end{array}\right) \overset{(b)}{\succeq} 0  \notag
\end{alignat}    
\end{small}
\hspace{-3pt}where the generalized inequality (a) follows from the definitions of $\overline{\gamma}^{n}$ and $\breve{\gamma}^{n}$ in (\ref{e3.3.2.2}), while the generalized inequality (b) arises from $p^{j} \in \mathbb{D}^{j}$. Therefore, $\widetilde{p}^{n}$ satisfy the first condition in $\widetilde{\mathbb{D}}^{n}$. For the second condition in $\widetilde{\mathbb{D}}_{\ell,t,k}^{n}$, we verify the corresponding sufficient condition in (\ref{eB12a}), which requires verifying the following inequality.
\begin{equation}    \label{eC9}
\begin{split}
    \widetilde{\Sigma}^{n} \preceq \hat{\widetilde{\Phi}}^{n} - \widetilde{\beta}^{n} \hat{\widetilde{\Sigma}}^{n}
\end{split}
\end{equation}
According to the definition of $\widetilde{\Sigma}^{n}$ in (\ref{eC5}), we first transform the $\left( \mu^{j} - \widetilde{\mu}^{n} \right) \left( \mu^{j} - \widetilde{\mu}^{n} \right)^T$ part as follows.
\begin{small}
\begin{alignat}{2}
    & \left( \mu^{j} - \widetilde{\mu}^{n} \right) \left( \mu^{j} - \widetilde{\mu}^{n} \right)^T   \label{eC10}\\
    =& \left\{ \left( \mu^{j} - \hat{\mu}^{j} \right) - \left( \widetilde{\mu}^{n} - \hat{\widetilde{\mu}}^{n} \right) + \left( \hat{\mu}^{j} - \hat{\widetilde{\mu}}^{n} \right) \right\} \notag\\
    & \left\{ \left( \mu^{j} - \hat{\mu}^{j} \right) - \left( \widetilde{\mu}^{n} - \hat{\widetilde{\mu}}^{n} \right) + \left( \hat{\mu}^{j} - \hat{\widetilde{\mu}}^{n} \right) \right\}^T \notag\\
    \overset{(a)}{\preceq}& 3\left( \mu^{j} - \hat{\mu}^{j} \right)\left( \mu^{j} - \hat{\mu}^{j} \right)^T + 3\left( \widetilde{\mu}^{n} - \hat{\widetilde{\mu}}^{n} \right)\left( \widetilde{\mu}^{n} - \hat{\widetilde{\mu}}^{n} \right)^T + \notag\\
    &3\left( \hat{\mu}^{j} - \hat{\widetilde{\mu}}^{n} \right)\left( \hat{\mu}^{j} - \hat{\widetilde{\mu}}^{n} \right)^T \notag\\
    \overset{(b)}{\preceq}& 3\beta^{j}\hat{\Sigma}^{j} + 3\widetilde{\beta}^{n} \hat{\widetilde{\Sigma}}^{n} + 3\left( \hat{\mu}^{j} - \hat{\widetilde{\mu}}^{n} \right)\left( \hat{\mu}^{j} - \hat{\widetilde{\mu}}^{n} \right)^T     \notag
\end{alignat}    
\end{small}
\hspace{-3pt}where the generalized inequality (b) follows from the inequality (\ref{eC8}) corresponding to the from in (\ref{eB11c}) and the fact $p^{j} \in \mathbb{D}^{j}$. Substituting (\ref{eC10}) into the definition of $\widetilde{\Sigma}^{n}$ in (\ref{eC5}), we can obtain
\begin{small}
\begin{alignat}{2}
    & \widetilde{\Sigma}^{n}        \label{eC11}\\
    \overset{(a)}{\preceq}&  \sum_{j \in \mathbb{M}_{n}} \frac{\gamma^{j}}{\widetilde{\gamma}^n} \left\{ \Sigma^{j} + 3\beta^{j}\hat{\Sigma}^{j} + 3( \hat{\mu}^{j} - \hat{\widetilde{\mu}}^{n} )( \hat{\mu}^{j} - \hat{\widetilde{\mu}}^{n} )^T \right\} + 3\widetilde{\beta}^{n} \hat{\widetilde{\Sigma}}^{n} \notag\\
    \overset{(b)}{\preceq}& \sum_{j \in \mathbb{M}_{n}} \frac{\gamma^{j}}{\widetilde{\gamma}^n} \left\{ ( \varepsilon^{j} + 3\beta^{j}) \hat{\Sigma}^{j} + 3( \hat{\mu}^{j} - \hat{\widetilde{\mu}}^{n} )( \hat{\mu}^{j} - \hat{\widetilde{\mu}}^{n} )^T \right\} + 3\widetilde{\beta}^{n} \hat{\widetilde{\Sigma}}^{n} \notag\\
    \overset{(c)}{\preceq}& \sum_{j \in \mathbb{M}_{n}} \bar{\gamma}_{\varphi}^{n} \left\{ ( \varepsilon^{j} + 3\beta^{j}) \hat{\Sigma}^{j} + 3( \hat{\mu}^{j} - \hat{\widetilde{\mu}}^{n} )( \hat{\mu}^{j} - \hat{\widetilde{\mu}}^{n} )^T \right\} + 3\widetilde{\beta}^{n} \hat{\widetilde{\Sigma}}^{n} \notag\\
    =& \hat{\widetilde{\Phi}}^{n} - \widetilde{\beta}^{n} \hat{\widetilde{\Sigma}}^{n}  \notag
\end{alignat}
\end{small}
\hspace{-4pt}where the generalized inequality (b) arises from $p^{j} \in \mathbb{D}_{\ell,t,k}^{j}$, and the generalized inequality (c) follows from the definitions of $\bar{\gamma}^{n}$. Therefore, $\widetilde{p}^{n}$ satisfies the second condition in $\widetilde{\mathbb{D}}^{n}$. This completes the proof that $\widetilde{p}^{n} \in \widetilde{\mathbb{D}}^{n}$.

Finally, the proof is complete. $\quad\square$ 
\end{pf}

\vspace{-0.3cm}
\section{Proof of Theorem \ref{thm4.2.1}}          \label{appendix D}
\vspace{-0.3cm}
\begin{pf}
By the definition of CVaR and $\mathcal{L}_{\mathcal{H}}(\varphi\oplus\mathcal{O})$, the DR-CVaR (\ref{e2.3.2.5}) can be rewritten below.
\begin{small}
\begin{align}
&\mathop{\sup}\nolimits_{\upsilon \in \widetilde{\mathbb{D}}} {\rm CVaR}_{\alpha^{u}}^{\upsilon} ( \mathcal{L}_{\mathcal{H}} (\varphi \oplus \mathcal{O}))    \label{eD1}\\
=&\sup\nolimits_{\upsilon \in \widetilde{\mathbb{D}}} \inf\nolimits_{z \in \mathbb{R}} \left( z + \tfrac{1}{1-\alpha^{u}} \mathbb{E}^{\upsilon} \left[ (-h^T \varphi -g +\mathcal{S}_{\mathcal{O}}(-h) -z )^+  \right]  \right)   \notag\\
=&\inf\nolimits_{z \in \mathbb{R}} \left( z +\tfrac{1}{1-\alpha^{u}} \sup\nolimits_{\upsilon \in \widetilde{\mathbb{D}}} \mathbb{E}^{\upsilon} \left[ (-h^T \varphi -g +\mathcal{S}_{\mathcal{O}}(-h) -z)^+ \right]  \right)      \notag
\end{align}        
\end{small}
\hspace{-3pt}where the second equality follows from the minimax theorem \cite{shapiro2002minimax}. Based on the definition of the ambiguity set (\ref{e3.3.2.2}), the worst-case expectation problem $\sup_{\upsilon \in  \widetilde{\mathbb{D}}} \mathbb{E}^{\upsilon} \left[ (-h^T \varphi -g +\mathcal{S}_{\mathcal{O}}(-h) -z)^+ \right]$ can be rewritten as follows.
\begin{align}
     & \sup\nolimits_{\upsilon \in \widetilde{\mathbb{D}}} \mathbb{E}^{\upsilon} \left[ (-h^T \varphi -g +\mathcal{S}_{\mathcal{O}}(-h) -z)^+ \right] \label{eD2}\\
    =& \sup_{\bm{\gamma} \in \widetilde{\Delta}} \sum\nolimits_{n=1}^{M} \gamma^{n} \sup_{p^{n} \in \widetilde{\mathbb{D}}^{n}} \mathbb{E}^{p^{n}} \left[ (-h^T \varphi -g +\mathcal{S}_{\mathcal{O}}(-h) -z)^+ \right]  \notag
\end{align}
The worst-case expectation problem (\ref{eD2}) includes two layers of robustness. We first focus on the inner maximization problem, which can be reformulated below.
\begin{small}
\begin{align}
     \sup_{p^{n}\in\mathcal{M}(\mathbb{W})} & \int_{\mathbb{W}} (-h^T \varphi -g +\mathcal{S}_{\mathcal{O}}(-h) -z)^+ p^{n}(d\varphi)   \label{eD3}\\
     s.t. & \int_{\mathbb{W}} -\left(\begin{array}{cc}
         \hat{\widetilde{\Sigma}}^{n} & \varphi - \hat{\widetilde{\mu}}^{n}     \notag\\
         (\varphi - \hat{\widetilde{\mu}}^{n})^T & \hat{\widetilde{\beta}}^{n}
     \end{array}  \right) p^{n}(d\varphi) \preceq 0     \notag\\
     & \int_{\mathbb{W}} \left(\varphi-\hat{\widetilde{\mu}}^{n} \right) \left( \varphi-\hat{\widetilde{\mu}}^{n} \right)^{T} p^{n}(d\varphi) \preceq\hat{\widetilde{\Phi}}^{n}     \notag
\end{align}    
\end{small}
\hspace{-3pt}where $\hat{\widetilde{\mu}}^n$, $\hat{\widetilde{\Sigma}}^n$, $\hat{\widetilde{\beta}}^n$, and $\hat{\widetilde{\Phi}}^n$ are the elements in the ambiguity set (\ref{e3.3.2.1}). We can obtain the dual problem of (\ref{eD3}) as follows \cite{boyd2004convex}.
\begin{small}
\begin{align}
    \inf_{\substack{\Lambda^{n}, \xi^{n},\\ \tau^{n}, \Omega^{n}}}& \left\{\begin{array}{l}
     \tau^{n}\hat{\widetilde{\beta}}^{n} - 2(\hat{\widetilde{\mu}}^{n})^T\xi^{n} - (\hat{\widetilde{\mu}}^{n})^T \Omega^{n} \hat{\widetilde{\mu}}^{n} + \hat{\widetilde{\Sigma}}^{n} \bullet \Lambda^{n} +    \\
     \hat{\widetilde{\Phi}}^{n} \bullet \Omega^{n}+ \sup\limits_{\varphi \in \mathbb{W}} \left\{\begin{array}{l}
          (-h^T \varphi -g +\mathcal{S}_{\mathcal{O}}(-h) -z)^+   \\
          + 2\varphi^T(\xi^{n} + \Omega^{n}\mu^{n}) - \varphi^T\Omega^{n}\varphi 
     \end{array}\right\}
    \end{array}\right.     \notag\\
    s.t.& \left(\begin{array}{cc}
        \Lambda^{n} & \xi^{n}    \\
        (\xi^{n})^T & \tau^{n}
    \end{array}\right) \succeq 0 \qquad \Omega^{n} \succeq 0  \label{eD5}
\end{align}    
\end{small}
\hspace{-3pt}where $\{  \Lambda^{n}, \xi^{n},\tau^{n} \}$ and $\Omega^{n}$ are the dual variables corresponding to the two constraints in (\ref{eD3}). For matrices $A$ and $B$, $A \bullet B$ represents their Frobenius inner product. We denote the optimal value of the above optimization problem by $\phi^{n}$. The outer maximization problem of (\ref{eD2}) can be reformulated below.
\begin{equation}    \label{eD6}
\begin{split}
    \sup\nolimits_{\bm{\gamma} \in \widetilde{\Delta}} \sum\nolimits_{n=1}^{M} \gamma^{n} \cdot \phi^{n}
\end{split}
\end{equation}
Applying the strong duality of linear programming, we have the following.
\begin{subequations}    \label{eD7}
\begin{alignat}{2}
    \inf_{\vartheta, \zeta, r_1^n, r_2^n}& \hat{\widetilde{\theta}} \zeta + \sum\nolimits_{n=1}^{M} \hat{\widetilde{\gamma}}^{n} (r_1^n-r_2^n) + \vartheta \label{eD7a}\\
    s.t. \quad & r_1^n - r_2^n + \vartheta \geq \phi^{n} \label{eD7b}\\
    & \zeta - r_1^n - r_2^n = 0     \label{eD7c}\\
    & \zeta \geq 0, \quad r_1^n \geq 0, \quad r_2^n \geq 0  \label{eD7d}
\end{alignat}
\end{subequations}
By substituting (\ref{eD5}) into (\ref{eD7b}) and introducing an auxiliary variable $s^{n}$, we obtain
\begin{subequations}    \label{eD8}
\begin{alignat}{2}
    & r_1^n - r_2^n + \vartheta \geq \tau^{n}\hat{\widetilde{\beta}}^{n} - 2(\hat{\widetilde{\mu}}^{n})^T\xi^{n} - (\hat{\widetilde{\mu}}^{n})^T \Omega^{n} \hat{\widetilde{\mu}}^{n}     \label{eD8a}\\
    & \quad\quad\quad\quad\quad\quad + \hat{\widetilde{\Sigma}}^{n} \bullet \Lambda^{n} + \hat{\widetilde{\Phi}}^{n} \bullet \Omega^{n} + s^{n} \notag\\
    & g -\mathcal{S}_{\mathcal{O}}(-h) +z + h^T \varphi +      \notag\\
    & s^{n} -2\varphi^T(\xi^{n} + \Omega^{n}\hat{\widetilde{\mu}}^{n}) +\varphi^T\Omega^{n}\varphi \geq 0,  \ \forall \varphi \in \mathbb{W}     \label{eD8b}\\
    & s^{n} -2\varphi^T(\xi^{n} + \Omega^{n}\hat{\widetilde{\mu}}^{n}) +\varphi^T\Omega^{n}\varphi \geq 0,  \ \forall \varphi \in \mathbb{W}     \label{eD8c}
\end{alignat}
\end{subequations}
Constraints (\ref{eD8b}) and (\ref{eD8c}) are semi-infinite constraints, so some further reformulations are needed. Constraints (\ref{eD8b}) can be rewritten as follows.
\begin{equation} \label{eD9}
\begin{split}
\min_{\varphi\in\mathbb{W}}\left\{\begin{array}{l}
    g -\mathcal{S}_{\mathcal{O}}(-h) +z + s^{n} + \\
    \varphi^T(h - 2\xi^{n} - 2\Omega^{n}\hat{\widetilde{\mu}}^{n}) + \varphi^T\Omega^{n}\varphi
\end{array} \right\} \geq 0
\end{split}
\end{equation}
Based on the duality of convex quadratic program and Schur complements theorem \cite{boyd2004convex}, (\ref{eD9}) is reformulated as the following LMI.
\begin{small}
\begin{equation} \label{eD12}
\begin{split}
\hspace{-2pt}  \left(\begin{array}{cc}
  \Omega^{n}   & * \\
  \dfrac{1}{2}\left(\begin{array}{c}
                           E^T\eta^{n} - 2\xi^{n}  \\
                           +h -2\Omega^{n}\hat{\widetilde{\mu}}^{n}
                      \end{array} \right)^T   & \left.\begin{array}{l}
                                                          g-\mathcal{S}_{\mathcal{O}}(-h)+  \\
                                                          z+s^{n}-f^T\eta^n 
                                                     \end{array}\right.
\end{array}\right) \succeq 0
\end{split}
\end{equation}
\end{small}
\hspace{-2pt}where $\eta^{n} \geq 0$ is the dual variables corresponding to the constraint $E\varphi \leq f$. Similarly, constraint (\ref{eD8c}) can be reformulated as LMI constraints below.
\begin{small}
\begin{equation} \label{eD13}
\begin{split}
\left(\begin{array}{cc}
  \Omega^{n}   & * \\
  \tfrac{1}{2}\left( E^T\lambda^{n} - 2\xi^{n} -2\Omega^{n}\hat{\widetilde{\mu}}^{n} \right)^T   & s^{n}-f^T\lambda^n 
\end{array}\right) \succeq 0
\end{split}
\end{equation}
\end{small}
\hspace{-3pt}where $\lambda^{n}\geq0$ is the dual variables corresponding to the constraint $E\varphi\leq f$. Substituting (\ref{eD7}), (\ref{eD8}), (\ref{eD12}), and (\ref{eD13}) into (\ref{eD1}), the theorem is proved. $\quad\square$ 
\end{pf}

\vspace{-0.2cm}
\section{Proof of Theorem \ref{thm4.4.1}}          \label{appendix E}
\vspace{-0.2cm}
\begin{pf}
Due to the spatial allocation protocol, if (\ref{e2.4.2.1}) is feasible at $t$, the collision between two robots will never occur. Next, we analyze the collision avoidance between robots and obstacles. According to Theorem \ref{thm3.2.1}, Lemma \ref{lem3.2.1} and \ref{lem3.3.1}, we know that
\begin{small}
\begin{equation}        \label{eE1}
\begin{split}
    \mathbb{P}\{ p_{\ell,t,k}^{*} \in \widetilde{\mathbb{D}}_{\ell,t,k}, \forall k=1,...,K-1 \} \geq \alpha_{\ell}
\end{split}
\end{equation}
\end{small}
\hspace{-3pt}At each control time $t$, for the obstacle $\ell \in \mathbb{O}_{i,t}$, the probability of collision avoidance at time $t+1$ is given by (\ref{e2.3.2.5}) as $\alpha^{u}$. As a result, the control trajectory $y_{i}(t+1) = y_{i}(1|t)$ is collision avoidance with obstacle $\ell$ with at least $\alpha_{\ell} \alpha^{u}$. For the obstacle $\ell \notin \mathbb{O}_{i,t}$, according to Section \ref{section4.3}, the formulation of the avoidance constraints with obstacle $\ell$ is either based on the ambiguity set from a past time step or based on the support set. Therefore, the trajectory also avoids the collision with obstacle $\ell$ with at last a probability of $\alpha_{\ell} \alpha^{u}$. In summary, the total collision avoidance probability of $y_{i}(t+1) = y_{i}(1|t)$ is at least $\prod_{\ell=1}^{L} (\alpha_{\ell} \alpha)^{L}$. Therefore, if (\ref{e2.4.2.1}) is feasible at each time $t = 0,...,T-1$, the trajectory $y_{i}(1),...,y_{i}(T)$ is collision-free with at least $\prod_{\ell=1}^{L}(\alpha_{\ell} \alpha^{u})^{T}$. This completes the proof. $\quad\square$ 
\end{pf}

\vspace{-0.2cm}
\section{Proof of Theorem \ref{thm4.4.2}}          \label{appendix F}
\vspace{-0.2cm}
\begin{pf}
We assume that there exists an optimal control input sequences $\bm{u_{i}^{*}(t)}$ and corresponding stage and output trajectories $\bm{x_{i}^{*}(t)}$, $\bm{y_{i}^{*}(t)}$ at time $t$ as follows.
\begin{equation}    \label{eF1}
\begin{split}
    &\bm{u_{i}^{*}(t)} = ( u_{i}^{*}(0|t),...,u_{i}^{*}(K-1|t) ) \\
    &\bm{x_{i}^{*}(t)} = ( x_{i}^{*}(0|t),...,x_{i}^{*}(K|t) ) \\
    &\bm{y_{i}^{*}(t)} = ( y_{i}^{*}(0|t),...,y_{i}^{*}(K|t) ) 
\end{split}
\end{equation}
We construct the candidate control sequence and corresponding stage and output trajectories at time $t+1$ using the terminal control law $\kappa_{i}^{f}(\cdot)$ as follows.
\begin{small}
\begin{alignat}{2}
    \bm{u_{i}(t+1)} &= ( u_{i}(0|t+1),...,u_{i}(K-1|t+1) )      \label{eF2}\\
    &= ( u_{i}^{*}(1|t),...,u_{i}^{*}(K-1|t), \kappa_{i}^{f}(x_{i}^{*}(K|t)) )     \notag\\ 
    &= ( u_{i}^{*}(1|t),...,u_{i}^{*}(K-1|t), u_{i}^{f}(K|t) )      \notag\\ 
    \bm{x_{i}(t+1)} &= ( x_{i}(0|t+1),...,x_{i}(K|t+1) ) \notag\\
    &= ( x_{i}^{*}(1|t),...,x_{i}^{*}(K|t), f(x_{i}^{*}(K|t), u_{i}^{f}(K|t)) ) \notag\\
    &= ( x_{i}^{*}(1|t),...,x_{i}^{*}(K|t), x_{i}^{f}(K+1|t) ) \notag\\
    \bm{y_{i}(t+1)} &= ( y_{i}(0|t+1),...,y_{i}(K|t+1) ) \notag\\
    &= ( y_{i}^{*}(1|t),...,y_{i}^{*}(K|t), h(x_{i}^{f}(K+1|t)) )  \notag\\
    &= ( y_{i}^{*}(1|t),...,y_{i}^{*}(K|t), y_{i}^{f}(K+1|t) ) \notag
\end{alignat}    
\end{small}
\hspace{-3pt}Next, we will verify the feasibility of the above candidate solution w.r.t. each constraint in (\ref{e2.4.2.1}) at time $t+1$.

First, for constraints~(\ref{e2.4.2.1b})-(\ref{e2.4.2.1f}), by the construction of the candidate solution~(\ref{eF2}) and the definition of the terminal set $\mathcal{X}_{i}^{f}$,  the feasibility of candidate solution (\ref{eF2}) is straightforward.

Second, for the collision avoidance affine constraints~(\ref{e2.4.2.1g}) and~(\ref{e2.4.2.1h}), the safe update rules for separating hyperplanes proposed in Section~\ref{section4.3} ensure that, at each time step, the adopted separating hyperplanes are compatible with the current committed trajectory (i.e., the candidate output). As a result, these affine constraints remain feasible for the candidate solution~(\ref{eF2}).

Finally, for the terminal constraint~(\ref{e2.4.2.1i}), by the construction of the terminal set $\mathcal{X}_{i}^{f}$ in Definition~\ref{defn4.3.1}, the terminal state corresponding to the candidate solution satisfies $\hat{x}_{i}(K|t+1) \in \mathcal{X}_{i}^{f}$. Therefore, the terminal constraint is also satisfied.

As a result, the candidate solution (\ref{eF2}) is feasible at time $t+1$. And we can obtain the following expression.
\begin{small}
\begin{alignat}{2}
    &J_{i}(\bm{u_{i}(t+1)})    \label{eF6}\\
    =& \sum\nolimits_{k=0}^{K-1} l_{i}(x_{i}(k|t+1), u_{i}(k|t+1))    \notag\\
     &+ q_{i}(x_{i}(K|t+1))     \notag\\
    =& J_{i}^{*}(\bm{u_{i}^{*}(t)}) - l_{i}(x_{i}^{*}(0|t), u_{i}^{*}(0|t)) - q_{i}(x_{i}^{f}(K|t)) \notag\\
     &+ l_{i}(x_{i}^{*}(K|t), u_{i}^{f}(K|t)) + q_{i}(x_{i}^{f}(K+1|t))   \notag\\
    \leq& J_{i}^{*}(\bm{u_{i}^{*}(t)}) - l_{i}(x_{i}^{*}(0|t), u_{i}^{*}(0|t))    \notag
\end{alignat}    
\end{small}
\hspace{-3pt}where the inequality follows from Assumption \ref{assum4.4.2}. Using the optimality and the control law $u_{i}(t) = u_{i}^{*}(0|t)$, we have the following inequality.
\begin{alignat}{2}
    &J_{i}^{*}(\bm{u_{i}^{*}(t+1)}) - J_{i}^{*}(\bm{u_{i}^{*}(t)}) \leq - l_{i}(x_{i}(t), u_{i}(t))    \notag
\end{alignat}
We sum both sides of the above inequality from $t=0$ to $t=T-1$ and divide by $T \geq 0$ to give
\begin{alignat}{2}
    \frac{1}{T} \sum_{t=0}^{T} l_{i}(x_{i}(t), u_{i}(t)) \leq \frac{J_{i}^{*}(\bm{u_{i}^{*}(T)}, T) - J_{i}^{*}(\bm{u_{i}^{*}(t)})}{T}    \notag
\end{alignat}
Considering the Assumption \ref{assum4.4.1} and taking $T \rightarrow \infty$, we can obtain
\begin{alignat}{2}
    \lim_{T \rightarrow \infty} \frac{1}{T} \sum\nolimits_{t=0}^{T} g(\Vert x_{i}(t)-r_{i}^{x}(t) \Vert) = 0   \label{eF8}
\end{alignat}
Since $g$ is a $\mathcal{K}_{\infty}$-function, the theorem is proved. $\quad\square$ 
\end{pf}

\bibliographystyle{cas-model2-names}        
\bibliography{root}           

\end{document}